%% file: 1_AIMS_0_main.tex
\documentclass{ipiEA} 
\usepackage{amsmath}
\usepackage{paralist}
\usepackage[misc]{ifsym}
\usepackage{epsfig} 
\usepackage{epstopdf} 
\usepackage[colorlinks=true]{hyperref}
\hypersetup{urlcolor=blue, citecolor=red}
\allowdisplaybreaks

\def\currentvolume{X}
 \def\currentissue{X}
  \def\currentyear{200X}
   \def\currentmonth{XX}
    \def\ppages{X-XX}
     \def\DOI{10.3934/ipi.2023040}

\theoremstyle{definition}

\usepackage{lipsum}
\usepackage{amsmath,amssymb,amsfonts}
\usepackage{bbm}

\usepackage{mathrsfs}

\usepackage{graphicx}
\usepackage{float}
\usepackage{caption}
\usepackage[labelformat=simple]{subcaption}

\usepackage[export]{adjustbox}

\usepackage{cprotect}

\usepackage[ruled, vlined, linesnumbered]{algorithm2e}
\SetKwRepeat{Do}{do}{while}
\SetKwComment{Comment}{}{}

\renewcommand{\thesection}{\arabic{section}}

\newcommand{\R}{\mathbb{R}}

\newcommand{\N}{\mathbb{N}}
\newcommand{\id}{\mathrm{d}}
\newcommand{\grad}{\nabla} 
\newcommand{\laplace}{\bigtriangleup} 

\usepackage{tikz}
\usetikzlibrary{calc, positioning, shapes, backgrounds, fit, arrows,angles,quotes}
\usepackage{pgf-spectra}
\usepackage{siunitx}
\usepackage{contour}

\usepackage{enumitem}
\setlist[enumerate]{leftmargin=.5in}
\setlist[itemize]{leftmargin=.5in}

\usepackage{cleveref}

\usepackage{setspace}
\usepackage{multicol}
\usepackage{tabularx}
\usepackage{float}
\usepackage{makecell}
\newcommand{\tabincell}[2]{\begin{tabular}{@{}#1@{}}#2\end{tabular}}

\usepackage{booktabs}

\DeclareMathAlphabet{\mathpzc}{OT1}{pzc}{m}{it}

\title[Surface reconstruction from few slices]
{
Super-resolution surface reconstruction\\
from few low-resolution slices
} 

\author[Yiyao Zhang, Ke Chen and Shang-Hua Yang]{}

\subjclass{Primary:
49Q20,
65K10,
65D18,
94A08;
Secondary:
68U10.
}
\keywords{Surface reconstruction,
  variational model,
  perimeter energy,
  willmore energy,
  Euler-Elastica energy,
  phase-field approximation,
  $\Gamma$-convergence,
  alternating direction method of multipliers,
  gaussian curvature,
  mean curvature,
  discrete geometry. }

\thanks{The first author is partially supported by the UoL-NTHU Dual PhD Programme. }

\thanks{$^*$Corresponding author: Ke Chen (\href{mailto:K.Chen@strath.ac.uk}
      {K.Chen@strath.ac.uk})
 [\href{https://www.liv.ac.uk/~cmchenke}{www.liv.ac.uk/\url{~}cmchenke}]. }

\begin{document}
\maketitle

\centerline{
\scshape
Yiyao Zhang$^{{\href{mailto:Yiyao.Zhang@liverpool.ac.uk; yiyaozhang@gapp.nthu.edu.tw}{\textrm{\Letter}}}1,2}$,
Ke Chen$^{{\href{mailto:K.Chen@strath.ac.uk}{\textrm{\Letter}}}*3,4}$ and
Shang-Hua Yang$^{{\href{mailto:shanghua@ee.nthu.edu.tw}{\textrm{\Letter}}}2,5}$
}

\medskip

{\footnotesize
 \centerline{$^1$
            Centre for Mathematical Imaging Techniques
            and Department of Mathematical Sciences
            }
 \centerline{
            University of Liverpool, Liverpool, UK
            }
} 

\medskip

{\footnotesize
 \centerline{$^2$
            Institute of Electronics Engineering, 
            National Tsing Hua University, Taiwan
            }
} 

\medskip

{\footnotesize
 \centerline{$^3$
            Department of Mathematics and Statistics, 
            University of Strathclyde, Glasgow, UK 
            }
} 

\medskip

{\footnotesize
 \centerline{$^4$
 	    Centre for Mathematical Imaging Techniques, 
	    University of Liverpool, Liverpool, UK
            }
} 

\medskip

{\footnotesize
 \centerline{$^5$
 	    Department of Electrical Engineering, 
	    National Tsing Hua University, Taiwan
            }
} 

\bigskip

 \centerline{(Communicated by Weihong Guo)}


\begin{abstract}
In many imaging applications where segmented features (e.g. blood vessels) are further used for other numerical simulations (e.g. finite element analysis), the obtained surfaces do not have fine resolutions suitable for the task.
Increasing the resolution of such surfaces becomes crucial.
This paper proposes a new variational model for solving this problem, based on an Euler-Elastica-based regulariser.
Further, we propose and implement two numerical algorithms for solving the model, a projected gradient descent method and the alternating direction method of multipliers.
Numerical experiments using real-life examples (including two from outputs of another variational model) have been illustrated for effectiveness.
The advantages of the new model are shown through quantitative comparisons by the standard deviation of Gaussian curvatures and mean curvatures from the viewpoint of discrete geometry.
\end{abstract}


\section{Introduction}
\label{sec:Intro}

\begin{figure}[htbp] 
    \centering
    \begin{subfigure}[t]{0.4\textwidth}
    \includegraphics[width=\textwidth]{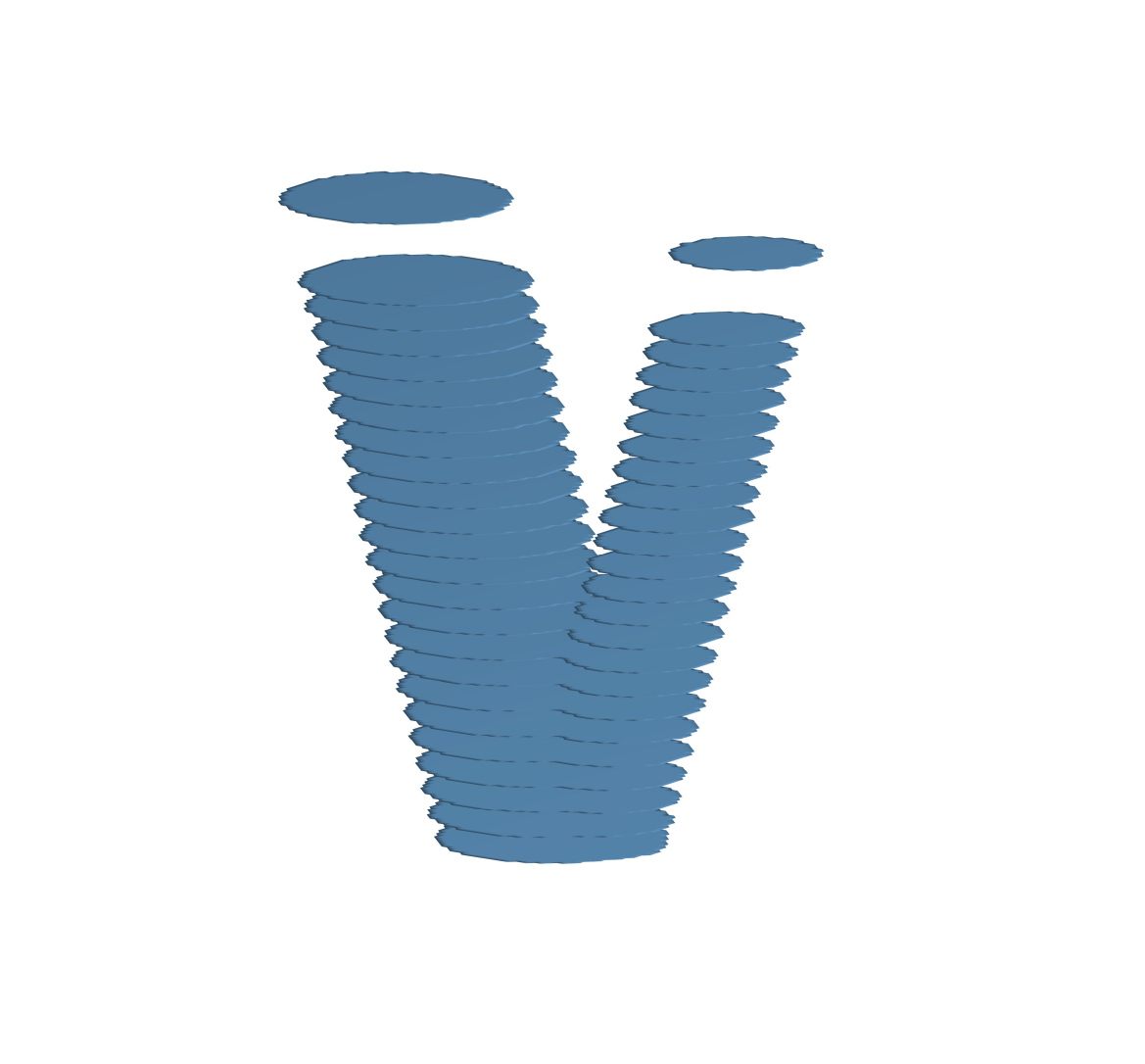}
    \label{p20220215_BC_N128S24_NoTITLE_GCMClim_PWEE_Slices}
    \end{subfigure}
\hspace{0.1mm}
    \begin{subfigure}[t]{0.4\textwidth}
    \includegraphics[width=\textwidth]{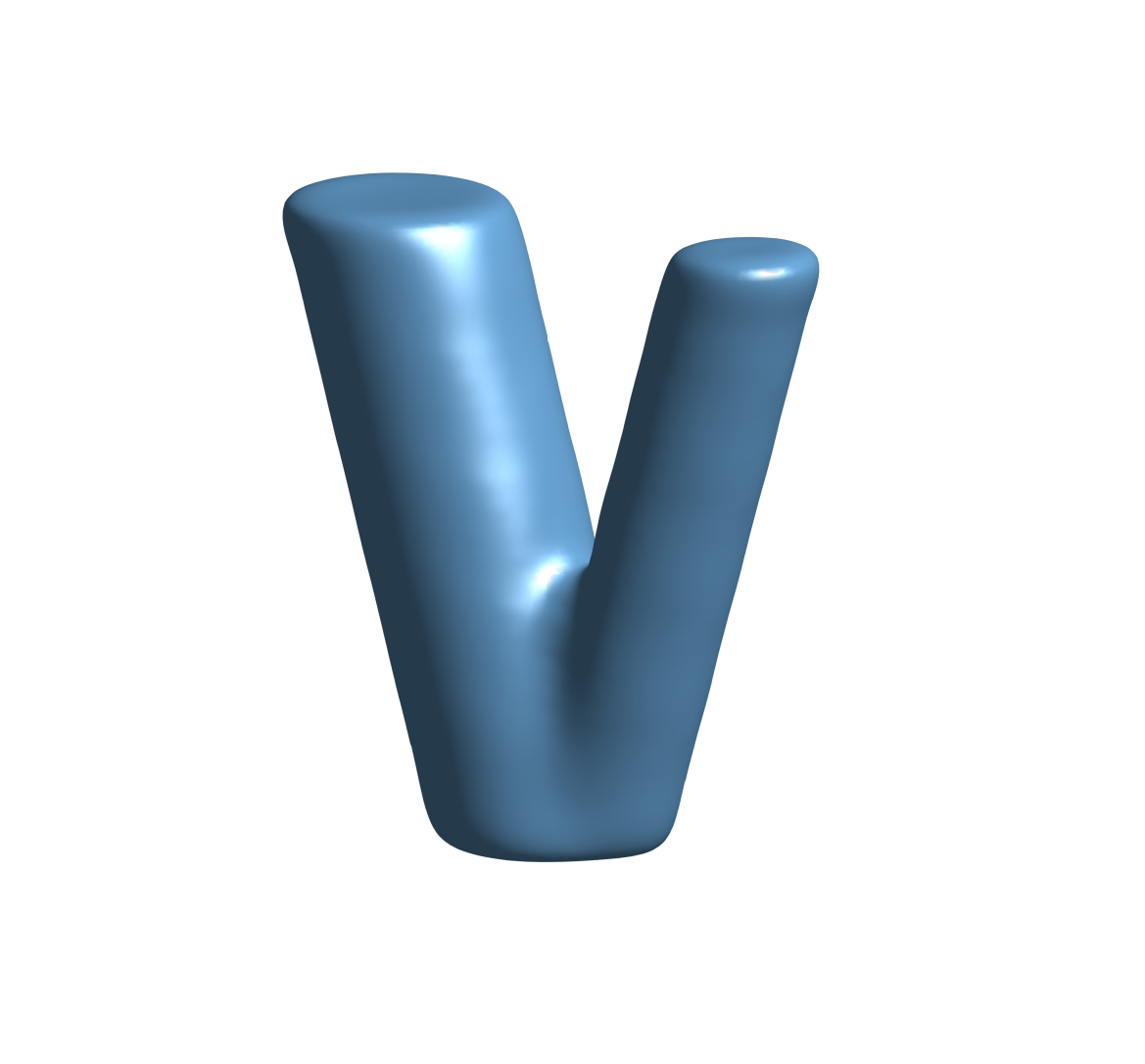}
    \label{p20220215_BC_N128S24_NoTITLE_GCMClim_PWEE_EE}
    \end{subfigure}
\cprotect \caption{Three-dimensional reconstruction of Branching Cylinders by the new proposed Euler-Elastica-based model from low-resolution inputs $N = 128$ with $24$ deliberately uneven slices. }
\label{p20220215_BC_Slices_EE}
\end{figure}

In this paper, we propose a new variational model for getting a faster and smoother three-dimensional (3D) surface reconstruction in high resolution from the collection of a few low-resolution images (input of cross-sections or slices), where the gaps of given slices are often large and uneven.
There are various practical or operational reasons why only low-resolution data (often with poor quality by way of artefacts or noise) are available.
For example, imaging equipment may have limitations,  such as short scanning times or low radiation doses to minimise harm or damage to patients.
A more subtle reason is that the region of interest (RoI) may be very small in a large and high-resolution (HR) image.
Hence, it is necessary to develop appropriate and suitable mathematical models to reconstruct HR surfaces as far as the RoI is concerned.
Once a HR datum (surface) is achieved, further measures or even new simulations based on the new geometry can be conducted, such as finite element analysis of blood flows in clinical imaging.
Though the surface reconstruction problem arises as a problem in discrete geometry (see \Cref{p20220215_BC_Slices_EE}), we shall develop a variational model for solving it.


The problem under study might be classified as an inpainting problem in which gaps are filled in some manner \cite{Carola15}.
Mathematical methods using variational models based on partial-differential-equation-based (PDE-based) have been widely exploited in image processing since the 1990s, as these approaches can accurately simulate our real, visible and physical world as evident from diverse application areas such as medicine, economics, and computer vision~\cite{Mathematical_problems_in_image_processing, Image_processing_and_analysis}.
Beyond imaging processing, the topic of 3D surface inpainting is much less studied in the variational framework.

In discrete geometry, many works exist for surface  reconstruction.
For retrieving the 3D shape, mathematical methods can be roughly divided into two categories depending on the surface representations:
(i) the explicit reconstruction
(e.g.
Delaunay  triangulation~\cite{Delaunay_Triangulation_Based_Surface_Reconstruction:_Ideas_and_Algorithms},
and Voronoi diagram~\cite{A_new_Voronoi-based_surface_reconstruction_algorithm, Surface_reconstruction_by_Voronoi_filtering, Voronoi-based_variational_reconstruction_of_unoriented_point_sets}),
and
(ii) the implicit reconstruction
(e.g.
radial basis functions~\cite{Reconstruction_and_representation_of_3D_objects_with_radial_basis_functions},
Poisson reconstruction~\cite{Poisson_surface_reconstruction, Screened_poisson_surface_reconstruction},
and level set method~\cite{01_Volume_Reconstruction_from_Slices, Fast_Surface_Reconstruction_Using_the_Level_Set_Method, A_Fast_Marching_Level_Set_Method_for_Monotonically_Advancing_Fronts}). 
In a nutshell, surfaces by category (i) are typically piecewise linear and these methods are relatively easy to implement. 
Although this kind of representation is suitable for graphics purposes, it is tough to cope with the non-uniform, potentially noisy, or incomplete data. 
For instance, when the number of vertices and faces in the surface is not sufficiently large, it can also have difficulties in tracking topological changes and large deformations. 
Category (ii) methods can produce reconstructions through forming a physical model based on differentiable distance functions over an implicit surface and are capable of coping better with the cases of non-uniform, potentially noisy, or incomplete data. 

\begin{figure}[H]
    \centering
    \vspace{-.5cm}
    \begin{subfigure}[t]{\textwidth}
    \begin{tikzpicture}
    \node
    {\includegraphics[width=\textwidth]
    {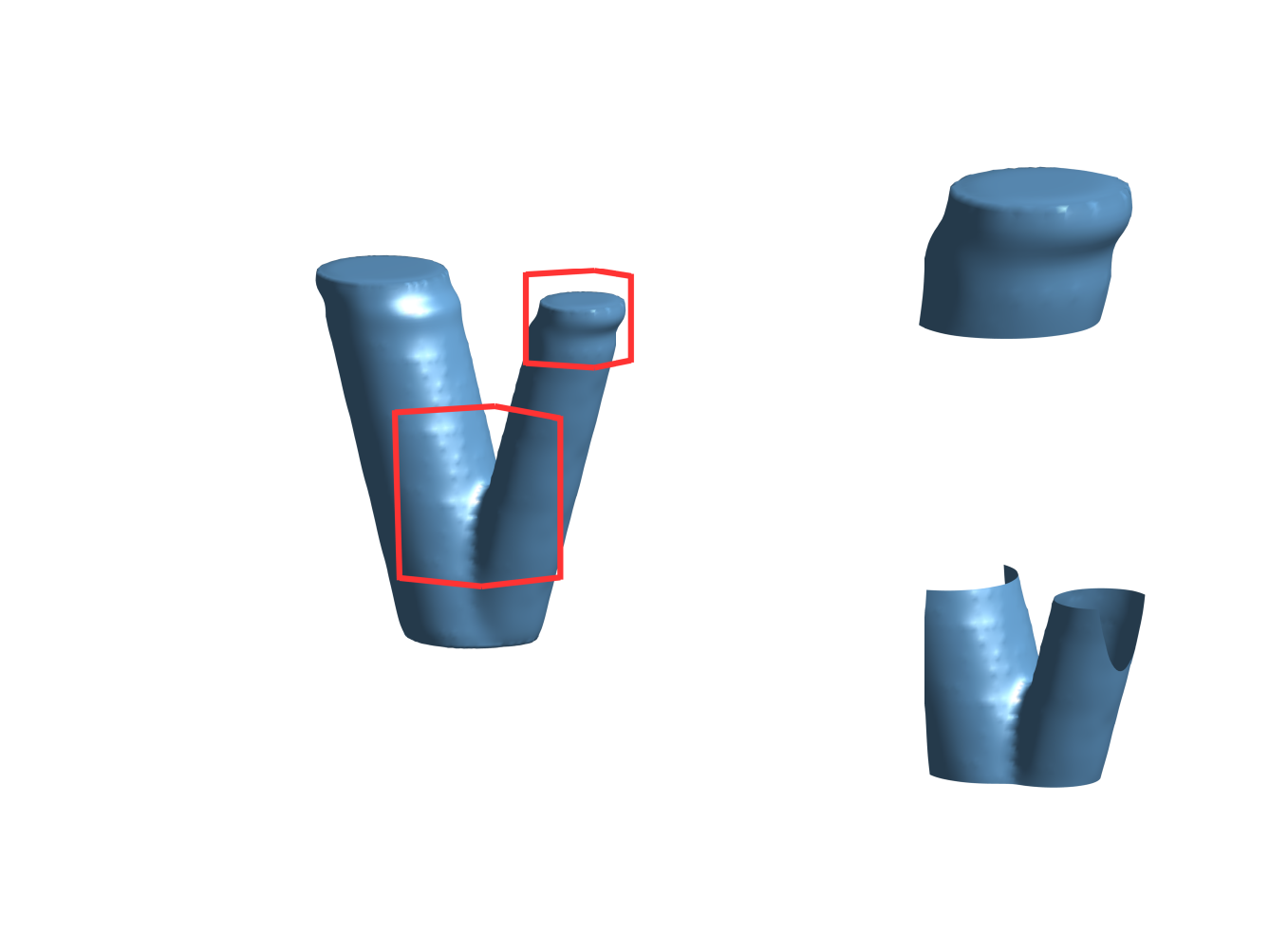}
    };

    \draw[red!80, ultra thick]
    (-0.1,2.0) -- (2,2.8);
    \draw[red!80, ultra thick]
    (-0.1,1.15) -- (2,1.15);

    \draw[red!80, ultra thick]
    (-0.8,0.6) -- (2,-0.4);
    \draw[red!80, ultra thick]
    (-0.8,-1.0) -- (2,-2.8);
    \end{tikzpicture}
    \end{subfigure}
\vspace{-1.8cm}
\cprotect \caption{Reconstruction results with defects (top right: flattened top with jagged edges and sunken shape; bottom: serious bumps) by the perimeter-based formulation of Branching Cylinders from low-resolution inputs $N = 128$ with $24$ uneven collected input slices (the left-hand side of~\Cref{p20220215_BC_Slices_EE}). }
\label{p20220731_BC_N128S24_PWEE_ZoomIn_Per}
\end{figure}
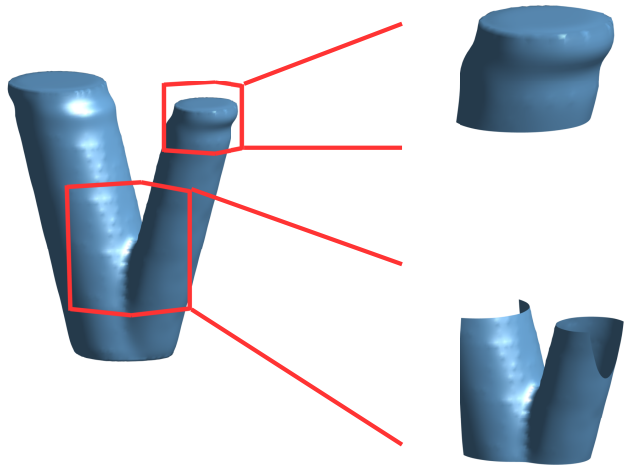

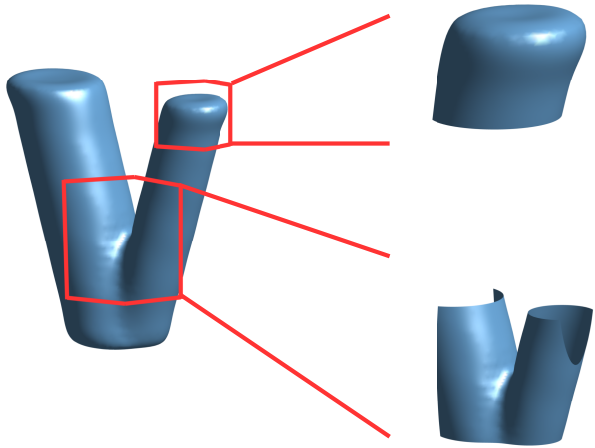
\begin{figure}[htbp]
    \centering
    \begin{subfigure}[t]{0.95\textwidth}
    \begin{tikzpicture}
    \node
    {\includegraphics[width=\textwidth]
    {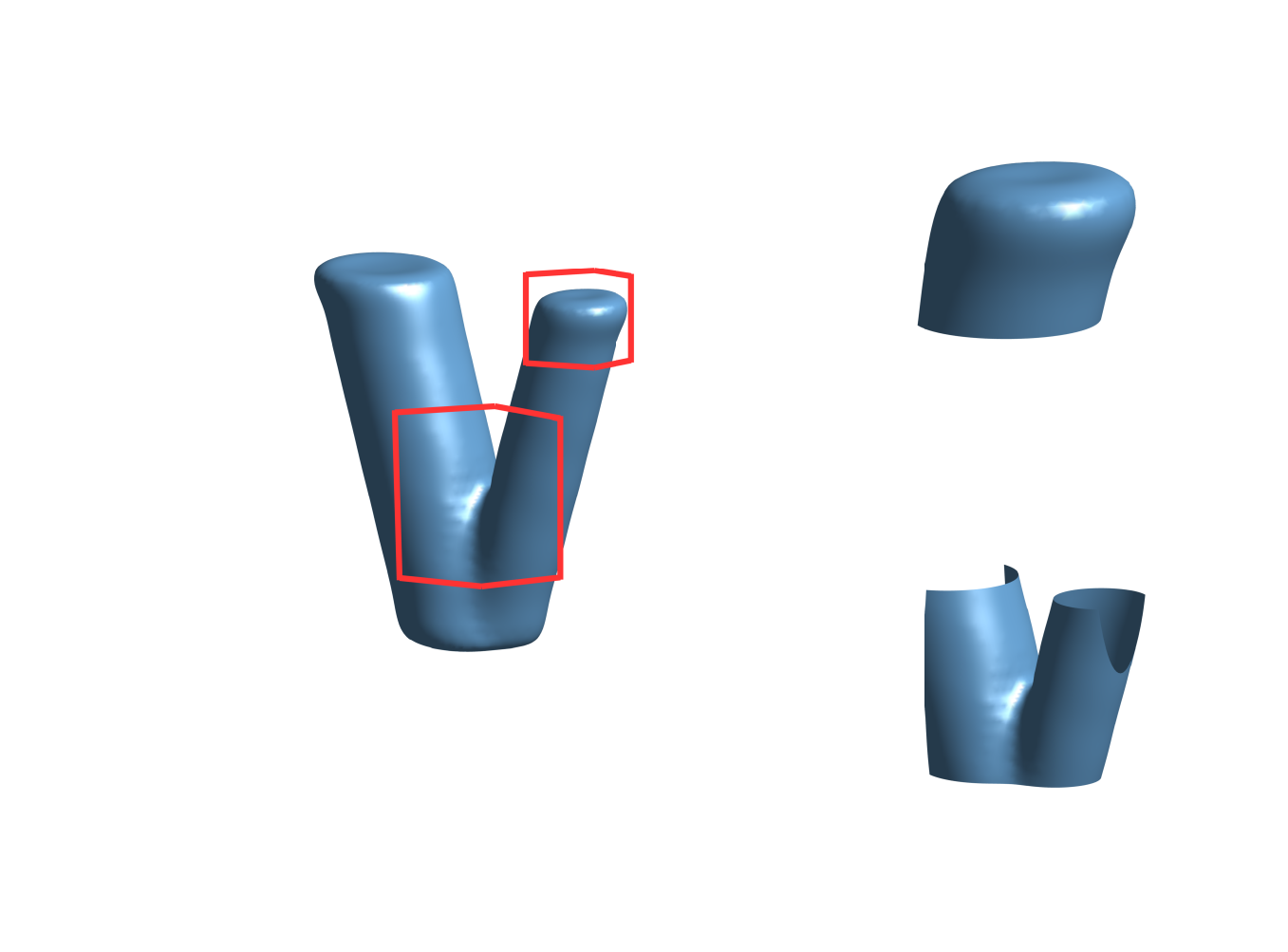}
    };

    \draw[red!80, ultra thick]
    (-0.1,1.9) -- (2,2.8);
    \draw[red!80, ultra thick]
    (-0.1,1.1) -- (2,1.1);

    \draw[red!80, ultra thick]
    (-0.8,0.55) -- (2,-0.4);
    \draw[red!80, ultra thick]
    (-0.8,-0.9) -- (2,-2.8);
    \end{tikzpicture}
    \end{subfigure}
\vspace{-1.8cm}
\cprotect \caption{Reconstruction results with defects (top right: bulging top with dented centre and jagged edges slightly; bottom: slight bumps) by the Willmore-based formulation of Branching Cylinders from low-resolution inputs $N = 128$ with $24$ uneven collected input slices (the left-hand side of~\Cref{p20220215_BC_Slices_EE}). }
\label{p20220731_BC_N128S24_PWEE_ZoomIn_Will}
\vspace{-0.2cm}
\end{figure}

This paper proposes a method in category (ii) inspired by the lucid reconstructed framework delineated by Bretin, Dayrens, and Masnou~\cite{01_Volume_Reconstruction_from_Slices} by applying geometric variational energies with the phase-field approximation, which are more pleasant to work with numerically than category (i) methods.
Two regularisation energies (perimeter-based energy and Willmore-based energy) were applied in the framework by~\cite{01_Volume_Reconstruction_from_Slices}.
Connecting with the phase-field approximation, the perimeter-based formulation (Van der Waals-Cahn-Hilliard energy) and the Willmore-based formulation are reformulated where the linear obstacle slack restrictions are constituted by the low-resolution input slices as the shape constraint. 
Nevertheless, there are some defects in these two proposed formulations as Figures~\ref{p20220731_BC_N128S24_PWEE_ZoomIn_Per}-\ref{p20220731_BC_N128S24_PWEE_ZoomIn_Will} show. 
\Cref{p20220731_BC_N128S24_PWEE_ZoomIn_Per} shows the final reconstruction result using the perimeter-based formulation, and it does not fully satisfy our visual criteria globally.
Although it successfully preserves edge features, particularly on the flat tops of the input slices, due to its energy being connected with the total variation, the overall result still contains defects with jagged edges and other undesirable artefacts.
On the other hand, \Cref{p20220731_BC_N128S24_PWEE_ZoomIn_Will} displays the result using the Willmore-based formulation, which has fewer defects compared to the perimeter-based formulation, resulting in smoother and more natural surfaces.
However, it tends to plump up the top plane due to the close relationship of Willmore energy with mean curvature, and there may be bumps if the number of input slices is not sufficiently large in the given lower resolution.
As such, the challenge becomes how to obtain a suitable variational model beyond the energy using mean curvature regularisation that can better address these issues.

Therefore, to improve on~\cite{01_Volume_Reconstruction_from_Slices}, we must employ some geometric regulariser that outperforms the mean curvature.
Motivated by other image processing works~\cite{Fast_iterative_algorithms_for_solving_the_minimization_of_curvature-related_functionals_in_surface_fairing, 20211230_Gamma-convergence_results_for_phase-field_approximations_of_the_2D-Euler_Elastica_Functional} where the Euler-Elastica-based formulations are better than the formulation related to mean curvature, our proposed work overcomes above deficiencies by minimising the Euler-Elastica-based energy.
Apart from the question of how to solve the new formulation, an interesting problem arises in comparing different models: how to deal with the different or inconsistent numbers of triangular meshes for final surfaces for a fair and objective comparison of results by different formulations?
We address this by computing the standard deviation of Gaussian curvatures (GC) and mean curvatures (MC) to indicate the corresponding level of smoothness from the viewpoint of discrete geometry, which is stimulated by~\cite{Discrete_Differential-Geometry_Operators_for_Triangulated_2-Manifolds, 20201014_A_novel_method_for_surface_mesh_smoothing:_Applications_in_biomedical_modeling}.

\vspace{-.5cm}

The rest of this paper is organised as follows:
\Cref{sec:mathspre} introduces some essential mathematical preliminaries, including definitions and notations for the reconstructed framework of a geometric variational method, linear obstacle slack restrictions, and phase-field approximation.
\Cref{sec:EE} proposes the new Euler-Elastica-based formulation first and derives the Euler–Lagrange equation, and then presents two numerical algorithms for solving the resulting optimisation problem where we extend the alternating direction method of multipliers (ADMM) to solve our model, which leads to faster numerical approximations.
Finally, simulated and realistic examples are depicted and compared in~\Cref{sec:Results} by different models (inpainting models and three formulations), where quantitative comparisons are also given to show the effectiveness of the new Euler-Elastica-based formulation.

\section{Mathematical preliminaries}
\label{sec:mathspre}

\begin{figure}[htbp]
    \centering
    \begin{subfigure}[b]{0.3\textwidth}
    \includegraphics[width=\textwidth]{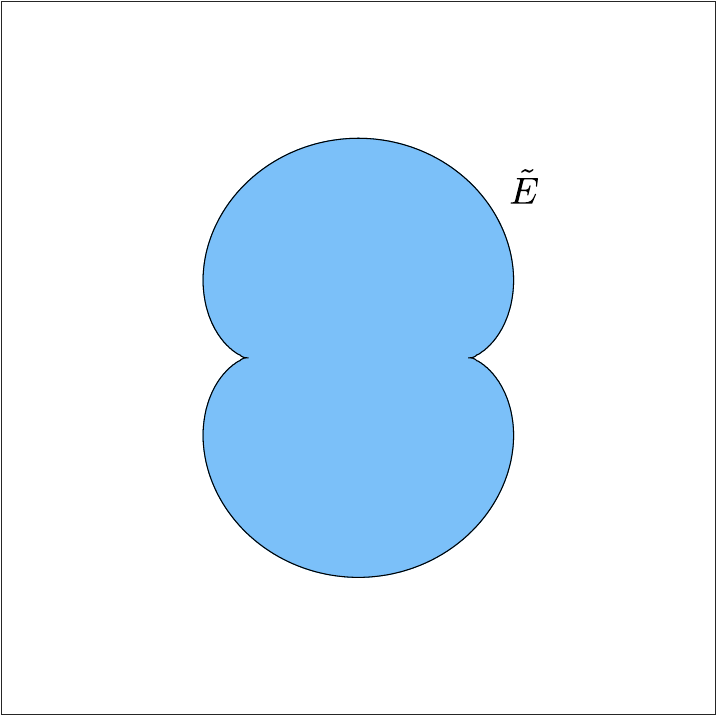}
    \caption{The target set $\tilde{E}$}
    \label{fig:3_E}
    \end{subfigure}
\quad
    \begin{subfigure}[b]{0.3\textwidth}
    \includegraphics[width=\textwidth]{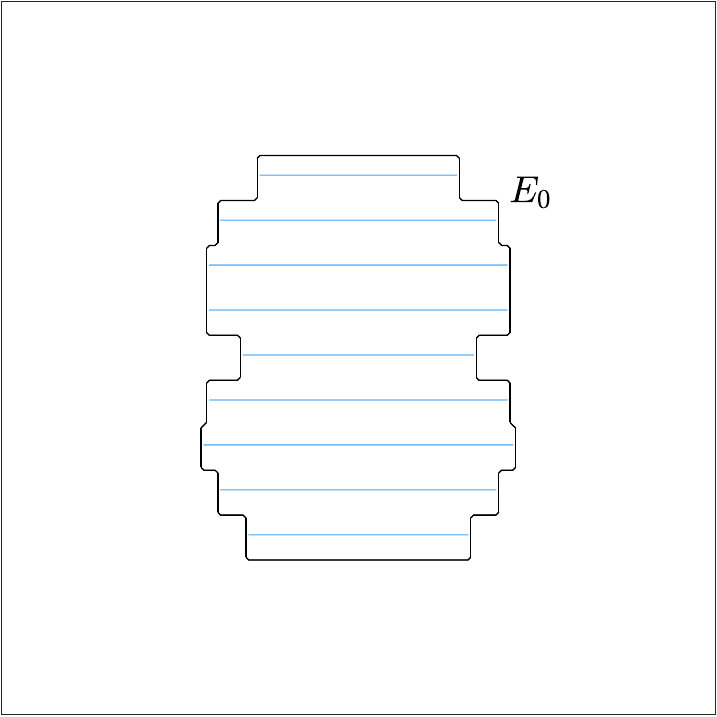}
    \caption{The initial set $E_{0}$}
    \label{fig:3_E0}
    \end{subfigure}
\quad
    \begin{subfigure}[b]{0.3\textwidth}
    \includegraphics[width=\textwidth]{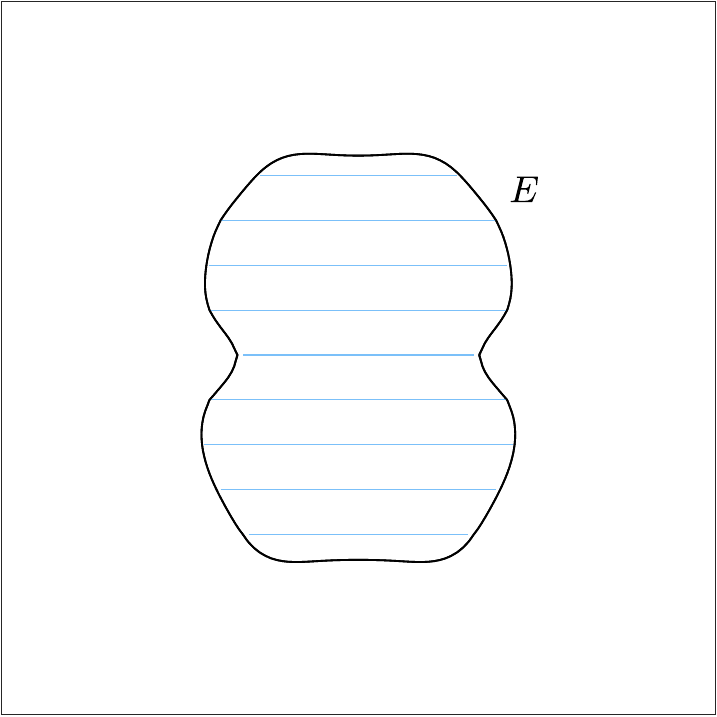}
    \caption{Possible result $E$}
    \label{fig:3_Estar}
    \end{subfigure}
\vspace{-0.4cm}
\caption{Exemplification of the target set $\tilde{E}$ (left), the initial set $E_{0}$ (centre) and the possible construction $E$ (right) in the two-dimensional scenario.}
\vspace{-0.5cm}
\label{fig:3_E_E0_Estar}
\end{figure}

The goal is to obtain a smooth $D$-dimensional reconstruction $E^{*} \subset \R^{D}$ from the initial set $E_{0}$, as~\Cref{fig:3_E_E0_Estar} exemplified, where $\tilde{E}$ is the desired target set and $E$ denotes as the possible reconstruction.
Here the initial set $E_{0}$ consists of given parallel cross-sections $\Pi_{i} \subset \R^{D}$ with the number of slices $s \in \N^{*}$ and $i = 1, \dots, s$.
We provide visual illustrations using the two-dimensional scenario ($D = 2$ for curve smoothing) to clarify some of the notations in this section (significant notations are summarised in~\Cref{sec:appendix}).
The primary implementations for the three-dimensional scenario ($D = 3$ for surface reconstruction) are presented in~\Cref{sec:Results}.
In this section, we first introduce a shape-preserving approach by interior and exterior restrictions to maintain the shape of the input set during reconstruction.
Later, we address the challenge of formulating and computing variational energies by employing a phase-field method.
This approach allows us to  represent the variational energies using a smooth and continuous function, which can be easily discretised and optimised numerically.

\subsection{Fidelity by interior and exterior restrictions}
\label{subsec:in_ex_restriction}

\begin{figure}[t]
    \centering
    \begin{subfigure}[b]{0.45\textwidth}
    \includegraphics[width=\textwidth]{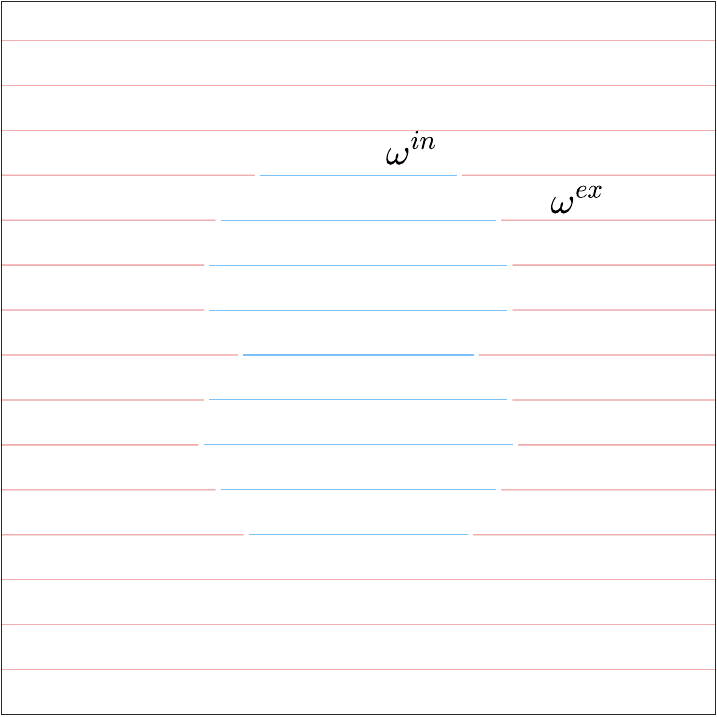}
    \caption{}
    \label{fig:0913_2_1_omega}
    \end{subfigure}
\hspace{0.1mm}
    \begin{subfigure}[b]{0.45\textwidth}
    \includegraphics[width=\textwidth]{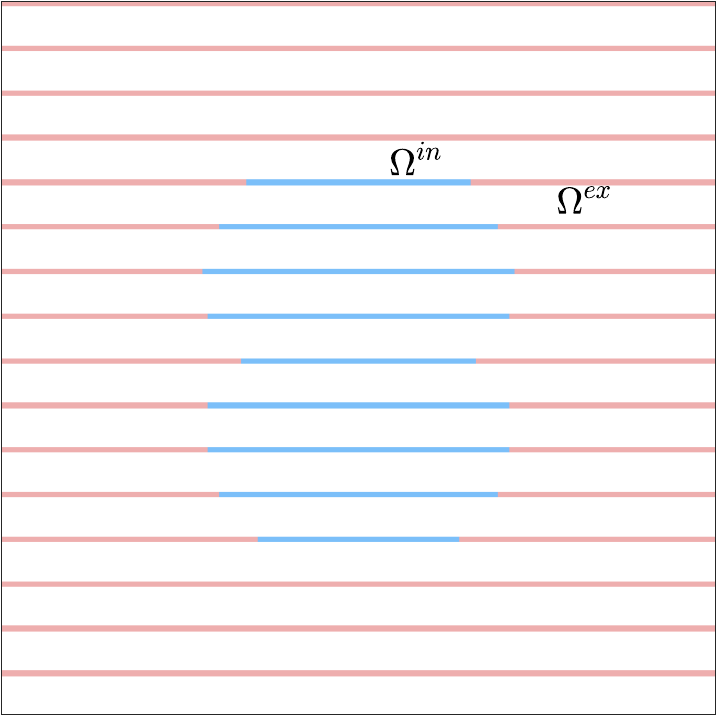}
    \caption{}
    \label{fig:0913_2_2_Omega}
    \end{subfigure}
\vspace{-0.2cm}
\caption{Illustrations of~\subref{fig:0913_2_1_omega} interior restrictions $\omega^{in}$ (cyan) and exterior restrictions $\omega^{ex}$ (pink), and of~\subref{fig:0913_2_2_Omega} fattened interior restrictions $\Omega^{in}$ (cyan) and exterior restrictions $\Omega^{ex}$ (pink) from the initial set $E_{0}$.}
\vspace{-0.3cm}
\label{fig:0913_2_omega_inex_Omega}
\end{figure}

To address the challenge of incorporating restrictions based on the initial set $E_{0}$ using a collection of $s$ finite hyperplanes $\Pi_{i}$ for $i  = 1, \dots, s$, we adopt a fidelity-driven approach to the reconstruction process.
It is worth noting that
    for any two distinct hyperplanes $i \neq j$, we assume that $\Pi_{i} \cap \Pi_{j} = \emptyset$, and
    hyperplanes can be expressed as
    $
    \Pi_{i} = \left\{ (\xi, 0): \xi \in \R^{D-1}, i = 1, \dots, s \right\}
    $
    if an appropriate orthonormal system of coordinates is chosen in $\R^{D}$ (this is possible due to slices being parallel and sitting in a low dimension).
To preserve the shape of the input data, we use a term in this approach similar to the fidelity term used in other variational models.
However, due to the discrete nature of the input data in $E_{0}$, we do not strictly use the least-squares fidelity term, such as $\int_{\Omega} |E^{*} - E_{0}|^{2} \, \id \Omega$, in the following setting.
Specifically, we first introduce interior and exterior restrictions $\omega^{in}:=\omega^{in}_{E_{0}}$ and $\omega^{ex}:=\omega^{ex}_{E_{0}}$ respectively, which are originated from the set of hyperplanes related to the initial set $E_{0}$
\begin{eqnarray*}
\omega^{in}
:=
\omega^{in}_{E_{0}}
=
\bigcup\limits_{i = 1}^{s} \omega_{i}^{in}
\quad \mbox{and} \quad
\omega^{ex}
:=
\omega^{ex}_{E_{0}}
=
\bigcup\limits_{i = 1}^{s} \omega_{i}^{ex}
\end{eqnarray*}
where $\omega_{i}^{in}$, $\omega_{i}^{ex}$ are preset based on $\Pi_{i}$.
Here,
\begin{itemize}
    \item
    $\omega^{in}, \omega^{ex} \subset \{\Pi_{i}\}$ as~\Cref{fig:0913_2_1_omega} illustrated, i.e. $\omega_{i}^{in}, \omega_{i}^{ex} \subset \Pi_{i}$,

    \item
    $\omega_{i}^{in}
    \subset E_{0} \cap \Pi_{i}
    \subset \Pi_{i} \setminus \omega_{i}^{ex}$ satisfies for every initial object $E_{0}$, and

    \item
    $\omega_{i}^{in} \cup \omega_{i}^{ex} \supset \Pi_{i}$ does not require as it enables to handle noisy inputs.
\end{itemize}
Then, the problem for reconstruction from slices can be formulated as a variational model, which is to find the (local) minimum $E^{*}$ subject to
\begin{eqnarray}
\label{equ:JE}
    E^{*} = \mathop{\mbox{argmin}}
    _{ \substack{
       \omega^{in} \subset E \\
       E \cap \omega^{ex} = \emptyset}
     } \mathcal{E}(E)
\end{eqnarray}
where $\omega^{in}, \omega^{ex}$ are interior and exterior restrictions from the initial set $E_{0}$, and $\mathcal{E}$ can be the perimeter-based $\mathscr{P}$, the Willmore-based $\mathscr{W}$, and the Euler-Elastica-based $\mathscr{E}$ energy, as shown shortly in~\Cref{subsec:approx_phasefield} to impose smoothness on the surface $\partial E$.

\begin{figure}[H]
\centering
\begin{tikzpicture}[scale=.9]
\tikzstyle{every node}=[scale=.9]
\fill[gray!60]  (-5,3.37) -- (-2,3.37) --(-2,0.63) -- (-5,0.63);

\draw[ultra thick] (-5,2) -- (-1.75,2);
\draw[ultra thick, dashed] (-1.75,2) -- (-1,2);
\filldraw (-2.5,2) circle (1.5pt);
\draw[thick] (-2.5,2) -- (-2.5,3.37);
\draw[<->] (-5,1.75) -- (-2.5,1.75);
\draw[dashed] (-2.5,3.37) -- (-1,3.37);
\draw[<->] (-2.25,2.08) -- (-2.25,3.29);
\node at (-2.3,1.7) {$\xi$};
\node at (-3.7,1.5) {$|\mathpzc{d}_{i}(\xi, \omega_{i})|$};
\node at (-0.65,2) {$\omega_{i}$};
\node at (-1.35,2.75) {$h|\mathpzc{d}_{i}(\xi, \omega_{i})|$};
\node at (-2.49,0.98) {$\Omega_{i, h}$};
\end{tikzpicture}
\vspace{-0.3cm}
\caption{Illustration of the fattened restriction of $\omega_{i}$ with the thickness parameter $h$ where part of the fattened region $\Omega_{i, h} = \left\{
(\xi, \zeta) \in (\R^{D-1} \times \R) \cap E_{0}:
\xi \in \omega_{i},
|\zeta| < h |\mathpzc{d}_{i}(\xi, \omega_{i})|
\right\}$ is shown in gray. }
\label{fig:fatten_restriction}
\end{figure}
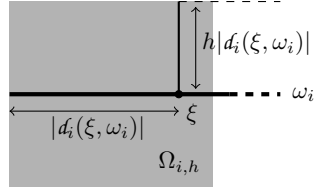

Next, we suggest two fattened restrictions depicted in~\Cref{fig:0913_2_2_Omega} to enhance more flexibility during the reconstruction process.
As we only have a limited number of slices, these restrictions are designed to enlarge the feasible region and improve the reconstruction.
Define the fattened interior and exterior restrictions as follows:
\begin{eqnarray*}
\Omega^{in}
:=\Omega^{in}_{E_{0}}
=
\bigcup\limits_{i = 1}^{s} \Omega_{i, h}^{in}
\quad \mbox{and} \quad
\Omega^{ex}
:=\Omega^{ex}_{E_{0}}
=
\bigcup\limits_{i = 1}^{s} \Omega_{i, h}^{ex}
\end{eqnarray*}
where the thickness parameter $h = \varepsilon^{\alpha}> 0$ is determined by given $\varepsilon > 0$ and $\alpha \in [0, 1]$.
Here, $h$ controls the thickness of the restrictions, and the value of $\alpha$ determines the rate at which the thickness increases with $\varepsilon$.
To define $\Omega_{i, h}^{in}$ and $\Omega_{i, h}^{ex}$, we use the signed distance function (or called oriented distance function) $\mathpzc{d}_{i}$ to an arbitrary subset $\pi_{i}$ of hyperplanes $\Pi_{i}$.
The restrictions are then defined as:
\begin{eqnarray*}
\Omega_{i, h}^{in}
=
\left\{
(\xi, \zeta) \in (\R^{D-1} \times \R) \cap E_{0}:
\xi \in \omega_{i}^{in},
|\zeta| < h |\mathpzc{d}_{i}(\xi, \omega_{i}^{in})|
\right\},
\\
\Omega_{i, h}^{ex}
=
\left\{
(\xi, \zeta) \in (\R^{D-1} \times \R) \cap E_{0}:
\xi \in \omega_{i}^{ex},
|\zeta| < h |\mathpzc{d}_{i}(\xi, \omega_{i}^{ex})|
\right\},
\end{eqnarray*}
as~\Cref{fig:fatten_restriction} illustrated where they are fixed from initial set $E_{0}$.
Here, the signed distance function $\mathpzc{d}_{i}$ is given by
$
\mathpzc{d}_{i}(\xi, \pi_{i})
=
dist(\xi, \pi_{i})
-
dist(\xi, \Pi_{i} \setminus \pi_{i})
$
for arbitrary $\xi \in \Pi_{i}$ and $\pi_{i} \subset \Pi_{i}$
where $dist(\cdot, \cdot)$ is the representative Euclidean distance in $\R^{D}$~\cite{09_Elementary_linear_algebra}.

\vspace{-0.3cm}
\vspace{-0.35cm}

\noindent
Therefore, incorporating the fattened restrictions, model~\eqref{equ:JE} becomes
\begin{eqnarray}
\label{equ:JE_fatten}
    E^{*} = \mathop{\mbox{argmin}}
    _{ \substack{
       \Omega^{in} \subset E \\
       E \cap \Omega^{ex} = \emptyset}
     } \mathcal{E}(E)
\end{eqnarray}
where fattened restrictions $\Omega^{in}, \Omega^{ex}$ are defined from $E_{0}$, and $\mathcal{E}$ as a general the energy function can be $\mathscr{P}, \mathscr{W}, \mathscr{E}$ (to be introduced shortly in~\Cref{subsec:approx_phasefield}).


\subsection{Relaxation by phase-field approximation}
\label{subsec:approx_phasefield}

The difficulty with model~\eqref{equ:JE_fatten} is that the unknown is the set $E$, not a function representing $E$.
We now introduce the phase-field method to represent the set by a function, and then turn the problem into a variational problem before considering discretisation and numerical solution.
The phase-field method is a widely used technique for modelling complex physical phenomena with sharp interfaces or discontinuities, such as phase transitions, fracture, and grain growth, among others.
In essence, the phase-field method replaces the sharp interface or discontinuity with a diffuse interface or transition zone of finite thickness, which is controlled by a scalar parameter known as the phase-field variable.
By doing so, the problem becomes amenable to standard numerical methods for solving partial differential equations, such as finite differences, finite elements, or spectral methods, among others.
In the following, we will describe how to use the phase-field method to approximate the binary indicator function and reformulate the inverse problem~\eqref{equ:JE_fatten} involving the phase-field approximation.

\begin{figure}[htbp]
\centering
    \begin{subfigure}[b]{0.4\textwidth}
    \hspace{-.8cm}
        \begin{tikzpicture}[scale=0.9]
        \tikzstyle{every node}=[scale=1]
        \draw  (-3,1.5) rectangle (3.4,-3);
        \draw[thick, fill=gray!60]  plot[smooth, tension=.7] coordinates {(-1.5,0.5) (0,1) (0,0) (1,0.5) (1.5,0) (2,-1) (2,-2) (1,-2.5) (-1,-2) (-2,-1.5) (-1.5,-1) (-1.5,-0.5) (-2.5,0) (-1.5,0.5)};
        \node at (0.5,-1.5) {$E$};
        \node at (2,-2.5) {$\partial E$};
        \node at (-0.5,-0.7) {$\mathpzc{d} < 0$};
        \node at (2,1) {$\mathpzc{d} > 0$};
        \node at (2.6,-2.1) {$\mathpzc{d} = 0$};
        \end{tikzpicture}
    \caption{}
    \label{fig:p2}
    \end{subfigure}
\quad
    \begin{subfigure}[b]{0.4\textwidth}
        \begin{tikzpicture}[scale=0.9]
        \tikzstyle{every node}=[scale=1]
        \draw  (-3,1.5) rectangle (3.4,-3);
        \draw[thick, fill=gray!60]  plot[smooth, tension=.7] coordinates {(-1.5,0.5) (0,1) (0,0) (1,0.5) (1.5,0) (2,-1) (2,-2) (1,-2.5) (-1,-2) (-2,-1.5) (-1.5,-1) (-1.5,-0.5) (-2.5,0) (-1.5,0.5)};
        \node at (0.5,-1.5) {$\Omega^{in}$};
        \node at (2.3,0.4) {$\Omega^{ex}$};
        \node at (2,-2.5) {$\partial \Omega$};
        \node at (-0.5,-0.7) {$u = 1$};
        \node at (2,1) {$u = 0$};
        \node at (2.6,-2.1) {$u = \frac{1}{2}$};
        \end{tikzpicture}
    \caption{}
    \label{fig:p2_next}
    \end{subfigure}
\vspace{-0.2cm}
\caption{Illustrations of~\subref{fig:p2} the signed distance function (negative inside $E$, zero on the boundary, and positive outside), and of~\subref{fig:p2_next} the values of $u$ for the interior region $\Omega^{in}$, boundary $\partial \Omega$, and the exterior region $\Omega^{ex}$ of the object $E$ using the profile function $q$. }
\end{figure}
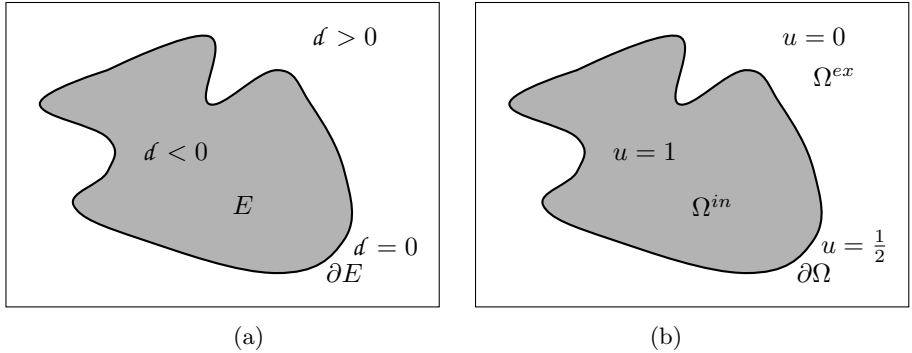

In order to describe the representation of the set $E$ by the phase-field function (recall that $E_{0}$ is the initial set and $E$ is the possible solution), Modica and Mortola in~\cite{06_Il_limite_nella_convergenza_di_una_famiglia_di_funzionali_ellittici}, as well as Bretin, Dayrens, and Masnou in~\cite{01_Volume_Reconstruction_from_Slices, 30_Approximation_par_champ_de_phase_de_mouvement_par_courbure_moyenne_anisotrope} have proposed a sequence $(u_{ \varepsilon})$ (i.e. phase-field function $u_{\varepsilon}(\cdot)$) defined as follows
\begin{eqnarray}
\label{equ:phase-field_function}
u_{\varepsilon}(\xi)
=
q
\left(
\frac{\mathpzc{d}(\xi, E)}{\varepsilon}
\right).
\label{equ:Uexfn}
\end{eqnarray}
This sequence is used to approximate the indicator function of set $E$ that characterises the interface between the target object and background region.
Here,
\begin{itemize}
    \item
    $\mathpzc{d}(\xi, E)$ is the signed distance function that measures the distance between the given point $\xi$ and the boundary $\partial E$ of the set $E$. Specifically,
    \begin{eqnarray*}
        \mathpzc{d}(\xi, E)
        \begin{cases}
        < 0
        & {\text{when }}\xi\in E\setminus \partial E
        \\
        = 0
        & {\text{when }}\xi\in \partial E \\
        > 0
        & {\text{otherwise}}
        \end{cases},
    \end{eqnarray*}
    as~\Cref{fig:p2} shows. Note that $|\mathpzc{d}(\xi, E)|$ gives the shortest distance from $\xi$ to the boundary $\partial E$ with the property $|\grad \mathpzc{d}| = 1$~\cite{01_Volume_Reconstruction_from_Slices, 30_Approximation_par_champ_de_phase_de_mouvement_par_courbure_moyenne_anisotrope, 07_Level_set_based_shape_prior_segmentation}.

    \item
    $\varepsilon$ is the phase-field variable related to the thickness parameter $h = \varepsilon^{\alpha}$ as previously mentioned.

    \item
    $q$ is called the profile function\index{Profile function}, which is required to be a piecewise function that takes the value of 1 inside $E$, half on $\partial E$, and 0 outside $E$ and is designed to enforce the continuity and smoothness of the phase-field function.
\end{itemize}
Then, the phase-field approximation method provides an equivalence between fattened restrictions and phase-field functions by
\begin{equation}
\label{equ:restrictions_equivalence}
\Omega^{in} \subset E \subset \R^{D} \setminus \Omega^{ex}
\Longleftrightarrow
u_{\varepsilon}^{in}
\le
u_{\varepsilon}
\le
u_{\varepsilon}^{ex}
\end{equation}
where the phase-field profiles $u_{\varepsilon}^{in}, u_{\varepsilon}^{ex}$ by~\eqref{equ:Uexfn} are defined by
\begin{eqnarray}
\label{equ:Uhinex}
u_{\varepsilon}^{in}
= q\left( \frac{\mathpzc{d}(\xi, \Omega^{in})}{\varepsilon} \right)
\quad \mbox{and} \quad
u_{\varepsilon}^{ex}
= 1 - q\left( \frac{\mathpzc{d}(\xi, \Omega^{ex})}{\varepsilon} \right),
\end{eqnarray}
and the latter of~\eqref{equ:restrictions_equivalence} is called linear obstacle restriction\index{Linear obstacle restriction}s on $u_{\varepsilon}$, which serves as an obstacle for the minimisation problem and allows the minimiser to satisfy the restrictions.


Next, to determine the profile function $q$ that satisfied the above requirements of being 1 inside $E$, half on $\partial E$, and 0 outside $E$, one suggestion is to utilise the double-well potential $W$~\cite{01_Volume_Reconstruction_from_Slices}.
The Euler equation $q'' = W'(q)$ with the initial condition $q(0)$ can be solved as the following Cauchy problem for $W \in \mathbb{C}^{2}$
\begin{eqnarray*}
       q' = -\sqrt{2W(q)}
       \quad \mbox{and} \quad
       q(0) = \frac{1}{2}.
\end{eqnarray*}
For the specific case of the double-well potential\index{Double-well potential} given by
$
    W(u) = \frac{1}{2}u^{2}(1-u)^{2}
$,
the profile function $q$ can be deduced via the separation of variables as
\begin{eqnarray}
    q(\xi)
    =
    \frac{1}{2}
    \left(
    1 - \tanh{
              \left(
              \frac{\xi}{2}
              \right)
              }
    \right)
    = \frac{1}{1+e^{\xi}}.
\label{equ:q_profilefn}
\end{eqnarray}
Then, by using the profile function $q$, the interior region $\Omega^{in}$ and exterior region $\Omega^{ex}$ of the initial object $E_{0}$ can be indicated as $1$ and $0$ respectively, and in fact, two phase-field profiles by~\eqref{equ:Uhinex} have the following convergence results
\begin{eqnarray}
\label{equ:UconvergenceResult}
    u_{\varepsilon}^{in}(\xi)
    \xrightarrow{\varepsilon \rightarrow 0}
    \begin{cases}
    1 & {\text{if }}\xi \in \Omega^{in} \\
    \frac{1}{2} & {\text{if }}\xi \in \partial \Omega^{in} \\
    0 & {\text{otherwise}}
    \end{cases}
    \quad \mbox{and} \quad
    u_{\varepsilon}^{ex}(\xi)
    \xrightarrow{\varepsilon \rightarrow 0}
    \begin{cases}
    0 & {\text{if }}\xi \in \Omega^{ex} \\
    \frac{1}{2} & {\text{if }}\xi \in \partial \Omega^{ex} \\
    1 & {\text{otherwise}}
    \end{cases}
\end{eqnarray}
as~\Cref{fig:p2_next} sketched.
Furthermore, these profiles can be indicated by
\begin{eqnarray}
u_{\varepsilon}^{in}
= \frac{1}{2} \mathbbm{1}_{\Omega^{in}}
\quad \mbox{and} \quad
u_{\varepsilon}^{ex}
= 1 - \frac{1}{2} \mathbbm{1}_{\Omega^{ex}}.
\label{equ:Uinout_appro}
\end{eqnarray}

Naturally, the variational framework~\eqref{equ:JE_fatten} with the phase-field relaxation is reformulated as:
\begin{eqnarray}
    u_{\varepsilon}^{*}
    =
    \mathop{\mbox{argmin}}
    _{ u^{in}_{E_{0}}
    \le
    u
    \le
    u^{ex}_{E_{0}}}
    \mathcal{E}_{\varepsilon}(u),
\label{equ:Jue}
\end{eqnarray}
where $u:= u_{\varepsilon}, u^{in}_{E_{0}}:= u_{\varepsilon}^{in}, u^{ex}_{E_{0}}:= u_{\varepsilon}^{ex}$ after omitting $\varepsilon$ for simplified, two restrictions are associated with the initial set $E_{0}$, and $\mathcal{E}_{\varepsilon}$ can be the perimeter-based $\mathscr{P}_{\varepsilon}$, Willmore-based $\mathscr{W}_{\varepsilon}$ and Euler-Elastica-based $\mathscr{E}_{\varepsilon}$ formulation from corresponding energies $\mathcal{E} = \mathscr{P}, \mathscr{W}, \mathscr{E}$, to be introduced in~\eqref{equ:perimeter_energy}-\eqref{equ:Euler-Elastica_energy} for energies, and~\eqref{eqn:P_model},~\eqref{eqn:W_model},~\eqref{eqn:EE_model} for formulations.


Our framework incorporates three classical energies: the perimeter-based $\mathscr{P}$, the Willmore-based $\mathscr{W}$, and the Euler-Elastica-based $\mathscr{E}$ energy~\cite{01_Volume_Reconstruction_from_Slices, 20211230_Gamma-convergence_results_for_phase-field_approximations_of_the_2D-Euler_Elastica_Functional, 20220104_On_a_Modified_Conjecture_of_De_Giorgi, 20220104_A_singular_perturbation_problem_with_integral_curvature_bound}.
These energies with their properties have been extensively studied in image processing and are now extended to our surface reconstruction problem.
\begin{description}
    \item[($\mathscr{P}$)]
    The perimeter-based energy is the first energy we consider, expressed as
    \begin{eqnarray}
    \mathscr{P}(E)
    =
    \int_{\partial E}
    \mathbbm{1}_{E} \, \id \mathcal{H}^{D-1}
    \label{equ:perimeter_energy}
    \end{eqnarray}
    where $\mathbbm{1}_{E}$ is alluded to the indicator function\index{Indicator function} or the characteristic function of the set $E$ that indicates elements in the set $E$~\cite{30_Approximation_par_champ_de_phase_de_mouvement_par_courbure_moyenne_anisotrope, 05_Surface_Reconstruction_from_Discrete_Indicator_Functions}, that is, $\mathbbm{1}_{E}(\xi) = 1$ if $\xi\in E$ and $0$ otherwise,
    $\mathcal{H}^{D-1}$ denotes the $(D-1)$-dimensional Hausdorff measure in $\mathbb{R}^{D}$ and $\partial E$ is the boundary of $E$.
    The perimeter-based energy is a simple and intuitive way to measure the length or surface area of a given set.
    One of the key properties of this energy is that it is scale-invariant, meaning that it does not depend on the size or location of the set being measured.
    Additionally, it has a well-defined gradient that can be used for optimisation purposes.
    However, it can be sensitive to noise and can be affected by irregularities in the boundary of the set.
    Besides, there is a relationship between the perimeter energy and the total variation (TV) functional.
    In fact, the perimeter energy can be viewed as a special case of the TV functional, where the signal or image is a binary indicator function of a set.
    More generally, the TV functional can be seen as a generalisation of the perimeter energy to functions that are not binary indicator functions.

    \item[($\mathscr{W}$)]
    For the second choice of energy: Willmore-based energy is defined by
    \begin{eqnarray}
    \mathscr{W}(E)
    =
    \int_{\partial E}
    \left|
    H
    \right|^{2}
    \, \id \mathcal{H}^{D-1}
    \label{equ:Willmore_energy}
    \end{eqnarray}
    where $H$ is the mean curvature on the boundary $\partial E$.
    The Willmore energy is a geometric energy that measures the bending or deviation of a surface from a plane.
    This energy is quadratic in the mean curvature of the surface, which captures both its local and global curvature information.
    The Willmore energy is also scale-invariant, which means that it is preserved under rescaling of the surface, making it robust to changes in size or orientation.

    \item[($\mathscr{E}$)]
    \hspace{0.03cm}
    Last, the Euler-Elastica-based energy is the key suggestion in our framework, being the third choice, which is the combination of two energies mentioned above
    \begin{eqnarray}
    \mathscr{E}(E)
    =
    \mathscr{P}(E)
    +
    \mathscr{W}(E)
    =
    \int_{\partial E}
    (
    \mathbbm{1}_{E}
    +
    \left|
    H
    \right|^{2}
    )
    \, \id \mathcal{H}^{D-1}
    \label{equ:Euler-Elastica_energy}
    \end{eqnarray}
    with the indicator function $\mathbbm{1}_{E}$ and the mean curvatures $H$ on the boundary $\partial E$ as~\eqref{equ:perimeter_energy} and~\eqref{equ:Willmore_energy}.
    Clearly, the Euler-Elastica-based energy combines the perimeter-based energy and the Willmore-based energy, which makes it an effective energy functional for shape analysis and geometric modelling.
    The first indicator term measures the smoothness of the boundary, while the second term measures the curvature.
    By combining these two terms, the Euler-Elastica energy can capture both the local and global properties of a shape.
    Moreover, the Euler-Elastica energy has been shown to have desirable properties, such as convexity and stability, which makes it an attractive option for optimisation-based approaches.
\end{description}

Furthermore, in recent years, new formulations of these energies have been explored using the approaches of $\Gamma$-convergence and phase-field approximation.
Notably, the $\Gamma$-convergence and phase-field approximation approaches have opened up new avenues for studying these energies and their properties, such as convergence, stability, and regularity.
Moreover, the use of regularised functions in these formulations can be beneficial in numerical optimisation, as they allow for efficient computation of minimisers and can avoid issues associated with non-smooth functionals.
These developments have led to significant progress in the field of geometric variational problems and their applications.
\begin{description}
\item[($\mathscr{P}_{\varepsilon}$)]
Modica and Mortola in~\cite{06_Il_limite_nella_convergenza_di_una_famiglia_di_funzionali_ellittici} along with additional literature~\cite{01_Volume_Reconstruction_from_Slices, 20211230_Gamma-convergence_results_for_phase-field_approximations_of_the_2D-Euler_Elastica_Functional} revealed that the perimeter-based energy $\mathscr{P}$ could be approximated using the classical phase-field formulation (Van der Waals-Cahn-Hilliard formulation, in accordance with other terminologies, say perimeter-based formulation $\mathscr{P}_{\varepsilon}$ instead)
\begin{eqnarray}
\label{eqn:P_model}
    \mathscr{P}_{\varepsilon}(u)
    =
    \int_{\Omega}
    \left(
    \frac{\varepsilon}{2} |\nabla u|^{2}
    +
    \frac{1}{\varepsilon} W(u)
    \right)
    \, \id \Omega
\end{eqnarray}
where $W(u) = \frac{1}{2}u^{2}(1-u)^{2}$ and $\varepsilon$ is the representative diffuse interface width.
More specifically, the $\Gamma(L^{1}(\Omega))$-limit is the $\Gamma$-convergence of $\mathscr{P}_{\varepsilon}$ to the area functional as $\varepsilon \to 0$.
In other words,
\begin{eqnarray}
    \Gamma(L^{1}(\Omega))
    -
    \lim_{\varepsilon \to 0} \mathscr{P}_{\varepsilon}(u)
    =
    \mathscr{P}(u)
    =
    \int_{\Omega} \, \mathrm{d} |\nabla u|
\end{eqnarray}
for $u \in \operatorname{BV}(\Omega, \{0, 1\})$  where $u = 2\chi_{\Omega} - 1 = \mathbbm{1}_{\Omega}$ with the characteristic function $\chi_{\Omega}$ of the finite perimeter set $\Omega$.
Then, there exists
$\mathscr{P}(u)
=
\int_{\Omega} \, \mathrm{d} |\nabla u|
=
\mathcal{H}^{D-1} (\partial \Omega)
=
\int_{\partial \Omega}
\mathbbm{1}_{\Omega} \, \id \mathcal{H}^{D-1}$.

\item[($\mathscr{W}_{\varepsilon}$)]
Further, they also demonstrated that the Willmore-based formulation $\mathscr{W}_{\varepsilon}$ can be used to approximate ($\leftrightsquigarrow$) the Willmore-based energy $\mathscr{W}$~\cite{01_Volume_Reconstruction_from_Slices, 06_Il_limite_nella_convergenza_di_una_famiglia_di_funzionali_ellittici, 30_Approximation_par_champ_de_phase_de_mouvement_par_courbure_moyenne_anisotrope}, i.e.
\begin{eqnarray}
\label{eqn:W_model}
\begin{aligned}
&
\mathscr{W}(E)
=
\int_{\partial E}
\left|
H
\right|^{2}
\, \id \mathcal{H}^{D-1}
\\
\leftrightsquigarrow \quad
&
\mathscr{W}_{\varepsilon}(u)
=
\frac{1}{2\varepsilon}
\int_{\Omega}
\left(
\varepsilon \laplace u
-
\frac{1}{\varepsilon} W'(u)
\right)^{2}
\, \id \Omega
\label{eqn:W_model_b}
\end{aligned}
\end{eqnarray}
where $W'(u) = u(u-1)(2u-1)$.
Remark that $\mathscr{W}_{\varepsilon}(u)$ stands for the rescaled norm of $L^{2}$-gradient flow of $\mathscr{P}_{\varepsilon}(u)$~\cite{20211230_Gamma-convergence_results_for_phase-field_approximations_of_the_2D-Euler_Elastica_Functional}.

\item[($\mathscr{E}_{\varepsilon}$)]
Eventually, we consider the Euler-Elastica-based formulation
\begin{eqnarray}
\begin{aligned}
\label{eqn:EE_model}
    & \,\,
    \mathscr{E}_{\varepsilon}(u)
    =
    \mathscr{P}_{\varepsilon}(u)
    +
    \mathscr{W}_{\varepsilon}(u) \\
    =
    &
    \int_{\Omega}
    \left(
    \frac{\varepsilon}{2} |\nabla u|^{2}
    +
    \frac{1}{\varepsilon} W(u)
    \right)
    \, \id \Omega
    +
    \frac{1}{2\varepsilon}
    \int_{\Omega}
    \left(
    \varepsilon \laplace u
    -
    \frac{1}{\varepsilon} W'(u)
    \right)^{2}
    \, \id \Omega
\end{aligned}
\end{eqnarray}
to approximate the Euler-Elastica-based energy where the double-well potential $W(u) = \frac{1}{2}u^{2}(1-u)^{2}$ has two minima and its first derivative is $W'(u) = u(u-1)(2u-1)$.
Röger and Schätzle in~\cite{20220104_On_a_Modified_Conjecture_of_De_Giorgi}, partially responding to the conjecture of De Giorgi~\cite{20220104_Some_remarks_on_Gamma-convergence_and_least_squares_method}, with additional literature~\cite{20211230_Gamma-convergence_results_for_phase-field_approximations_of_the_2D-Euler_Elastica_Functional, 20220104_A_singular_perturbation_problem_with_integral_curvature_bound} proved that the approximation is established with respect to $\Gamma$-convergence of $\mathscr{E}_{\varepsilon}$
\begin{eqnarray}
\label{eqn:Euler_Elastica_Functional}
    \Gamma(L^{1}(\Omega))
    -
    \lim_{\varepsilon \to 0} \mathscr{E}_{\varepsilon}(u)
    =
    \int_{\partial \Omega}
    (\mathbbm{1}_{\Omega} + |H_{\partial \Omega}|^{2})
    \, \mathrm{d} \mathcal{H}^{D-1}
\end{eqnarray}
where $\mathbbm{1}_{E}$ is the indicator function, $H_{\partial \Omega}$ is the mean curvatures vector of $\partial \Omega$.
\end{description}

\section{The new model and its numerical algorithms}
\label{sec:EE}

After establishing the necessary mathematical framework in~\Cref{sec:mathspre}, we propose a new model related to our objective functional~\eqref{equ:Jue}, which utilises the Euler-Elastica-based formulation~\eqref{eqn:EE_model} by minimising
\begin{eqnarray}
\begin{aligned}
\label{new_model}
    &
    \mathscr{E}_{\varepsilon}(u)
    =
    \int_{\Omega}
    \left(
    \frac{\varepsilon}{2} |\nabla u|^{2}
    +
    \frac{1}{\varepsilon} W(u)
    \right)
    \, \id \Omega
    +
    \frac{1}{2\varepsilon}
    \int_{\Omega}
    \left(
    \varepsilon \laplace u
    -
    \frac{1}{\varepsilon} W'(u)
    \right)^{2}
    \, \id \Omega \\
    &
    \mbox{ s.t. }
    u^{in}_{E_{0}}
    \le
    u
    \le
    u^{ex}_{E_{0}}
\end{aligned}
\end{eqnarray}
subject to linear obstacle restrictions as~\eqref{equ:restrictions_equivalence} related to the initial set $E_{0}$ where the double-well potential $W(u) = \frac{1}{2}u^{2}(1-u)^{2}$ has two minima and its first derivative is $W'(u) = u(u-1)(2u-1)$.

\subsection{Derivation for the Euler-Elastica-based formulation}
\label{sec:derivation}

In the following, we derive the Euler–Lagrange (E-L) PDE for the Euler-Elastica-based formulation~\eqref{eqn:EE_model} in order to compute the (local) minimum of our model~\eqref{new_model}. For other two formulations~\eqref{eqn:P_model} and~\eqref{eqn:W_model_b}, the derivation is analogous.

By G\^ateaux derivative,
for $\forall \tau$, there exists
\begin{subequations}
 \label{equ:limJ4}
\begin{align}
    \delta \mathscr{E}
     =
    & \, \frac{1}{\varepsilon}
    \int_{\Omega}
    W'(u) \tau
    \, \id \Omega
    +
    \varepsilon
    \int_{\Omega}
    \grad u \cdot \grad \tau
    \, \id \Omega \label{equ:limJ3} \\
    & + \varepsilon \int_{\Omega}
    \laplace u \cdot \laplace \tau
    \, \id \Omega
     + \frac{1}{\varepsilon^{3}}
     \int_{\Omega}
     W'(u) \cdot W''(u) \tau
     \, \id \Omega \label{equ:limJ4a} \\
    & -
    \frac{1}{\varepsilon}
    \left(
    \int_{\Omega}
    \laplace \tau \cdot W'(u)
    \, \id \Omega
    +
    \int_{\Omega}
    \laplace u \cdot W''(u) \tau
    \, \id \Omega
    \right). \label{equ:limJ4b}
\end{align}
\end{subequations}
Next by Green's formulae, for the second term of~\eqref{equ:limJ3}, the first term of~\eqref{equ:limJ4a} and of~\eqref{equ:limJ4b}, there exists
\begin{eqnarray}
\label{equ:GreenJ3}
\int_{\Omega}
\grad u \cdot \grad \tau
\, \id \Omega
=
\oint_{\Gamma}
\grad u \cdot \tau \cdot \boldsymbol{n}
\, \id \Gamma
-
\int_{\Omega}
\laplace u \cdot \tau
\, \id \Omega,
\end{eqnarray}

\begin{eqnarray}
\begin{aligned}
\label{equ:GreenJ4a}
     \int_{\Omega}
    \laplace u \cdot \laplace \tau
    \, \id \Omega
    = & \int_{\Omega}
        \laplace u \cdot \grad \cdot \grad \tau
          \id \Omega
    =   \oint_{\Gamma}
        \laplace u \cdot \grad \tau \cdot \boldsymbol{n}
          \id \Gamma
     - \int_{\Omega}
        \grad \cdot (\laplace u) \cdot \grad \tau
        \id \Omega \\
    = & \oint_{\Gamma}
        \laplace u \cdot \grad \tau \cdot \boldsymbol{n}
        \, \id \Gamma
     - \oint_{\Gamma}
       \grad \cdot (\laplace u) \cdot \tau \cdot \boldsymbol{n}
       \, \id \Gamma
     + \int_{\Omega}
       \laplace^{2} u \cdot \tau
       \, \id \Omega,
\end{aligned}
\end{eqnarray}
and
\begin{eqnarray}
\begin{aligned}
\label{equ:GreenJ4b}
   & \quad \int_{\Omega}
      \laplace \tau \cdot W'(u)
       \id \Omega
    =  \int_{\Omega}
        \grad \cdot \grad \tau \cdot W'(u)
         \id \Omega  \\
   & =   \oint_{\Gamma}
        \grad \tau \cdot W'(u) \cdot \boldsymbol{n}
      \id \Gamma
    - \int_{\Omega}
      \grad \tau \cdot \grad W'(u)
      \, \id \Omega \\
   & =  \oint_{\Gamma}
        \grad \tau \cdot W'(u) \cdot \boldsymbol{n}
        \id \Gamma
    - \oint_{\Gamma}
      \tau \cdot \grad W'(u) \cdot \boldsymbol{n}
      \, \id \Gamma
    + \int_{\Omega}
      \tau \cdot \laplace W'(u)
      \, \id \Omega.
\end{aligned}
\end{eqnarray}

Therefore, assembling~\eqref{equ:limJ4},~\eqref{equ:GreenJ3},~\eqref{equ:GreenJ4a} and~\eqref{equ:GreenJ4b},
\begin{eqnarray*}
\begin{aligned}
\delta \mathscr{E}
 & =
    \frac{1}{\varepsilon}
     \int_{\Omega}
     W'(u) \tau
     \, \id \Omega
    + \varepsilon
    \left(
    \oint_{\Gamma}
    \grad u \cdot \tau \cdot \boldsymbol{n}
    \, \id \Gamma
    -
    \int_{\Omega}
    \laplace u \cdot \tau
    \, \id \Omega
    \right) \\
 & \quad
     + \varepsilon
     \left(
        \oint_{\Gamma}
        \laplace u \cdot \grad \tau \cdot \boldsymbol{n}
        \, \id \Gamma
     - \oint_{\Gamma}
       \grad \cdot (\laplace u) \cdot \tau \cdot \boldsymbol{n}
       \, \id \Gamma
     + \int_{\Omega}
       \laplace^{2} u \cdot \tau
       \, \id \Omega
       \right) \\
 & \quad
     + \frac{1}{\varepsilon^{3}}
     \int_{\Omega}
     W'(u) \cdot W''(u) \tau
     \, \id \Omega \\
 & \quad
     -
     \frac{1}{\varepsilon}
     \left(
     \oint_{\Gamma}
        \grad \tau \cdot W'(u) \cdot \boldsymbol{n}
        \, \id \Gamma
    - \oint_{\Gamma}
      \tau \cdot \grad W'(u) \cdot \boldsymbol{n}
      \, \id \Gamma
      \right. \\
 & \quad
      \left.
    + \int_{\Omega}
      \tau \cdot \laplace W'(u)
      \, \id \Omega
     +
     \int_{\Omega}
     \laplace u \cdot W''(u) \tau
     \, \id \Omega
     \right) \\
 & = \int_{\Omega}
       \left(
          \frac{1}{\varepsilon} W'(u)
          -
          \varepsilon \laplace u
          + \varepsilon \laplace^{2} u
          + \frac{1}{\varepsilon^{3}}
            W'(u) \cdot W''(u)
            \right. \\
 & \quad
        \left.
          - \frac{1}{\varepsilon}
            \laplace W'(u)
          - \frac{1}{\varepsilon}
            \laplace u \cdot W''(u)
       \right) \tau
       \, \id \Omega  +
      \oint_{\Gamma}
         \left(
              \varepsilon \laplace u \cdot \boldsymbol{n}
            - \frac{1}{\varepsilon} W'(u)  \cdot \boldsymbol{n}
         \right) \grad \tau \, \id \Gamma \\
 & \quad +
      \oint_{\Gamma}
         \left(
            \varepsilon \grad u \cdot \boldsymbol{n}
            - \varepsilon \grad \cdot (\laplace u) \cdot \boldsymbol{n}
            + \frac{1}{\varepsilon} \grad W'(u)  \cdot \boldsymbol{n}
         \right) \tau \, \id \Gamma
  = 0
\end{aligned}
\end{eqnarray*}
is permitted by the following E-L equation
\begin{eqnarray*}
\frac{1}{\varepsilon} W'(u)
-
\varepsilon \laplace u
+
\varepsilon \laplace^{2} u
+ \frac{1}{\varepsilon^{3}} W'(u) \cdot W''(u)
- \frac{1}{\varepsilon} \laplace W'(u)
- \frac{1}{\varepsilon} \laplace u \cdot W''(u) = 0,
\end{eqnarray*}
that is,
\begin{eqnarray}
    \frac{1}{\varepsilon} W'(u)
    -
    \varepsilon \laplace u
    +
    \laplace \left( \varepsilon \laplace u - \frac{1}{\varepsilon} W'(u) \right)
    - \frac{1}{\varepsilon^{2}} W''(u) \left( \varepsilon \laplace u - \frac{1}{\varepsilon} W'(u) \right) = 0.
    \label{equ:willeq1}
\end{eqnarray}

\subsection{Numerical algorithm I}
\label{sec:numappro}
Following the above derivation of an E-L PDE, consider how to construct a numerical algorithm  for our models with respect to $\mathcal{E}_{\varepsilon}$.
A Cauchy problem is first recalled~\cite{01_Volume_Reconstruction_from_Slices}
\begin{eqnarray}
\label{equ:Cauchy_Problem}
\left\{
    \begin{array}{ccl}
       u_{t}
       & = &
       -\grad \mathcal{E}_{\varepsilon}(u), \\
       u(\xi, 0)
       & = &
       u_{0}(\xi).
    \end{array}
\right.
\end{eqnarray}
Here, $\mathcal{E}_{\varepsilon}$ can be the perimeter-based $\mathscr{P}_{\varepsilon}$, Willmore-based $\mathscr{W}_{\varepsilon}$ and Euler-Elastica-based $\mathscr{E}_{\varepsilon}$ formulation, and $u_{0}(\xi)$ is from the initial set $E_{0}$ by~\eqref{equ:phase-field_function}.
Then, by the Euler semi-implicit discretisation scheme in time~\cite{11_Convexity_splitting_in_a_phase_field_model_for_surface_diffusion}, the approximate numerical scheme can be expressed with the presetting synthetic time step $\tau$
\begin{eqnarray}
    \frac{u_{k+1} - u_{k}}{\tau}
    =
    - \grad \mathcal{E}_{\varepsilon}(u_{k+1}).
\label{equ:Euler_semiimp_dis}
\end{eqnarray}
In essence, the iterative solution $u^{k+1}$ meets the regularisation with~\eqref{equ:Jue}, that is,
\begin{eqnarray*}
    u_{k+1}
    =
    \mathop{\mbox{argmin}}_{u}
    \left\{
    \frac{1}{2 \tau}
    \int_{\Omega} (u - u_{k})^{2} \, \id \Omega
    + \mathcal{E}_{\varepsilon}(u)
    \right\},
\end{eqnarray*}
so that the restrictive conditions are promised as~\eqref{equ:restrictions_equivalence} aforesaid
\begin{eqnarray*}
    \mathcal{E}_{\varepsilon}(u_{k+1})
    \leqslant
    \frac{1}{2 \tau}
    \int_{\Omega} (u_{k+1} - u_{k})^{2} \, \id \Omega
    +
    \mathcal{E}_{\varepsilon}(u_{k+1})
    \leqslant
    \mathcal{E}_{\varepsilon}(u_{k}).
\end{eqnarray*}
In addition, we enforce the linear obstacle restriction~\eqref{equ:restrictions_equivalence} by applying the orthogonal projection to handle the inequality with increased relaxation, i.e.
\begin{eqnarray}
\label{eqn:linear_obstacle_restriction_apprximate}
    u^{in}_{E_{0}}
    \leq
    u
    \leq
    u^{ex}_{E_{0}}
    \quad
    \leftrightsquigarrow
    \quad
    \max(\min(u, u^{ex}_{E_{0}}), u^{in}_{E_{0}}).
\end{eqnarray}

Accordingly, the numerical scheme of the E-L equation~\eqref{equ:willeq1} for the Euler-Elastica-based formulation with respect to the time step $\tau$ is proposed as
\begin{eqnarray*}
    u_{t}
    =
    \varepsilon \laplace u
    -
    \frac{1}{\varepsilon} W'(u)
    +
    \laplace \left(
    \frac{1}{\varepsilon} W'(u)
    -
    \varepsilon \laplace u
    \right)
    +
    \frac{1}{\varepsilon^{2}} W''(u) \left( \varepsilon \laplace u - \frac{1}{\varepsilon} W'(u) \right).
\end{eqnarray*}
Following that, its numerical Euler semi-implicit discretisation scheme in time $\tau$ is expressed
\begin{eqnarray}
\begin{aligned}
    \frac{u_{k+1} - u_{k}}{\tau}
    =
    &
    \varepsilon \laplace u_{k+1}
    -
    \frac{1}{\varepsilon} W'(u_{k+1})
    + \laplace \left(\frac{1}{\varepsilon} W'(u_{k+1}) - \varepsilon \laplace u_{k+1} \right) \\
    &
    +
    \frac{1}{\varepsilon^{2}} W''(u_{k+1}) \left( \varepsilon \laplace u_{k+1} - \frac{1}{\varepsilon} W'(u_{k+1}) \right).
    \label{equ:euler_elastica_num}
\end{aligned}
\end{eqnarray}
Clearly, reorganising above equation~\eqref{equ:euler_elastica_num}, there exists
\begin{eqnarray}
\begin{aligned}
\label{eqn:EE_iterative}
    u_{k+1}
    =
    & \left(
    I_{D}
    - \varepsilon \tau \laplace
    + \tau \varepsilon \laplace^{2}
    \right)^{-1}
    \\
    & \left(
    u_{k}
    -
    \frac{\tau}{\varepsilon} W'(u_{k+1})
    +
    \frac{\tau}{\varepsilon} \laplace W'(u_{k+1})
    \right.
    \\
    &
    \left.
    -
    \frac{\tau}{\varepsilon} W''(u_{k+1}) \laplace u_{k+1}
    +
    \frac{\tau}{\varepsilon^{3}} W'(u_{k+1}) W''(u_{k+1})
    \right)
\end{aligned}
\end{eqnarray}
that is, the point $u_{k+1}$ is also a fixed point of the function
\begin{eqnarray}
\begin{aligned}
    \varPhi_{\mathscr{E}}(x)
    =
    & \left(
    I_{D}
    - \varepsilon \tau \laplace
    + \tau \varepsilon \laplace^{2}
    \right)^{-1}
    \\
    & \left(
    u_{k}
    - \frac{\tau}{\varepsilon} W'(x)
    + \frac{\tau}{\varepsilon} \laplace W'(x)
    - \frac{\tau}{\varepsilon} W''(x) \laplace x
    + \frac{\tau}{\varepsilon^{3}} W'(x) W''(x)
    \right).
\label{eqn:fix_point_f_EE}
\end{aligned}
\end{eqnarray}
Moreover, the prepositional operator of~\eqref{eqn:fix_point_f_EE}
\begin{eqnarray}
\begin{aligned}
    \rho_{\mathscr{E}}(\xi)
    =
    \left(
    I_{D}
    - \varepsilon \tau \laplace
    + \tau \varepsilon \laplace^{2}
    \right)^{-1}
\end{aligned}
\end{eqnarray}
is able to implement via fast Fourier transform (FFT) and its associated symbol of a differential operator
\begin{eqnarray}
\begin{aligned}
    \rho_{\mathscr{E}}(\xi)
   =
    \left(
    I_{D}
    - \varepsilon \tau \laplace
    + \tau \varepsilon \laplace^{2}
    \right)^{-1}  =
    \frac
    {1}
    {1 + 4 \tau \varepsilon \pi^{2} |\xi|^{2} + 16\tau \varepsilon \pi^{4} |\xi|^{4}}.
\end{aligned}
\end{eqnarray}

Overall, each iteration of the Projected Gradient Descent Method (PGDM) is summarised in the following~\Cref{alg:Projected_Gradient_Descent_Method_GDM} for estimating the numerical solution of the Euler-Elastica-based formulation where another two options for perimeter-based formulation and Willmore-based formulation can be referred to~\cite{01_Volume_Reconstruction_from_Slices}.

\begin{algorithm}[htbp]
\caption{Projected Gradient Descent Method (PGDM)}
\label{alg:Projected_Gradient_Descent_Method_GDM}
    \KwIn{
    Initial set $E_{0}$;
    Parameters $\tau, \varepsilon$.
    }

    \KwOut{Numerical solution $u_{k+1}$.
    }
    Initial input:
    $u_{0} = q
            \left(
            \frac{\mathpzc{d}(\xi, E_{0})}{\varepsilon}
            \right)$
            \Comment*[r]{\eqref{equ:phase-field_function}}
    Interior restriction:
    $u^{in}_{E_{0}} = q
            \left(
            \frac{\mathpzc{d}(\xi, \Omega^{in})}{\varepsilon}
            \right)$
            \Comment*[r]{\eqref{equ:Uhinex}}
    Exterior restriction:
    $u^{ex}_{E_{0}} = 1 - q
            \left(
            \frac{\mathpzc{d}(\xi, \Omega^{ex})}{\varepsilon}
            \right)$
            \Comment*[r]{\eqref{equ:Uhinex}}
    \For{$k = 0, 1, \ldots$}
        {
        $u_{k+\frac{1}{2}}
        =
        \max(\min(u_{k}, u^{ex}_{E_{0}}), u^{in}_{E_{0}})$
        \Comment*[r]{\eqref{eqn:linear_obstacle_restriction_apprximate}}

        $
        \begin{aligned}
        u_{k+1}
        =
        & \left(
        I_{D}
        - \varepsilon \tau \laplace
        + \tau \varepsilon \laplace^{2}
        \right)^{-1}
        \\
        & \left(
        u_{k+\frac{1}{2}}
        -
        \frac{\tau}{\varepsilon} W'(u_{k+\frac{1}{2}})
        +
        \frac{\tau}{\varepsilon} \laplace W'(u_{k+\frac{1}{2}})
        \right.
        \\
        &
        \left.
        -
        \frac{\tau}{\varepsilon} W''(u_{k+\frac{1}{2}}) \laplace u_{k+\frac{1}{2}}
        +
        \frac{\tau}{\varepsilon^{3}} W'(u_{k+\frac{1}{2}}) W''(u_{k+\frac{1}{2}})
        \right)
        \end{aligned}
        $
        \Comment*[r]{\eqref{eqn:EE_iterative}}

        \textbf{End} till some stopping criteria are met.
        }
\end{algorithm}

\subsection{Numerical algorithm II}
\label{sec:numappro2}

In recent years, there has been a lot of progress in developing fast alternating direction method of multipliers (ADMM) for various applications~\cite{ADMM_Tai11, ADMM_Kang2015}.
Here we have extended the method to solve our model~\eqref{new_model}.

Expecting new variables $\mathbf{w} = \nabla u$ to have a faster numerical approximation, the formulation is converted to
\begin{eqnarray}
\begin{aligned}
&
\arg \min_{u, \mathbf{w}}
\int_{\Omega}
    \left(
    \frac{\varepsilon}{2} |\mathbf{w}|^{2}
    +
    \frac{1}{\varepsilon} W(u)
    \right)
    \, \id \Omega
    +
    \frac{1}{2\varepsilon}
    \int_{\Omega}
    \left(
    \varepsilon \operatorname{div} \mathbf{w}
    -
    \frac{1}{\varepsilon} W'(u)
    \right)^{2}
    \, \id \Omega
    \\
    &
    \mbox{ s.t. }
    \mathbf{w} = \nabla u
    \quad \mbox{and} \quad
    u^{in}_{E_{0}}
    \le
    u
    \le
    u^{ex}_{E_{0}}.
\end{aligned}
\end{eqnarray}
Then, the augmented Lagrangian functional for the above is expressed as
\begin{eqnarray}
\begin{aligned}
\mathcal{L}^{\rho}
(u, \mathbf{w}; \boldsymbol{\lambda})
&
=
\int_{\Omega}
    \left[
    \left(
    \frac{\varepsilon}{2} |\mathbf{w}|^{2}
    +
    \frac{1}{\varepsilon} W(u)
    \right)
    +
    \frac{1}{2\varepsilon}
    \left(
    \varepsilon \operatorname{div} \mathbf{w}
    -
    \frac{1}{\varepsilon} W'(u)
    \right)^{2}
    \right.
    \\
    &
    \qquad \quad
    +
    \left \langle
    \boldsymbol{\lambda}, \nabla u - \mathbf{w}
    \right \rangle
    +
    \left.
    \frac{\rho}{2}
    \left|
    \nabla u - \mathbf{w}
    \right|^{2}
    \right]
\, \id \Omega
\\
&
=
\int_{\Omega}
    \left[
    \left(
    \frac{\varepsilon}{2} |\mathbf{w}|^{2}
    +
    \frac{1}{\varepsilon} W(u)
    \right)
    +
    \frac{1}{2\varepsilon}
    \left(
    \varepsilon \operatorname{div} \mathbf{w}
    -
    \frac{1}{\varepsilon} W'(u)
    \right)^{2}
    \right.
    \\
    &
    \qquad \quad
    +
    \left.
    \frac{\rho}{2}
    \left|
    \nabla u - \mathbf{w}
    + \rho^{-1} \boldsymbol{\lambda}
    \right|^{2}
    -
    \frac{\boldsymbol{\lambda}^{2}}{2 \rho}
    \right]
\, \id \Omega
\end{aligned}
\end{eqnarray}
with the penalty parameter $\rho > 0$ and the Lagrange multiplier $\boldsymbol{\lambda}$.
Therefore, the problem in the ADMM can be considered as solving two subproblems and updating one multiplier:
\begin{eqnarray}
\label{eqn:ADMMpro}
\left\{
\begin{array}{l}
u_{k+1}
=
\arg \min \limits_{u}
\int_{\Omega}
\left[
\frac{1}{\varepsilon} W(u)
-
\frac{1}{\varepsilon} (\operatorname{div} \mathbf{w}_{k})  W'(u)
+
\frac{1}{2\varepsilon^{3}}
(W'(u))^{2}
\right.
\\
\left.
\qquad \qquad \qquad \qquad
+
\frac{\rho}{2}
\left|
\nabla u - \mathbf{w}_{k}
+ \rho^{-1} \boldsymbol{\lambda}_{k}
\right|^{2}
\right]
\, \id \Omega
\\
\mathbf{w}_{k+1}
=
\arg \min \limits_{\mathbf{w}}
\int_{\Omega}
\left[
\frac{\varepsilon}{2} |\mathbf{w}|^{2}
+
\frac{\varepsilon}{2}
(\operatorname{div} \mathbf{w})^{2}
-
\frac{1}{\varepsilon} (\operatorname{div} \mathbf{w})  W'(u_{k+1})
\right.
\\
\left.
\qquad \qquad \qquad \qquad
+
\frac{\rho}{2}
\left|
\mathbf{w}
- \nabla u_{k+1}
- \rho^{-1} \boldsymbol{\lambda}_{k}
\right|^{2}
\right]
\, \id \Omega
\\
\boldsymbol{\lambda}_{k+1}
=
\boldsymbol{\lambda}_{k}
+
\rho
\left(
\nabla u_{k+1}-\mathbf{w}_{k+1}
\right)
\end{array}
\right..
\end{eqnarray}

For the $u$-subproblem of~\eqref{eqn:ADMMpro}, by G\^ateaux derivative, for $\forall \tau$, there exists
\begin{subequations}
 \label{equ:u_sub}
\begin{align}
    \delta \mathcal{E}_{u_{k+1}}
    =
    &
    \int_{\Omega}
    \left(
    \frac{1}{\varepsilon}
    W'(u) \tau
    -
    \frac{1}{\varepsilon}
    \operatorname{div} \mathbf{w}_{k}
    W''(u) \tau
    +
    \frac{1}{\varepsilon^{3}}
    W'(u) W''(u) \tau
    \right)
    \, \id \Omega \label{equ:u_sub_a} \\
    &
    +
    \rho
    \int_{\Omega}
    \left(
    \nabla u - \mathbf{w}_{k}
    + \rho^{-1} \boldsymbol{\lambda}_{k}
    \right)
    \nabla \tau
    \, \id \Omega.
    \label{equ:u_sub_b}
\end{align}
\end{subequations}
Next by Green's formulae, for~\eqref{equ:u_sub_b}, there exists
\begin{eqnarray}
 \label{equ:u_sub_b_green}
\begin{aligned}
    & \int_{\Omega}
    \left(
    \nabla u - \mathbf{w}_{k}
    + \rho^{-1} \boldsymbol{\lambda}_{k}
    \right)
    \nabla \tau
    \, \id \Omega \\
    = &
    \oint_{\Gamma}
    \left(
    \nabla u - \mathbf{w}_{k}
    + \rho^{-1} \boldsymbol{\lambda}_{k}
    \right) \tau \cdot \boldsymbol{n}
    \, \id \Gamma
    -
    \int_{\Omega}
    \nabla \cdot
    \left(
    \nabla u - \mathbf{w}_{k}
    + \rho^{-1} \boldsymbol{\lambda}_{k}
    \right)
    \tau
    \, \id \Omega.
\end{aligned}
\end{eqnarray}

Therefore, assembling~\eqref{equ:u_sub} and~\eqref{equ:u_sub_b_green},
\begin{eqnarray*}
\begin{aligned}
    \delta \mathcal{E}_{u_{k+1}}
    =
    &
    \int_{\Omega}
    \left(
    \frac{1}{\varepsilon}
    W'(u) \tau
    -
    \frac{1}{\varepsilon}
    \operatorname{div} \mathbf{w}_{k}
    W''(u) \tau
    +
    \frac{1}{\varepsilon^{3}}
    W'(u) W''(u) \tau
    \right)
    \, \id \Omega \\
    &
    +
    \rho
    \left(
    \oint_{\Gamma}
    \left(
    \nabla u - \mathbf{w}_{k}
    + \rho^{-1} \boldsymbol{\lambda}_{k}
    \right) \tau \cdot \boldsymbol{n}
    \, \id \Gamma
    \right. \\
    & \qquad
    \left.
    -
    \int_{\Omega}
    \nabla \cdot
    \left(
    \nabla u - \mathbf{w}_{k}
    + \rho^{-1} \boldsymbol{\lambda}_{k}
    \right)
    \tau
    \, \id \Omega
    \right)
    = 0
\end{aligned}
\end{eqnarray*}
is permitted by the following E-L equation
\begin{eqnarray}
\label{equ:EL_u_sub}
    \frac{1}{\varepsilon} W'(u)
    -
    \frac{1}{\varepsilon}
    \operatorname{div} \mathbf{w}_{k} W''(u)
    +
    \frac{1}{\varepsilon^{3}}
    W'(u) W''(u)
    -
    \rho
    \nabla \cdot
    \left(
    \nabla u - \mathbf{w}_{k}
    + \rho^{-1} \boldsymbol{\lambda}_{k}
    \right)
    = 0.
\end{eqnarray}

Analogously, for the $\mathbf{w}$-subproblem of~\eqref{eqn:ADMMpro}, by G\^ateaux derivative, there exists
\begin{eqnarray}
 \label{equ:w_sub}
\begin{aligned}
    \delta \mathcal{E}_{\mathbf{w}_{k+1}}
    =
    &
    \int_{\Omega}
    \left(
    \varepsilon
    \mathbf{w}
    \tau
    +
    \varepsilon
    \operatorname{div} \mathbf{w}
    \cdot
    \nabla \tau
    -
    \frac{1}{\varepsilon}
    W'(u_{k+1})
    \laplace \mathbf{w}
    \cdot
    \tau
    \right.
    \\
    &
    \left.
    \qquad
    +
    \rho
    \left(
    \mathbf{w} - \nabla u_{k+1}
    - \rho^{-1} \boldsymbol{\lambda}_{k}
    \right)
    \tau
    \right)
    \, \id \Omega
\end{aligned}
\end{eqnarray}
for $\forall \tau$.
Next by Green's formulae, for the second term of~\eqref{equ:w_sub}, there exists
\begin{eqnarray}
 \label{equ:w_sub_green}
\begin{aligned}
    \int_{\Omega}
    \operatorname{div} \mathbf{w}
    \cdot
    \nabla \tau
    \, \id \Omega
    =
    \oint_{\Gamma}
    \operatorname{div} \mathbf{w}
    \cdot
    \tau
    \cdot
    \boldsymbol{n}
    \, \id \Gamma
    -
    \int_{\Omega}
    \laplace \mathbf{w}
    \cdot
    \tau
    \, \id \Omega.
\end{aligned}
\end{eqnarray}

Therefore, assembling~\eqref{equ:w_sub} and~\eqref{equ:w_sub_green},
\begin{eqnarray*}
\begin{aligned}
    \delta \mathcal{E}_{\mathbf{w}_{k+1}}
    =
    &
    \int_{\Omega}
    \left(
    \varepsilon
    \mathbf{w}
    -
    \frac{1}{\varepsilon}
    W'(u_{k+1})
    \laplace \mathbf{w}
    +
    \rho
    \left(
    \mathbf{w} - \nabla u_{k+1}
    - \rho^{-1} \boldsymbol{\lambda}_{k}
    \right)
    \right)
    \tau
    \, \id \Omega \\
    &
    +
    \varepsilon
    \left(
    \oint_{\Gamma}
    \operatorname{div} \mathbf{w}
    \cdot
    \tau
    \cdot
    \boldsymbol{n}
    \, \id \Gamma
    -
    \int_{\Omega}
    \laplace \mathbf{w}
    \cdot
    \tau
    \, \id \Omega
    \right)
\end{aligned}
\end{eqnarray*}
is permitted by the following E-L equation
\begin{eqnarray}
\label{equ:EL_w_sub}
    \varepsilon \mathbf{w}
    -
    \frac{1}{\varepsilon}
    W'(u_{k+1})
    \laplace \mathbf{w}
    +
    \rho
    \left(
    \mathbf{w}
    -
    \nabla u_{k+1}
    -
    \rho^{-1} \boldsymbol{\lambda}_{k}
    \right)
    -
    \varepsilon
    \laplace \mathbf{w}
    = 0.
\end{eqnarray}

Before proceeding with the numerical solutions, we first still enforce the linear obstacle restriction~\eqref{equ:restrictions_equivalence} as~\eqref{eqn:linear_obstacle_restriction_apprximate} in~\Cref{alg:Projected_Gradient_Descent_Method_GDM} of~\Cref{sec:numappro}, i.e.
\begin{eqnarray*}
    u^{in}_{E_{0}}
    \leq
    u
    \leq
    u^{ex}_{E_{0}}
    \quad
    \leftrightsquigarrow
    \quad
    \max(\min(u, u^{ex}_{E_{0}}), u^{in}_{E_{0}}).
\end{eqnarray*}
Then, to progress the numerical solution of~\eqref{eqn:ADMMpro}, the minimising solution of $u$-subproblem is given by
\begin{eqnarray}
\begin{aligned}
    \frac
    {u_{k+1}-u_{k}}
    {\tau}
    =
    &
    \rho
    \nabla \cdot
    \left(
    \nabla u_{k+1} - \mathbf{w}_{k}
    + \rho^{-1} \boldsymbol{\lambda}_{k}
    \right)
    -
    \frac{1}{\varepsilon} W'(u_{k+1}) \\
    &
    +
    \frac{1}{\varepsilon}
    \operatorname{div} \mathbf{w}_{k} W''(u_{k+1})
    +
    \frac{1}{\varepsilon^{3}}
    W'(u_{k+1}) W''(u_{k+1}),
\end{aligned}
\end{eqnarray}
which leads to the update
\begin{eqnarray}
\label{eqn:EE_ADMM_u}
\begin{aligned}
    u_{k+1}
    =
    (I - \tau \rho \laplace)^{-1}
    &
    \left[
    u_{k}
    +
    \tau
    \left(
    - \rho \nabla \cdot \mathbf{w}_{k}
    + \nabla \cdot \boldsymbol{\lambda}_{k}
    -
    \frac{1}{\varepsilon} W'(u_{k})
    \right.
    \right.
    \\
    &
    \left.
    \left.
    +
    \frac{1}{\varepsilon}
    \operatorname{div} \mathbf{w}_{k} W''(u_{k})
    +
    \frac{1}{\varepsilon^{3}}
    W'(u_{k}) W''(u_{k})
    \right)
    \right].
\end{aligned}
\end{eqnarray}
Next, for the minimising solution of $\mathbf{w}$-subproblem, it is discretised by
\begin{eqnarray}
\begin{aligned}
    \frac
    {\mathbf{w}_{k+1}-\mathbf{w}_{k}}
    {\tau}
    =
    &
    \frac{1}{\varepsilon}
    W'(u_{k+1})
    \laplace \mathbf{w}_{k+1}
    -
    \varepsilon \mathbf{w}_{k+1} \\
    &
    +
    \rho
    \left(
    \mathbf{w}_{k+1}
    -
    \nabla u_{k+1}
    -
    \rho^{-1} \boldsymbol{\lambda}_{k}
    \right)
    +
    \varepsilon
    \laplace \mathbf{w}_{k+1},
\end{aligned}
\end{eqnarray}
which issues in the update
\begin{eqnarray}
\label{eqn:EE_ADMM_w}
\begin{aligned}
    \mathbf{w}_{k+1}
    =
    (I + \varepsilon \tau - \varepsilon \tau \laplace)^{-1}
    &
    \left[
    \mathbf{w}_{k}
    +
    \frac{\tau}{\varepsilon}
    W'(u_{k+1})
    \laplace \mathbf{w}_{k}
    \right.
    \\
    &
    \left.
    +
    \tau
    \rho
    \left(
    \mathbf{w}_{k}
    -
    \nabla u_{k+1}
    -
    \rho^{-1} \boldsymbol{\lambda}_{k}
    \right)
    \right].
\end{aligned}
\end{eqnarray}
Lastly, recall that the multiplier $\boldsymbol{\lambda}$ will be updated by
\begin{eqnarray}
\label{eqn:EE_ADMM_lambda}
\begin{aligned}
    \boldsymbol{\lambda}_{k+1}
    =
    \boldsymbol{\lambda}_{k}
    +
    \rho
    \left(
    \nabla u_{k+1}-\mathbf{w}_{k+1}
    \right).
\end{aligned}
\end{eqnarray}
This completes one step of the ADMM method, with the final algorithm shown in~\Cref{alg:Alternating_Direction_Method_of_Multipliers}.

\begin{algorithm}[htbp]
\caption{Alternating Direction Method of Multipliers (ADMM)}
\label{alg:Alternating_Direction_Method_of_Multipliers}
    \KwIn{Initial set $E_{0}$;
    Parameters $\tau, \varepsilon, \rho$.
    }

    \KwOut{Numerical solution $u_{k+1}$.
    }
    Initial input:
    $u_{0} = q
            \left(
            \frac{\mathpzc{d}(\xi, E_{0})}{\varepsilon}
            \right)$
            \Comment*[r]{\eqref{equ:phase-field_function}}
    \hspace{1.95 cm} $\mathbf{w}_{0} = \grad u_{0}$ \;
    \hspace{2 cm} $\boldsymbol{\lambda}_{0} = \mathbf{w}_{0}$ \;
    Interior restriction:
    $u^{in}_{E_{0}} = q
            \left(
            \frac{\mathpzc{d}(\xi, \Omega^{in})}{\varepsilon}
            \right)$
            \Comment*[r]{\eqref{equ:Uhinex}}
    Exterior restriction:
    $u^{ex}_{E_{0}} = 1 - q
            \left(
            \frac{\mathpzc{d}(\xi, \Omega^{ex})}{\varepsilon}
            \right)$
            \Comment*[r]{\eqref{equ:Uhinex}}
    \For{$k = 0, 1, \ldots$}
        {
        $u_{k+\frac{1}{2}}
        =
        \max(\min(u_{k}, u^{ex}_{E_{0}}), u^{in}_{E_{0}})$
        \Comment*[r]{\eqref{eqn:linear_obstacle_restriction_apprximate}}

        $
        \begin{aligned}
            u_{k+1}
            =
            &
            (I - \tau \rho \laplace)^{-1}
            \left[
            u_{k+\frac{1}{2}}
            +
            \tau
            \left(
            - \rho \nabla \cdot \mathbf{w}_{k}
            + \nabla \cdot \boldsymbol{\lambda}_{k}
            \right.
            \right.
            \\
            &
            \left.
            \left.
            -
            \frac{1}{\varepsilon} W'(u_{k+\frac{1}{2}})
            +
            \frac{1}{\varepsilon}
            \operatorname{div} \mathbf{w}_{k} W''(u_{k+\frac{1}{2}})
            \right.
            \right.
            \\
            &
            \left.
            \left.
            +
            \frac{1}{\varepsilon^{3}}
            W'(u_{k+\frac{1}{2}}) W''(u_{k+\frac{1}{2}})
            \right)
            \right]
        \end{aligned}
        $
        \Comment*[r]{\eqref{eqn:EE_ADMM_u}}

        $
        \begin{aligned}
            \mathbf{w}_{k+1}
            =
            (I + \varepsilon \tau - \varepsilon \tau \laplace)^{-1}
            &
            \left[
            \mathbf{w}_{k}
            +
            \frac{\tau}{\varepsilon}
            W'(u_{k+1})
            \laplace \mathbf{w}_{k}
            \right.
            \\
            &
            \left.
            +
            \tau
            \rho
            \left(
            \mathbf{w}_{k}
            -
            \nabla u_{k+1}
            -
            \rho^{-1} \boldsymbol{\lambda}_{k}
            \right)
            \right]
        \end{aligned}
        $
        \Comment*[r]{\eqref{eqn:EE_ADMM_w}}

        $
        \boldsymbol{\lambda}_{k+1}
        =
        \boldsymbol{\lambda}_{k}
        +
        \rho
        \left(
        \nabla u_{k+1}-\mathbf{w}_{k+1}
        \right)
        $
        \Comment*[r]{\eqref{eqn:EE_ADMM_lambda}}

        \textbf{End} till some stopping criteria are met.
        }
\end{algorithm}

\section{Experimental results and quantitative comparisons with analysis}
\label{sec:Results}

\begin{figure}[H]
    \centering
    \begin{subfigure}[t]{0.32\textwidth}
    \includegraphics[width=\textwidth]{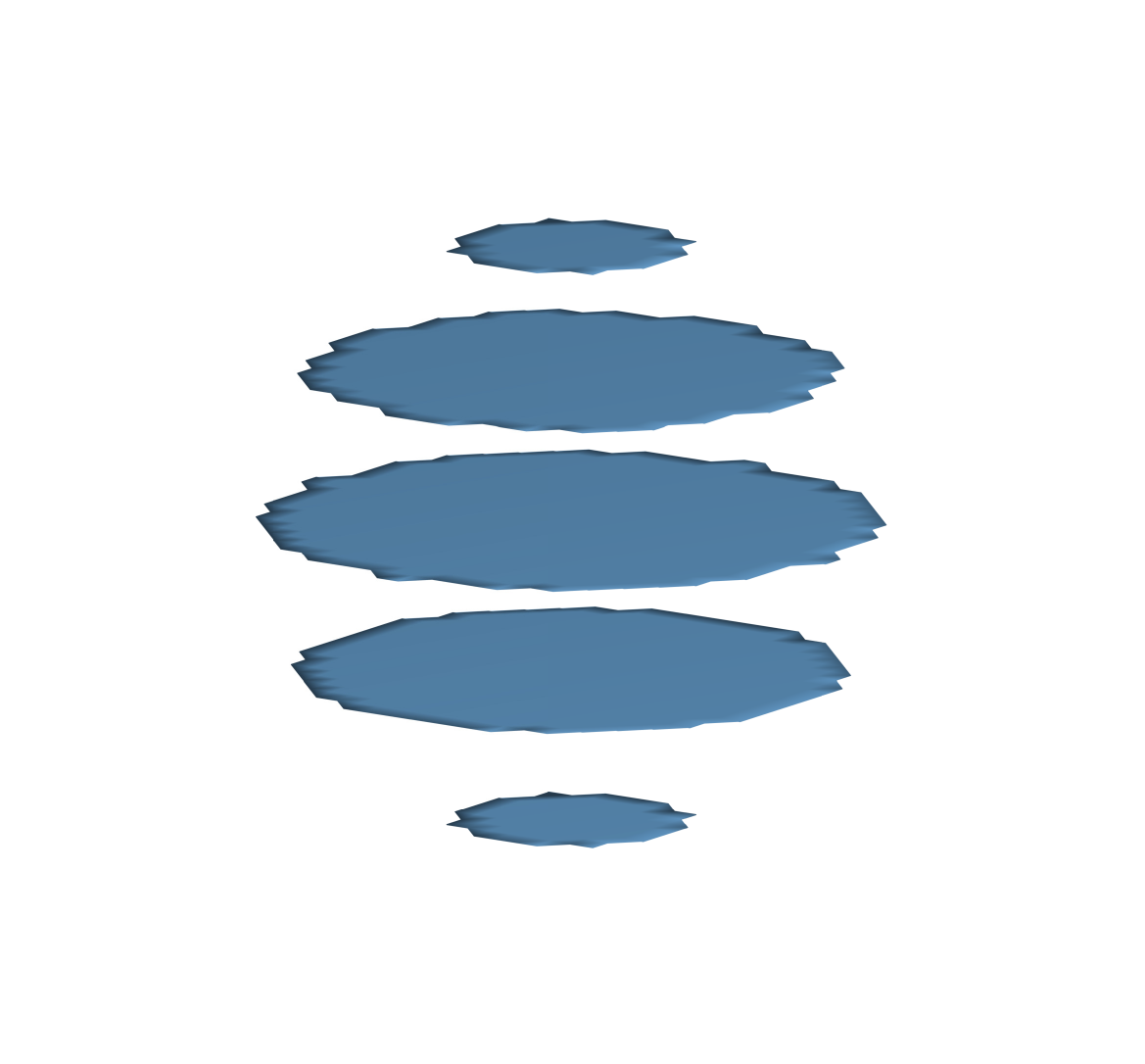}
    \caption{}
    \label{p20220215_S_N32S5_NoTITLE_GCMClim_PWEE_2}
    \end{subfigure}
\hspace{0.1mm}
    \begin{subfigure}[t]{0.32\textwidth}
    \includegraphics[width=\textwidth]{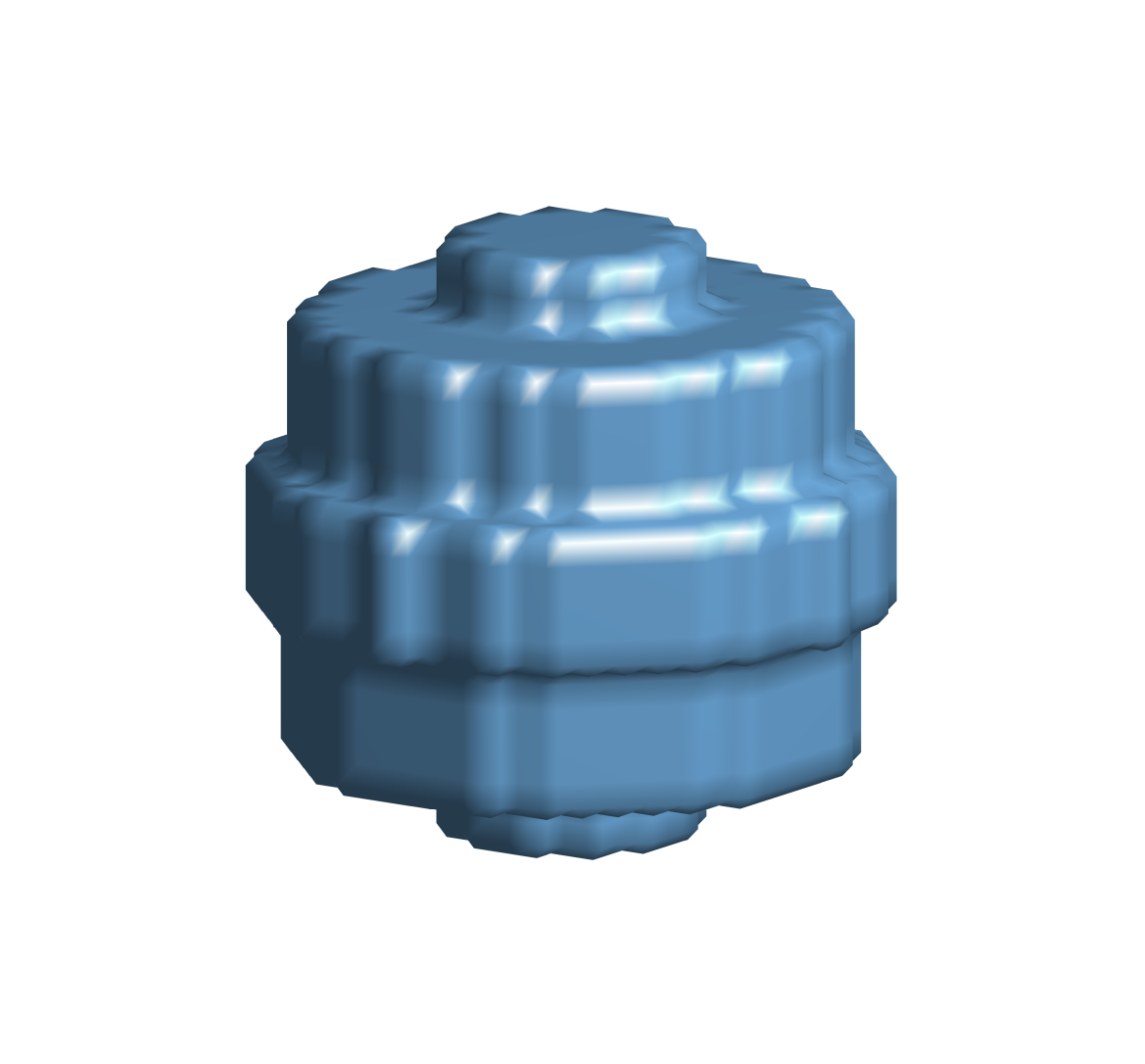}
    \caption{}
    \label{p20220215_S_N32S5_NoTITLE_GCMClim_PWEE_3}
    \end{subfigure}
\cprotect \caption{Visualisation of~\subref{p20220215_S_N32S5_NoTITLE_GCMClim_PWEE_3} the initial rough surface of {\bf Example 1} Sphere by duplicating the slices from~\subref{p20220215_S_N32S5_NoTITLE_GCMClim_PWEE_2} the given $5$ slices under the low resolution $N = 32$. }
\vspace{-.2cm}
\label{p20220215_S_N32S5_NoTITLE_GCMClim_PWEE_2-3}
\end{figure}

\begin{figure}[H]
    \centering
    \begin{subfigure}[t]{0.32\textwidth}
    \includegraphics[width=\textwidth]{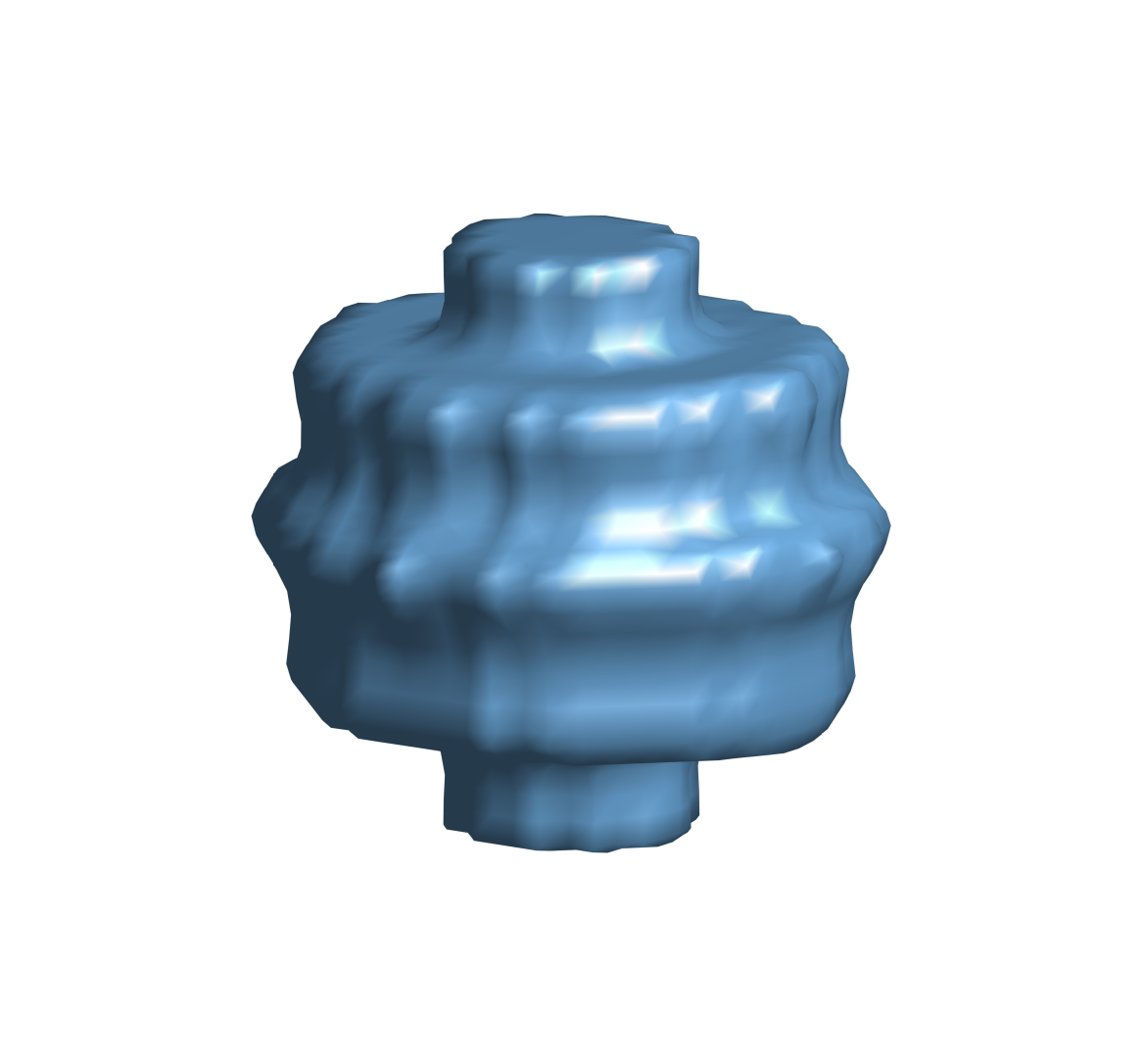}
    \vspace{-.8cm}
    \caption{}
    \label{p20220715_S_N32S5_CH_MS_AMLE_T_H_inpainting_CH_2}
    \end{subfigure}
\hspace{0.1mm}
    \begin{subfigure}[t]{0.32\textwidth}
    \includegraphics[width=\textwidth]{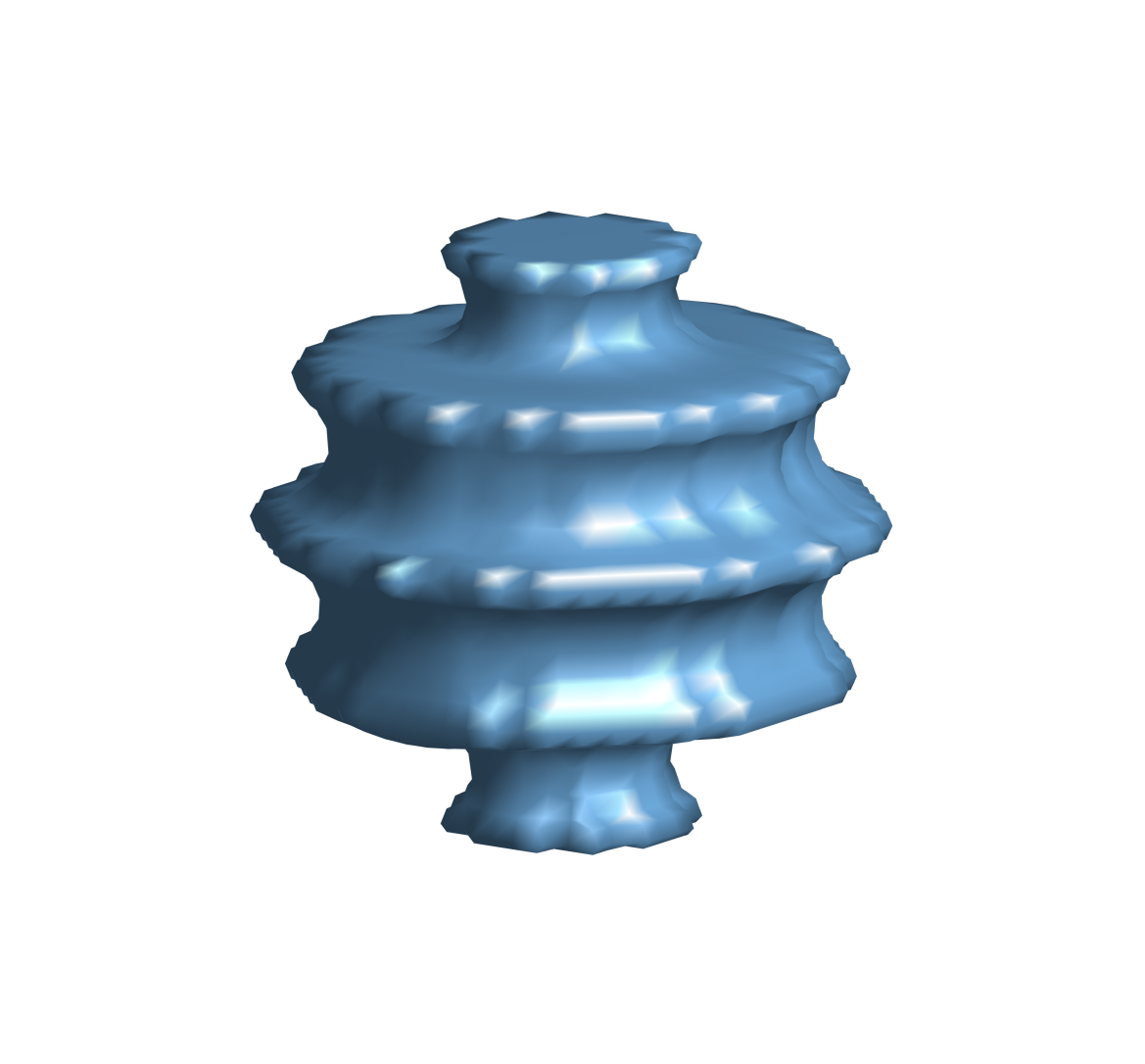}
    \vspace{-.8cm}
    \caption{}
    \label{p20220715_S_N32S5_CH_MS_AMLE_T_H_inpainting_MS_3}
    \end{subfigure}
\\
    \begin{subfigure}[t]{0.32\textwidth}
    \includegraphics[width=\textwidth]{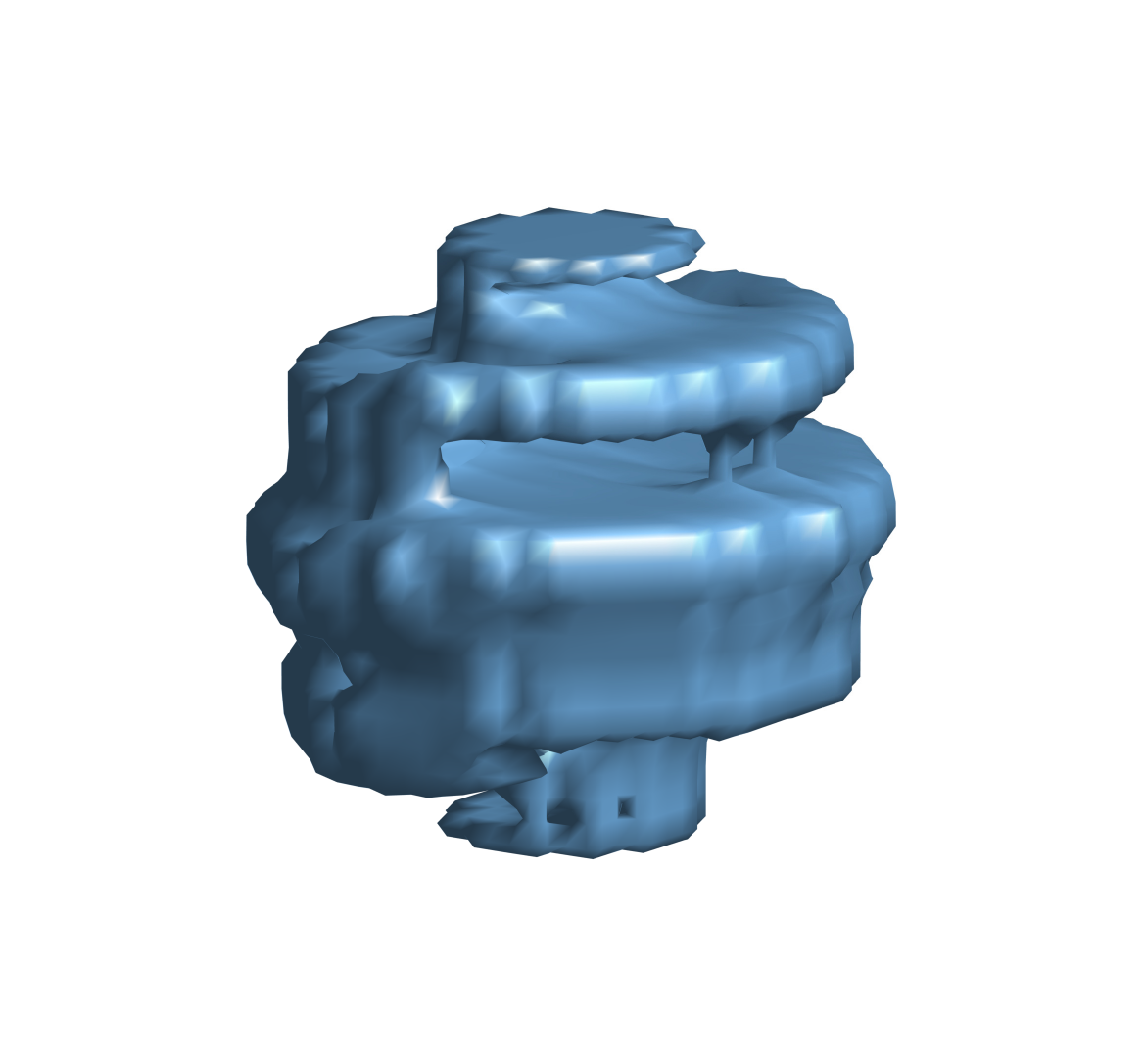}
    \vspace{-.8cm}
    \caption{}
    \label{p20220715_S_N32S5_CH_MS_AMLE_T_H_inpainting_T_5}
    \end{subfigure}
\hspace{0.1mm}
    \begin{subfigure}[t]{0.32\textwidth}
    \includegraphics[width=\textwidth]{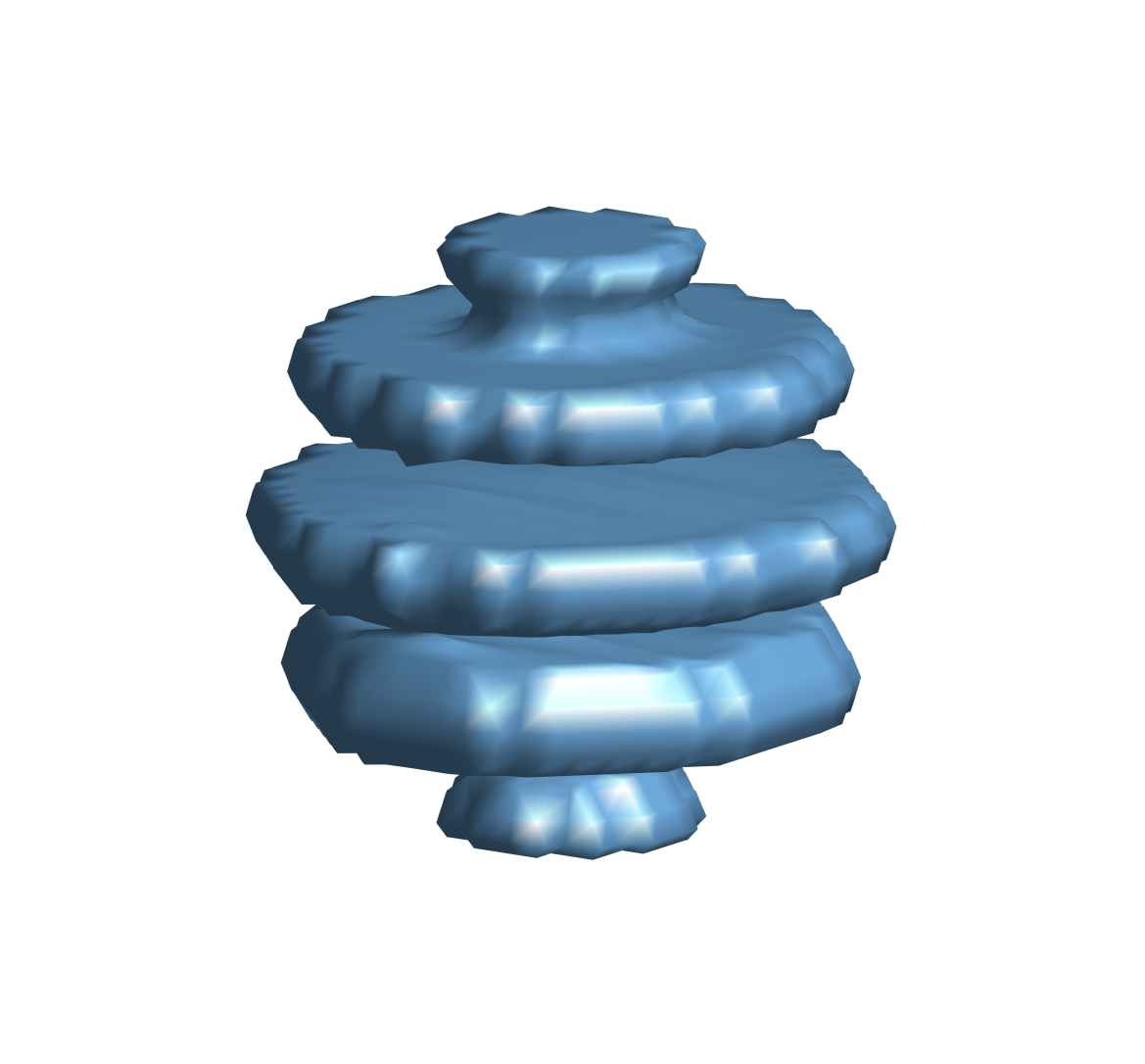}
    \vspace{-.8cm}
    \caption{}
    \label{p20220715_S_N32S5_CH_MS_AMLE_T_H_inpainting_AMLE_4}
    \end{subfigure}
\cprotect \caption{Inpainting results of {\bf Example 1} Sphere by four models:
~\subref{p20220715_S_N32S5_CH_MS_AMLE_T_H_inpainting_CH_2} the Cahn-Hilliard model;
~\subref{p20220715_S_N32S5_CH_MS_AMLE_T_H_inpainting_MS_3} the Mumford-Shah model;
~\subref{p20220715_S_N32S5_CH_MS_AMLE_T_H_inpainting_T_5} the transport model and~\subref{p20220715_S_N32S5_CH_MS_AMLE_T_H_inpainting_AMLE_4} the absolute minimising Lipschitz extensions from~\Cref{p20220215_S_N32S5_NoTITLE_GCMClim_PWEE_2} the given $5$ input slices under the low resolution $N = 32$. (Clearly, inpainting methods do not work well if given only a few slices. )}
\label{p20220715_S_N32S5_CH_MS_AMLE_T_H_inpainting_2-5}
\end{figure}

\begin{figure}[H]
    \centering
    \vspace{-.8cm}
    \begin{subfigure}[t]{0.32\textwidth}
    \includegraphics[width=\textwidth]{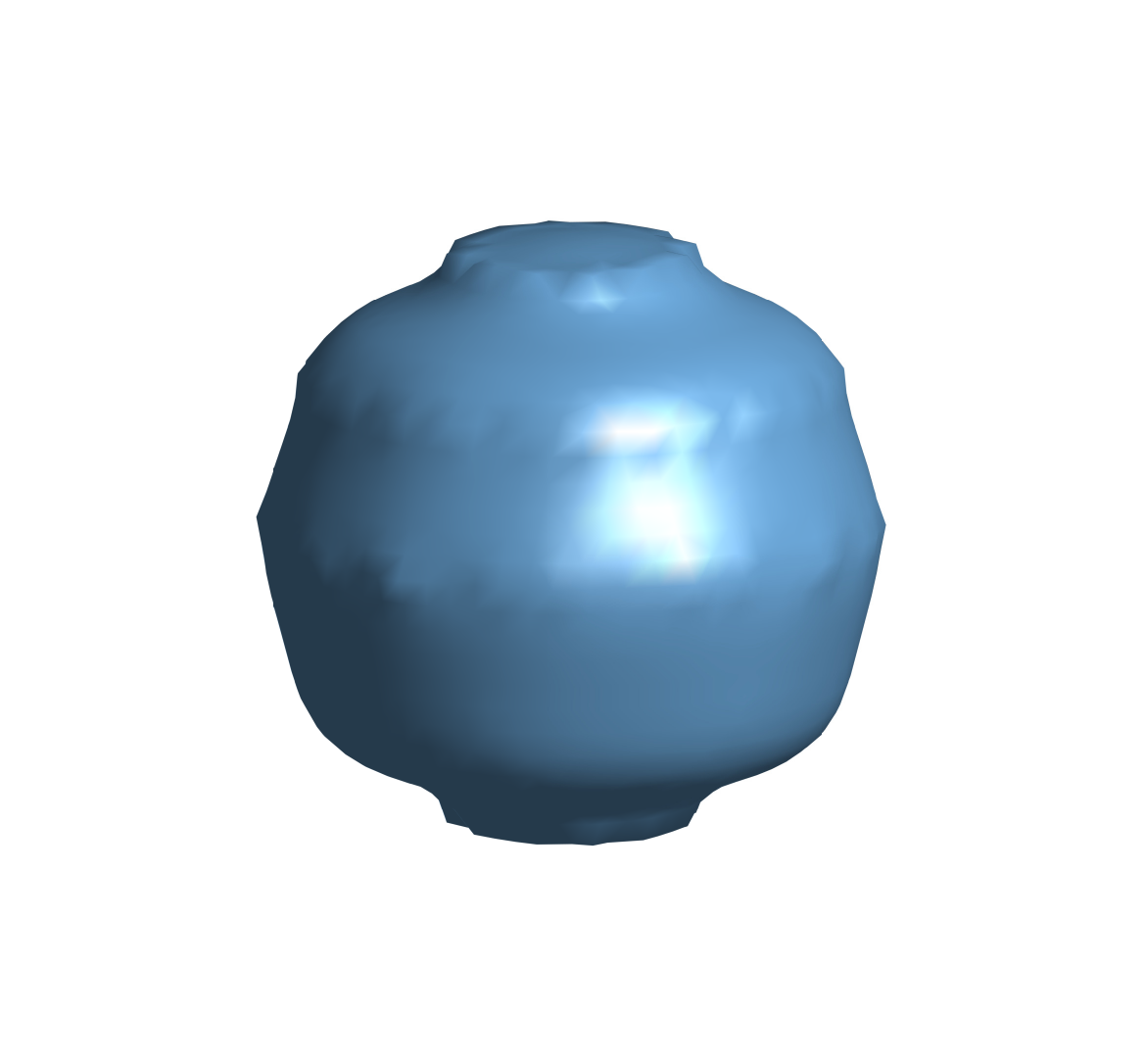}
    \vspace{-.8cm}
    \caption{}
    \label{p20220215_S_N32S5_NoTITLE_GCMClim_PWEE_4}
    \end{subfigure}
\hspace{0.1mm}
    \begin{subfigure}[t]{0.32\textwidth}
    \includegraphics[width=\textwidth]{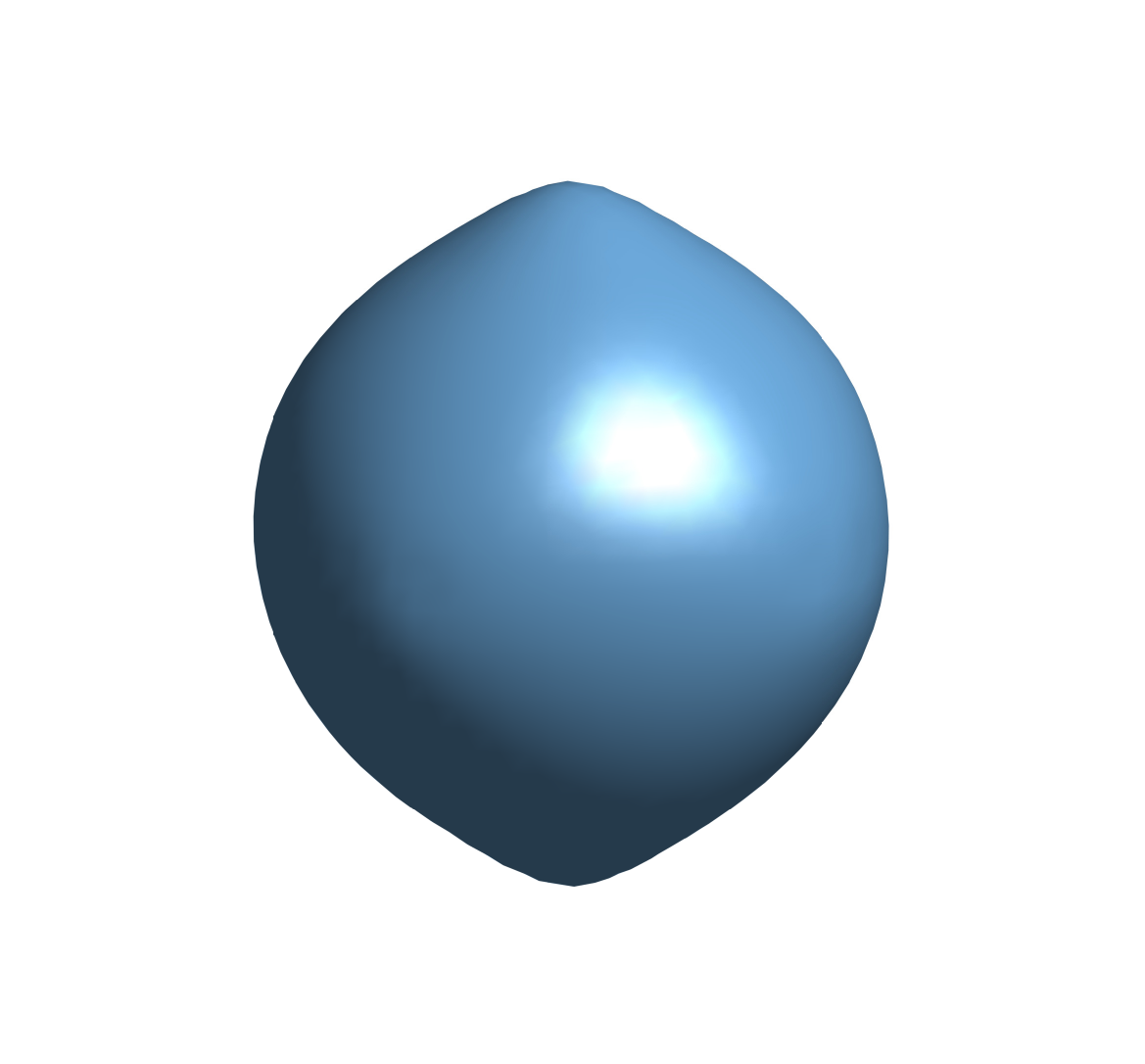}
    \vspace{-.8cm}
    \caption{}
    \label{p20220215_S_N32S5_NoTITLE_GCMClim_PWEE_5}
    \end{subfigure}
\hspace{0.1mm}
    \begin{subfigure}[t]{0.32\textwidth}
    \includegraphics[width=\textwidth]{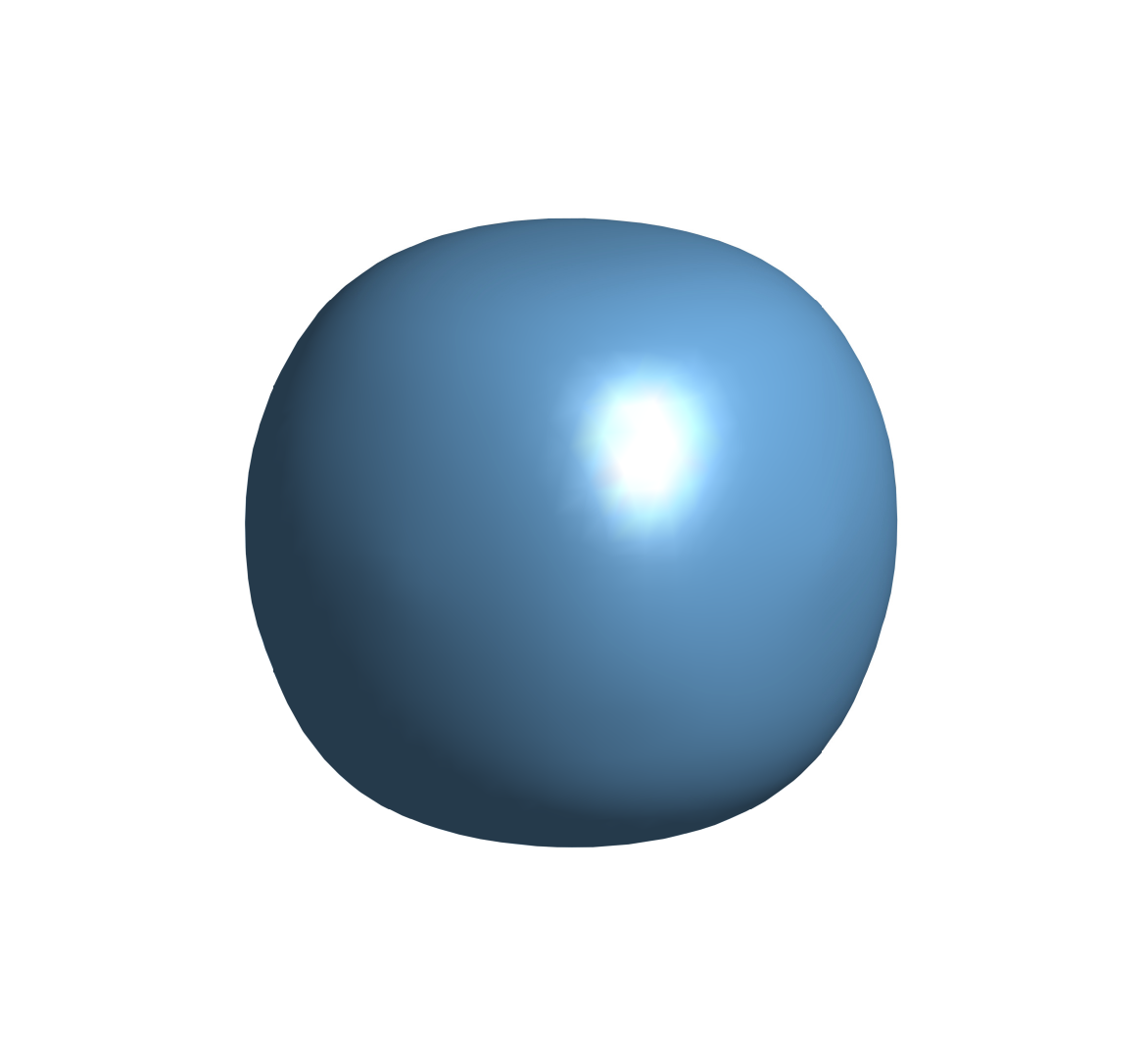}
    \vspace{-.8cm}
    \caption{}
    \label{p20220215_S_N32S5_NoTITLE_GCMClim_PWEE_6}
    \end{subfigure}
\cprotect \caption{Final reconstructed results of {\bf Example 1} Sphere by three formulations:~\subref{p20220215_S_N32S5_NoTITLE_GCMClim_PWEE_4} perimeter-based formulation;~\subref{p20220215_S_N32S5_NoTITLE_GCMClim_PWEE_5} Willmore-based formulation; and~\subref{p20220215_S_N32S5_NoTITLE_GCMClim_PWEE_6} Euler-Elastica-based formulation; from low-resolution inputs $N = 32$ with $5$ slices. (Visually, the result~\subref{p20220215_S_N32S5_NoTITLE_GCMClim_PWEE_6} by the new model is the best. )}
\label{p20220215_S_4_5_6}
\end{figure}

\begin{figure}[H]
    \centering
    \begin{subfigure}[t]{0.32\textwidth}
    \includegraphics[width=\textwidth]{00_Fig/DATE20220215_BC_N128S24_NoTITLE_GCMClim_PWEE/p20220215_BC_N128S24_NoTITLE_GCMClim_PWEE_2-eps-converted-to.pdf}
    \vspace{-1cm}
    \caption{}
    \label{p20220215_BC_N128S24_NoTITLE_GCMClim_PWEE_2}
    \end{subfigure}
\hspace{0.1mm}
    \begin{subfigure}[t]{0.32\textwidth}
    \includegraphics[width=\textwidth]{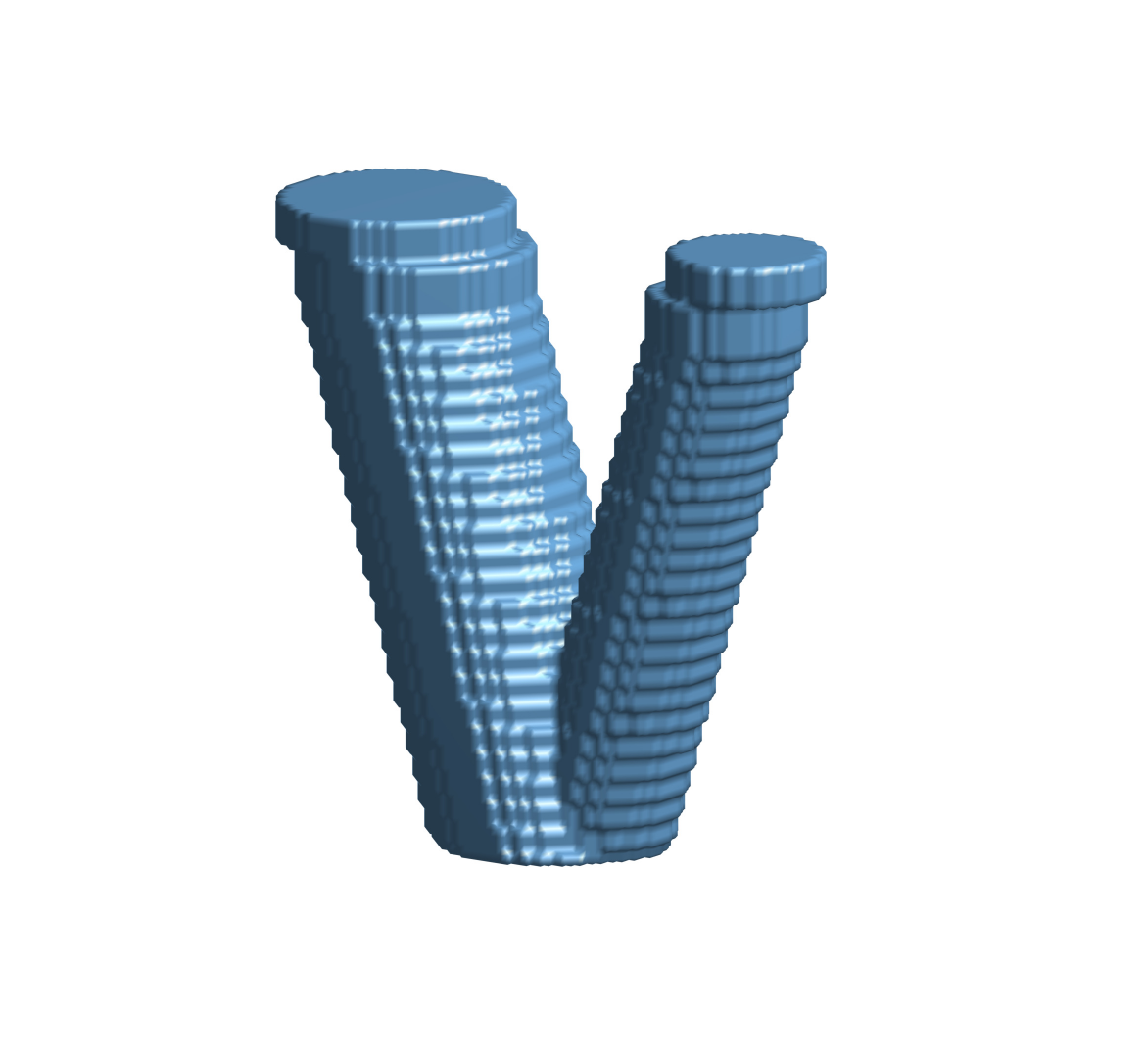}
    \vspace{-1cm}
    \caption{}
    \label{p20220215_BC_N128S24_NoTITLE_GCMClim_PWEE_3}
    \end{subfigure}
\vspace{-.3cm}
\cprotect \caption{Visualisation of~\subref{p20220215_BC_N128S24_NoTITLE_GCMClim_PWEE_3} the initial rough surface of {\bf Example 2} Branching Cylinders from~\subref{p20220215_BC_N128S24_NoTITLE_GCMClim_PWEE_2} the given $24$ slices under the low resolution $N = 128$. }
\label{p20220215_BC_2_3}
\end{figure}

\begin{figure}[H]
    \centering
    \begin{subfigure}[t]{0.32\textwidth}
    \includegraphics[width=\textwidth]{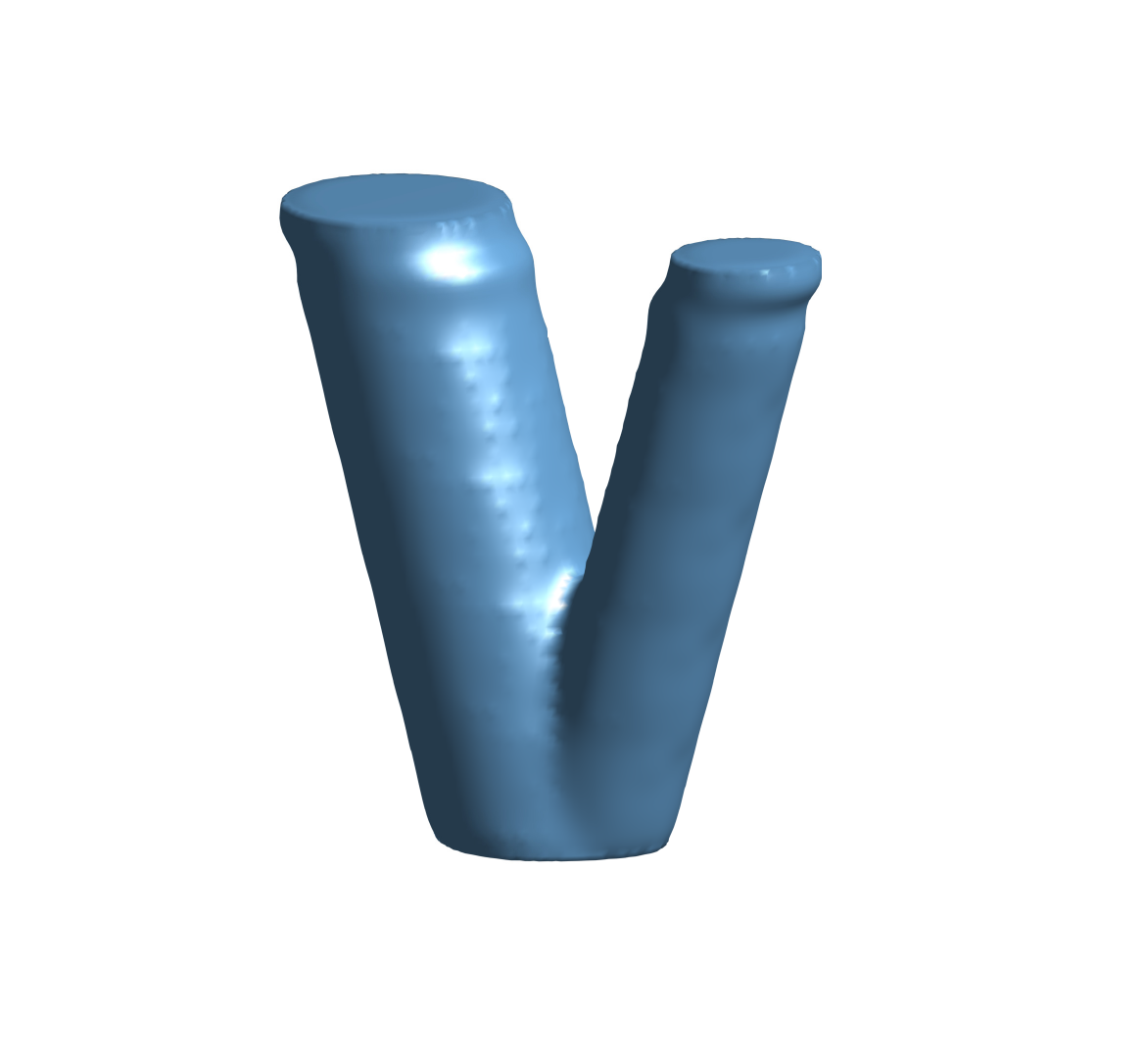}
    \vspace{-1cm}
    \caption{}
    \label{p20220215_BC_N128S24_NoTITLE_GCMClim_PWEE_4}
    \end{subfigure}
\hspace{0.1mm}
    \begin{subfigure}[t]{0.32\textwidth}
    \includegraphics[width=\textwidth]{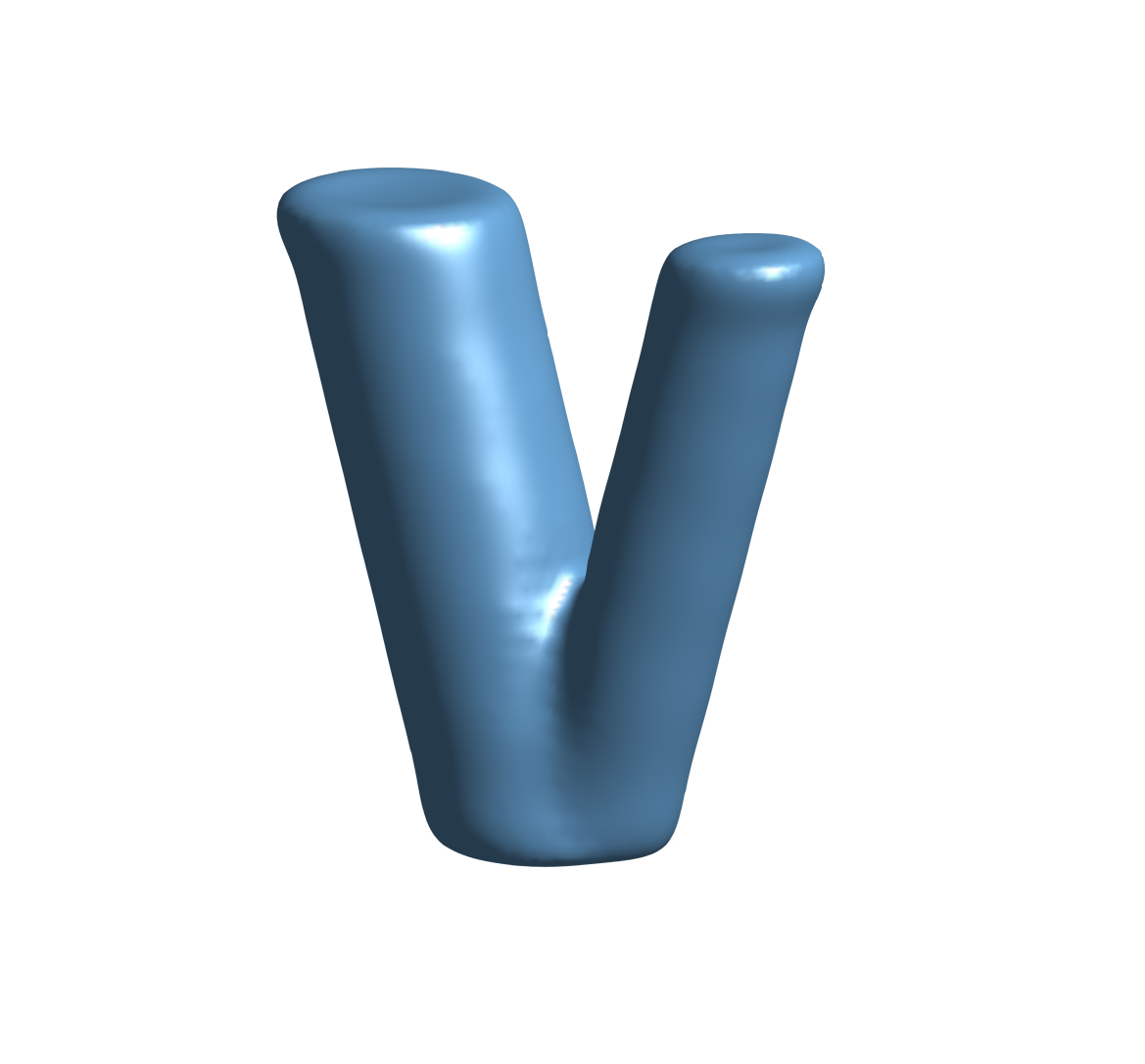}
    \vspace{-1cm}
    \caption{}
    \label{p20220215_BC_N128S24_NoTITLE_GCMClim_PWEE_5}
    \end{subfigure}
\hspace{0.1mm}
    \begin{subfigure}[t]{0.32\textwidth}
    \includegraphics[width=\textwidth]{00_Fig/DATE20220215_BC_N128S24_NoTITLE_GCMClim_PWEE/p20220215_BC_N128S24_NoTITLE_GCMClim_PWEE_6-eps-converted-to.pdf}
    \vspace{-1cm}
    \caption{}
    \label{p20220215_BC_N128S24_NoTITLE_GCMClim_PWEE_6}
    \end{subfigure}
\vspace{-.3cm}
\cprotect \caption{Final reconstructed results of {\bf Example 2} Branching Cylinders by three formulations:~\subref{p20220215_BC_N128S24_NoTITLE_GCMClim_PWEE_4} perimeter-based formulation;~\subref{p20220215_BC_N128S24_NoTITLE_GCMClim_PWEE_5} Willmore-based formulation; and~\subref{p20220215_BC_N128S24_NoTITLE_GCMClim_PWEE_6} Euler-Elastica-based formulation from low-resolution inputs $N = 128$ with $24$ slices. }
\label{p20220215_BC_4_5_6}
\end{figure}

\begin{figure}[htbp]
\centering
\input{00_Fig/0_Tikz/20220716_ViewRing}
\caption{Computing curvatures -- Illustration of the patch within the 2-ring neighbouring vertices $_{2}NV_{i}$ of the centre vertex $\mathbf{v}_{i}$ (blue) where the first ring neighbouring vertices are red, and the second ring neighbouring vertices are green. }
\label{20220716_ViewRing}
\end{figure}
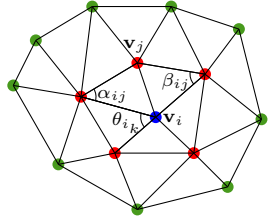

\begin{figure}[H]
    \centering
    \begin{subfigure}[t]{0.32\textwidth}
    \includegraphics[width=\textwidth]{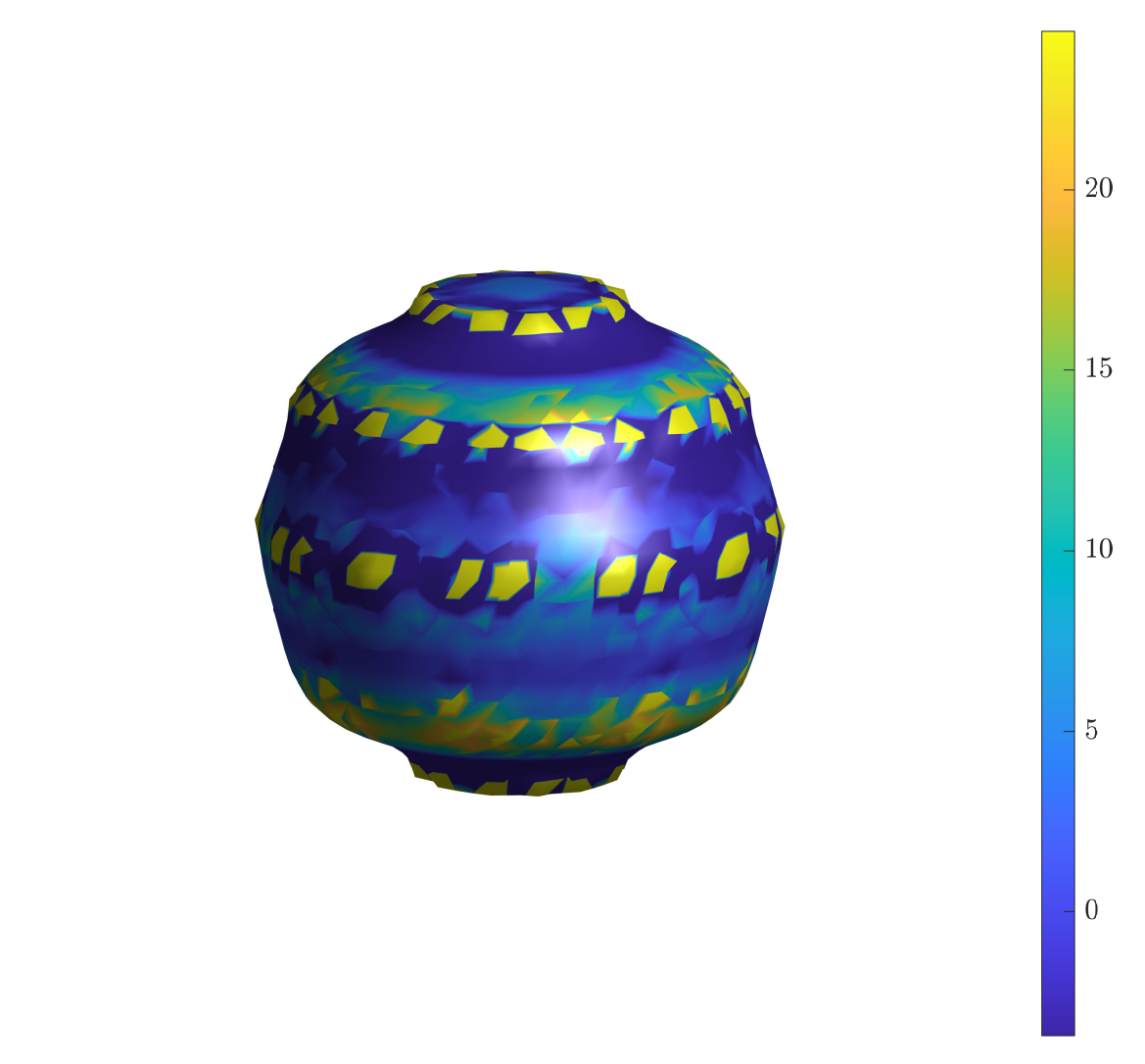}
    \label{p20220215_S_N32S5_NoTITLE_GCMClim_PWEE_9}
    \end{subfigure}
\hspace{0.1mm}
    \begin{subfigure}[t]{0.32\textwidth}
    \includegraphics[width=\textwidth]{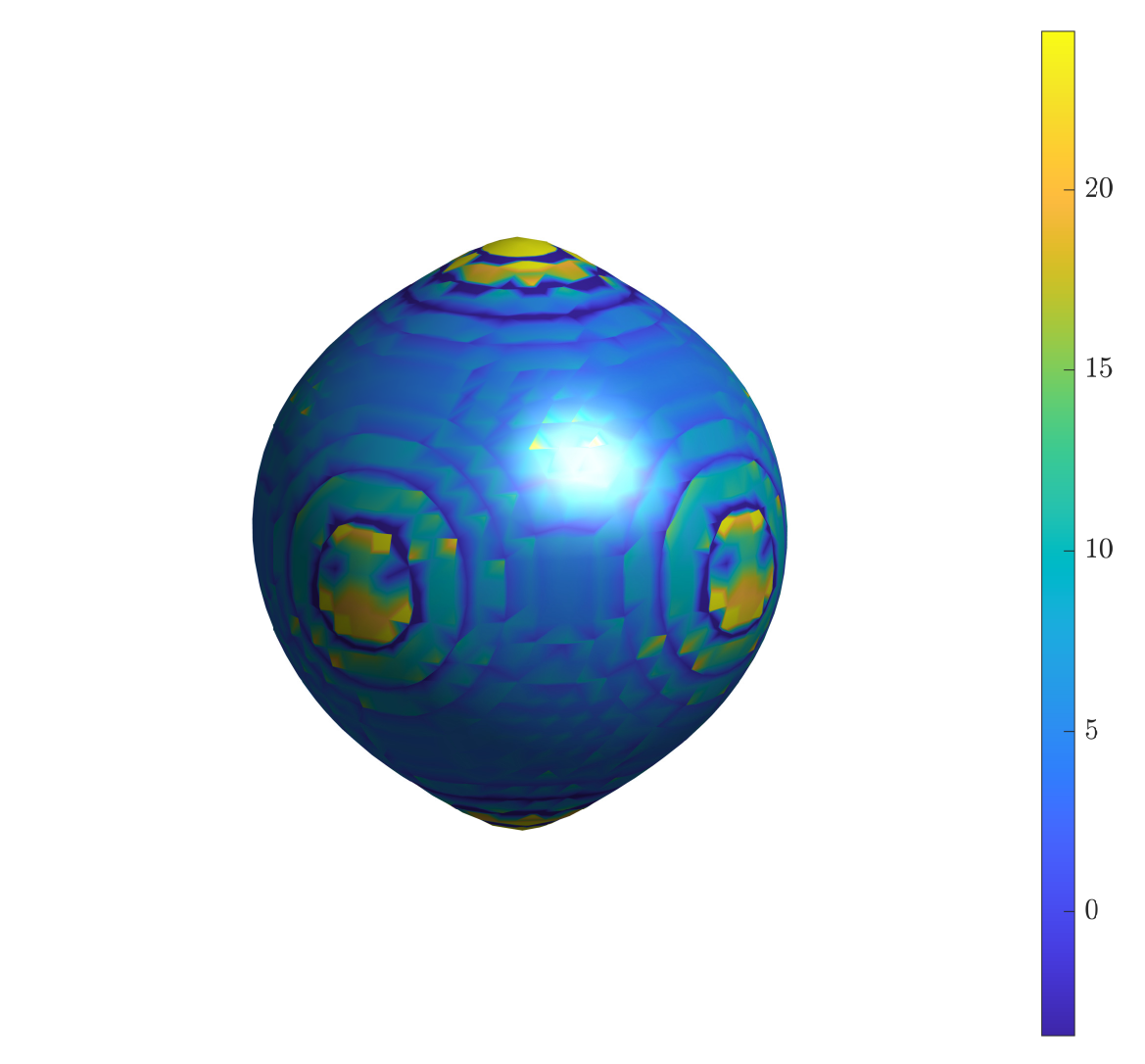}
    \label{p20220215_S_N32S5_NoTITLE_GCMClim_PWEE_11}
    \end{subfigure}
\hspace{0.1mm}
    \begin{subfigure}[t]{0.32\textwidth}
    \includegraphics[width=\textwidth]{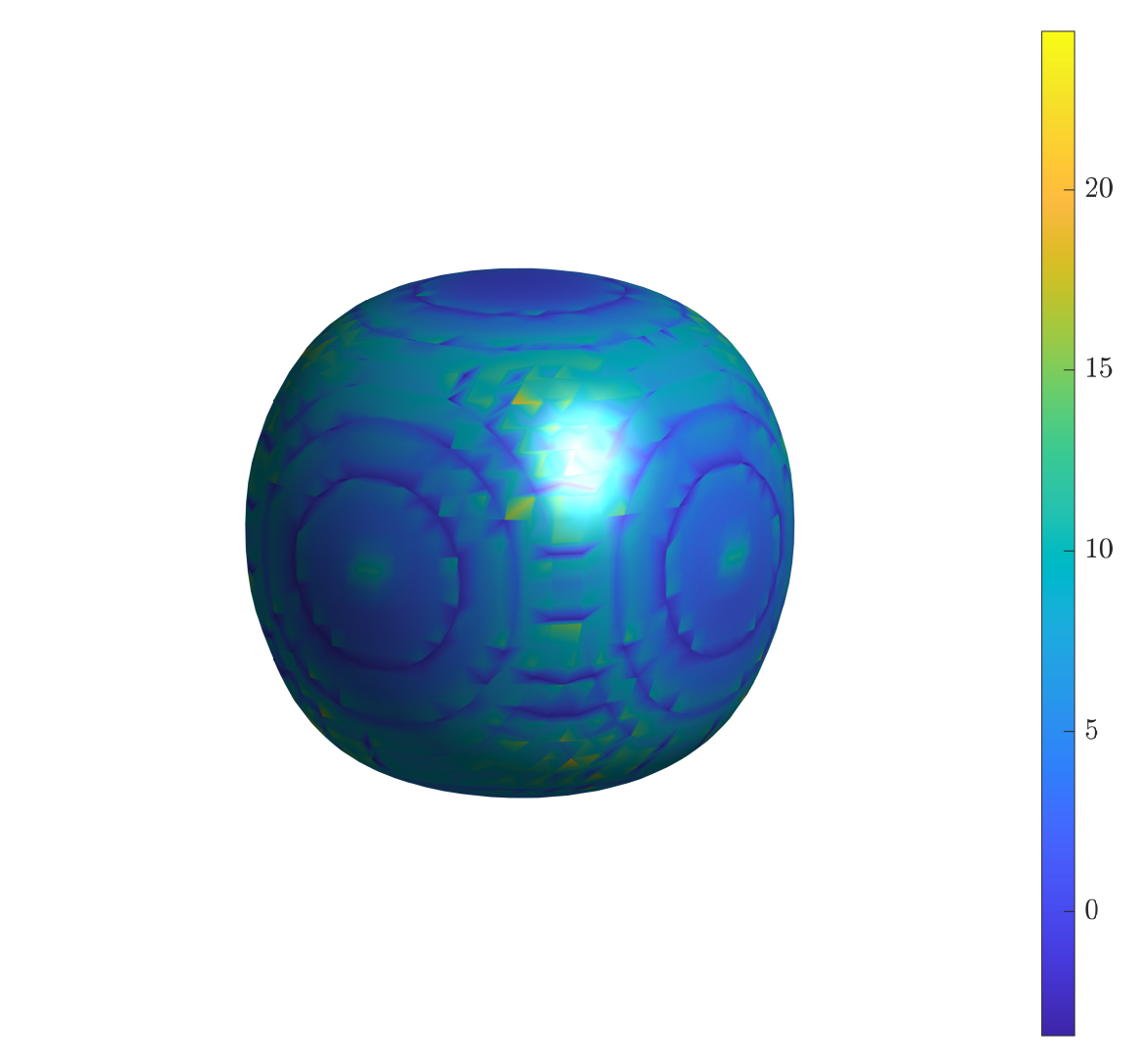}
    \label{p20220215_S_N32S5_NoTITLE_GCMClim_PWEE_13}
    \end{subfigure}
\\
    \begin{subfigure}[t]{0.32\textwidth}
    \includegraphics[width=\textwidth]{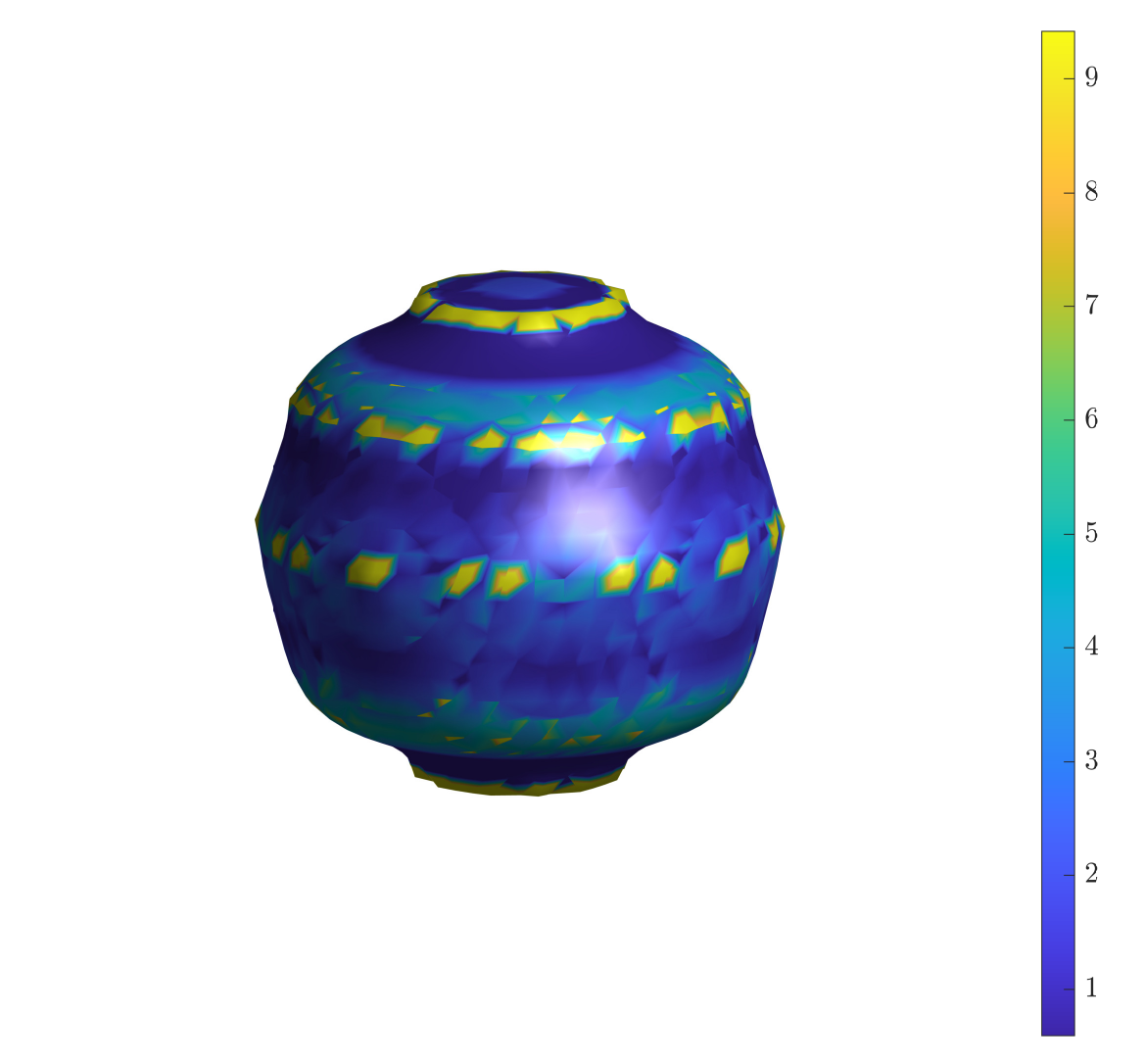}
    \label{p20220215_S_N32S5_NoTITLE_GCMClim_PWEE_10}
    \end{subfigure}
\hspace{0.1mm}
    \begin{subfigure}[t]{0.32\textwidth}
    \includegraphics[width=\textwidth]{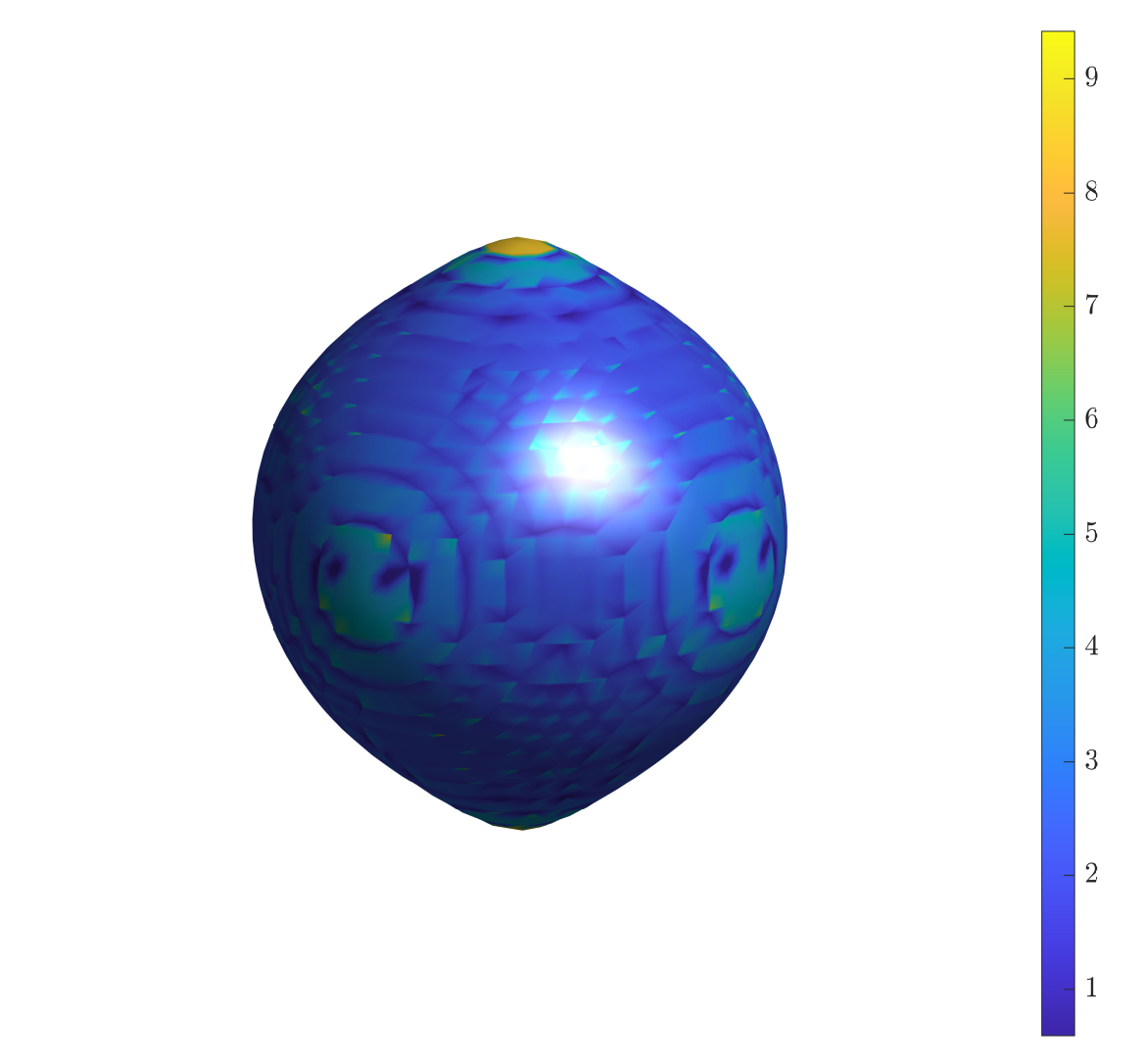}
    \label{p20220215_S_N32S5_NoTITLE_GCMClim_PWEE_12}
    \end{subfigure}
\hspace{0.1mm}
    \begin{subfigure}[t]{0.32\textwidth}
    \includegraphics[width=\textwidth]{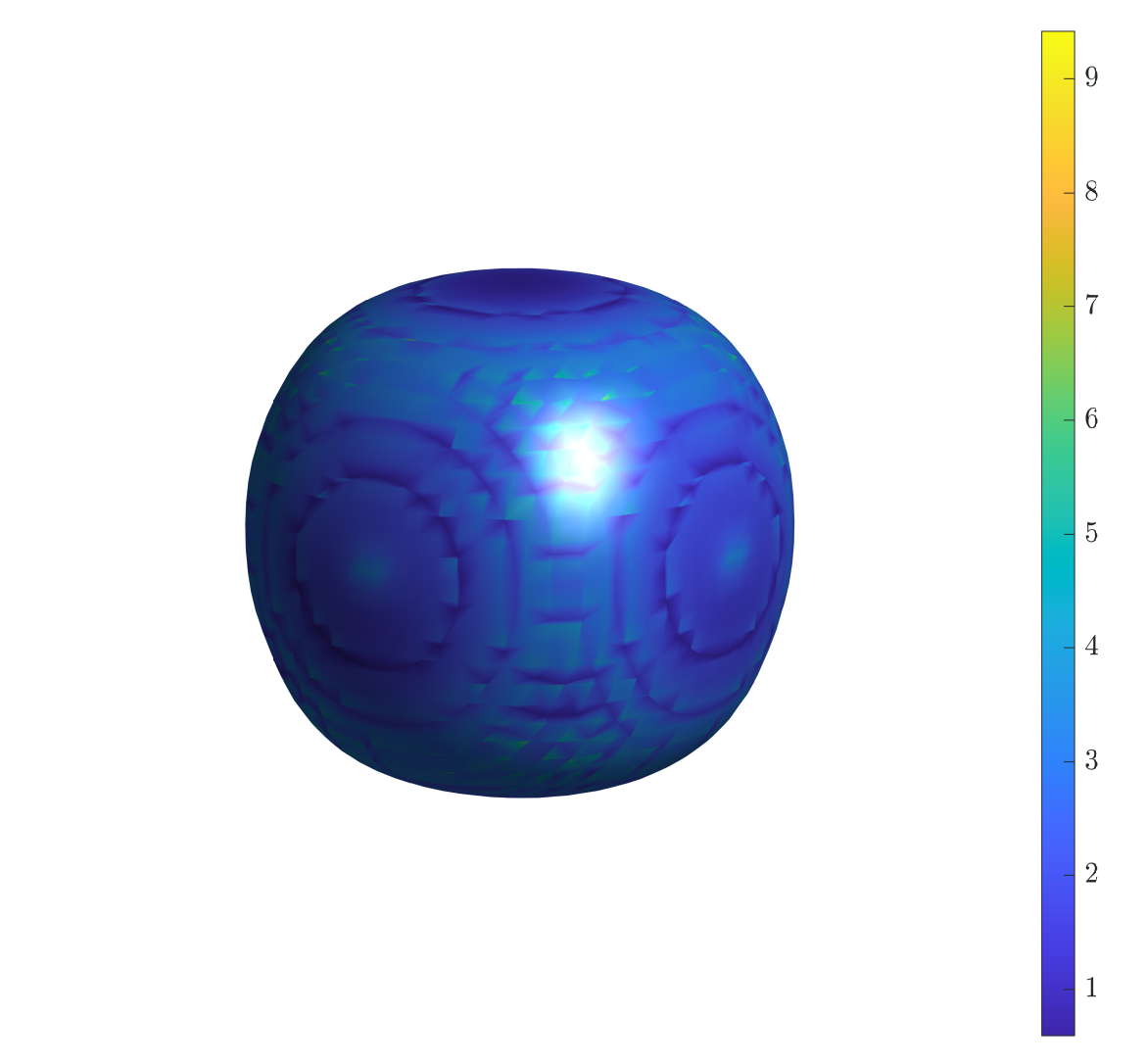}
    \label{p20220215_S_N32S5_NoTITLE_GCMClim_PWEE_14}
    \end{subfigure}
\vspace{-0.8cm}
\cprotect \caption{Visualisation of Gaussian curvatures (top line) and mean curvatures (bottom line) in final reconstructed results of {\bf Example 1} Sphere by three formulations: perimeter-based formulation (left); Willmore-based formulation (middle); and Euler-Elastica-based formulation (right). }
\label{p20220215_S_9_11_13_10_12_14}
\end{figure}

\begin{figure}[H]
    \centering
    \begin{subfigure}[t]{0.48\textwidth}
    \includegraphics[width=\textwidth]{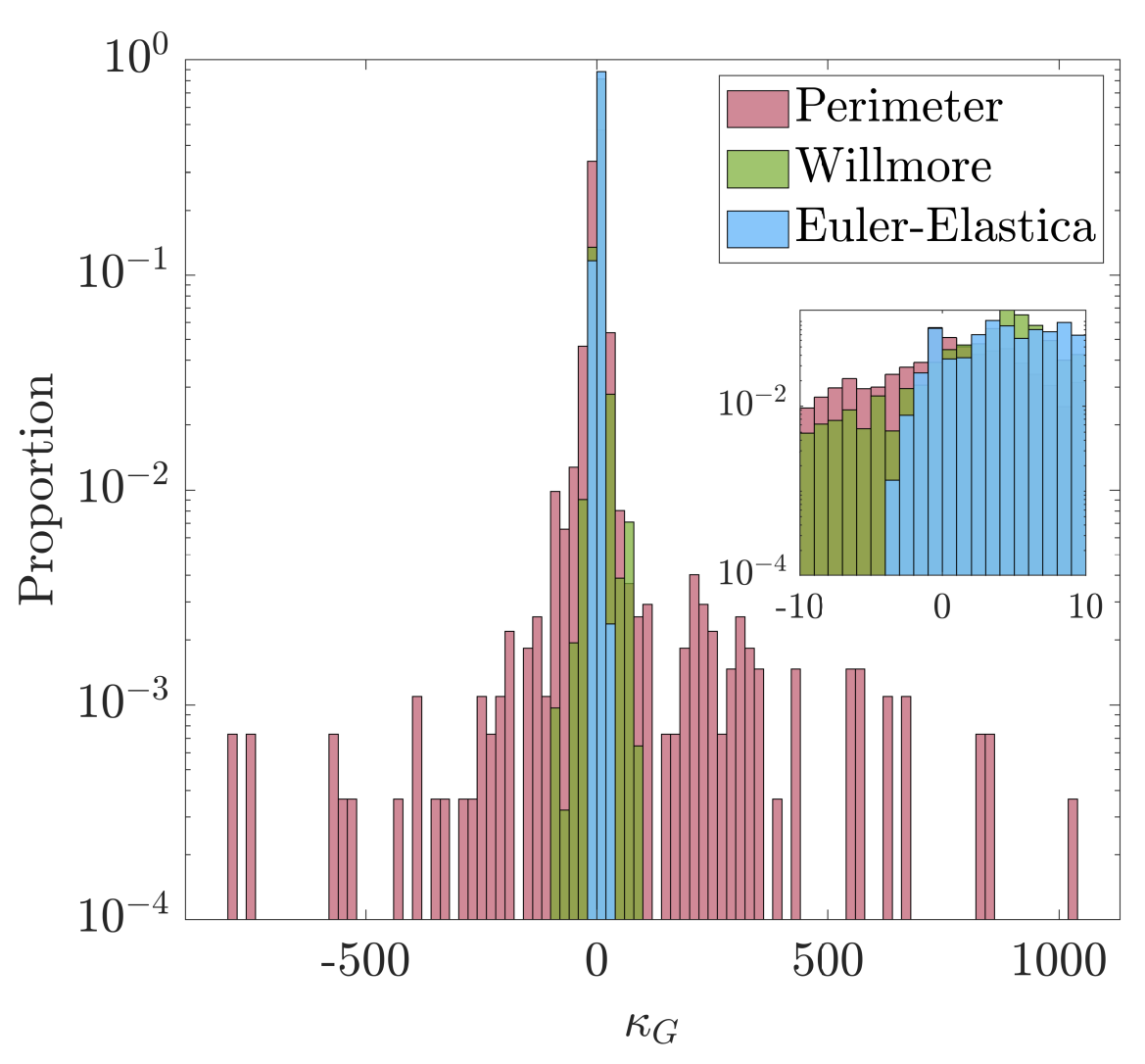}
    \label{p20230319_S_N32S5_Histogram_GCMClim_PWEE_1}
    \end{subfigure}
\hspace{0.1mm}
    \begin{subfigure}[t]{0.48\textwidth}
    \includegraphics[width=\textwidth]{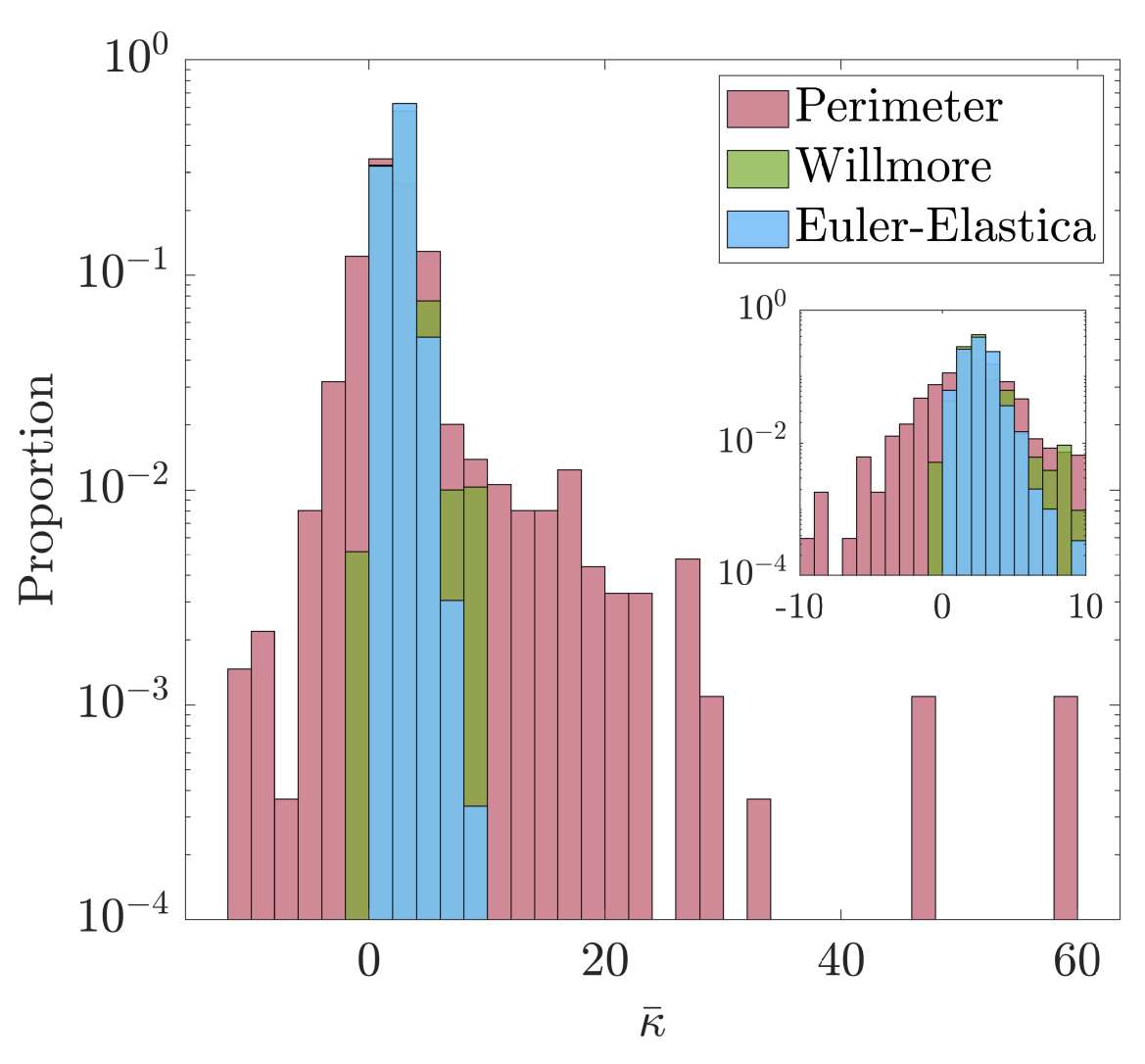}
    \label{p20230319_S_N32S5_Histogram_GCMClim_PWEE_2}
    \end{subfigure}
\vspace{-0.8cm}
\cprotect \caption{Histograms of all Gaussian curvatures $\kappa_{G}$ (left) and mean curvatures $\bar{\kappa}$ (right) in final reconstructed results of {\bf Example 1} Sphere by three formulations: perimeter-based formulation (pink); Willmore-based formulation (green); and Euler-Elastica-based formulation (cyan). }
\label{p20230319_S_N32S5_Histogram_GCMClim_PWEE}
\end{figure}

\begin{figure}[H]
    \centering
    \begin{subfigure}[t]{0.32\textwidth}
    \includegraphics[width=\textwidth]{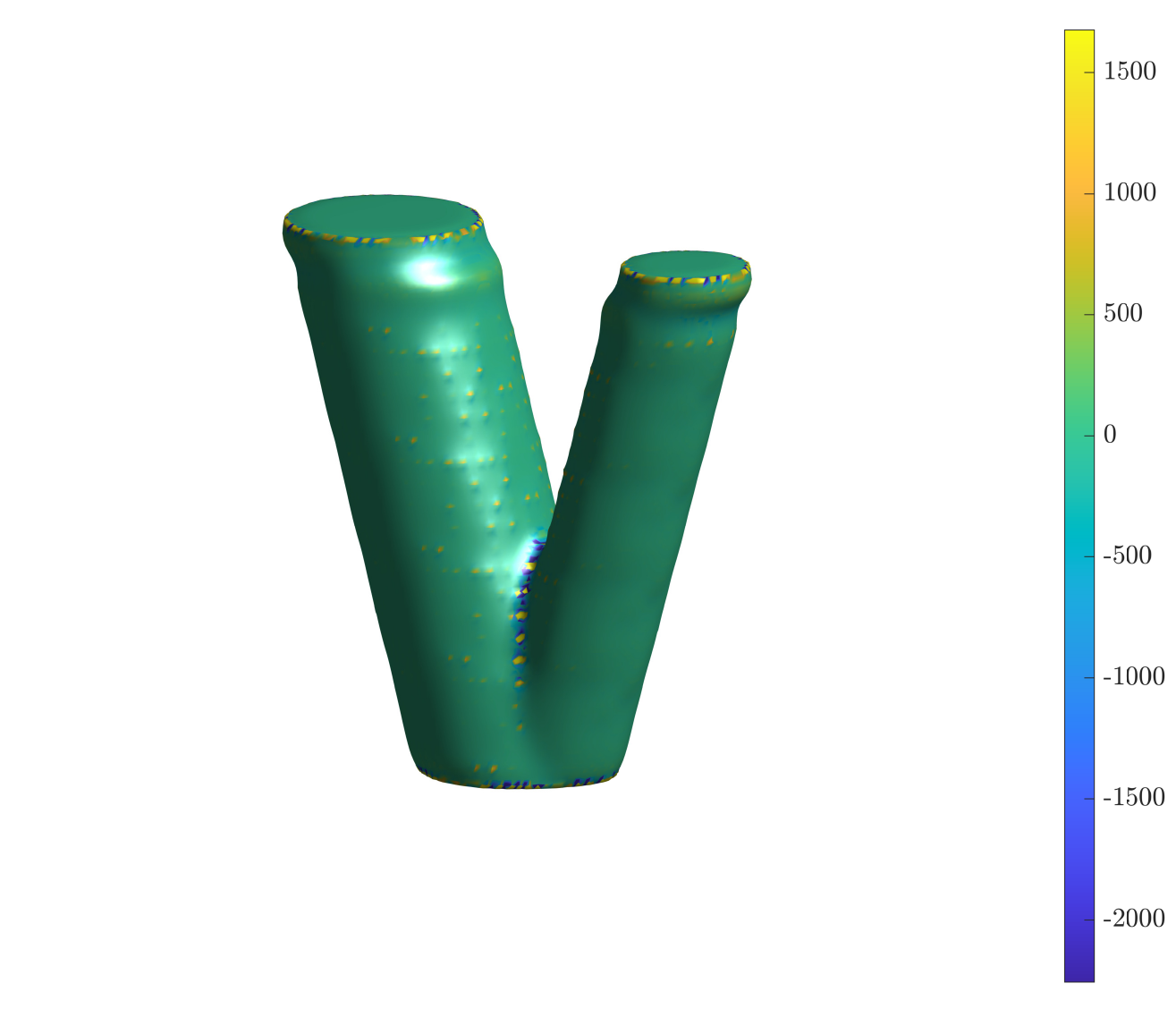}
    \label{p20220215_BC_N128S24_NoTITLE_GCMClim_PWEE_9}
    \end{subfigure}
\hspace{0.1mm}
    \begin{subfigure}[t]{0.32\textwidth}
    \includegraphics[width=\textwidth]{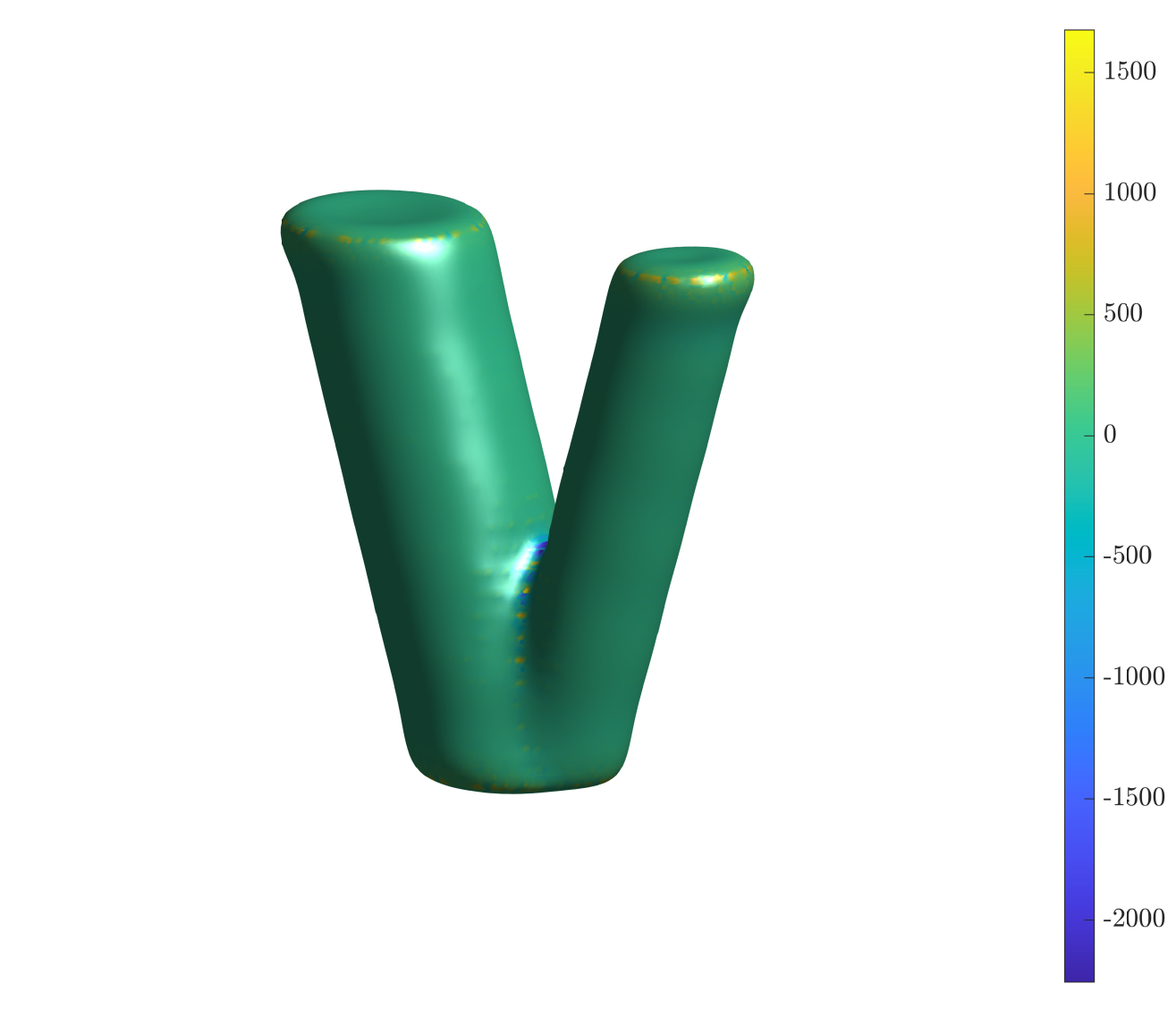}
    \label{p20220215_BC_N128S24_NoTITLE_GCMClim_PWEE_11}
    \end{subfigure}
\hspace{0.1mm}
    \begin{subfigure}[t]{0.32\textwidth}
    \includegraphics[width=\textwidth]{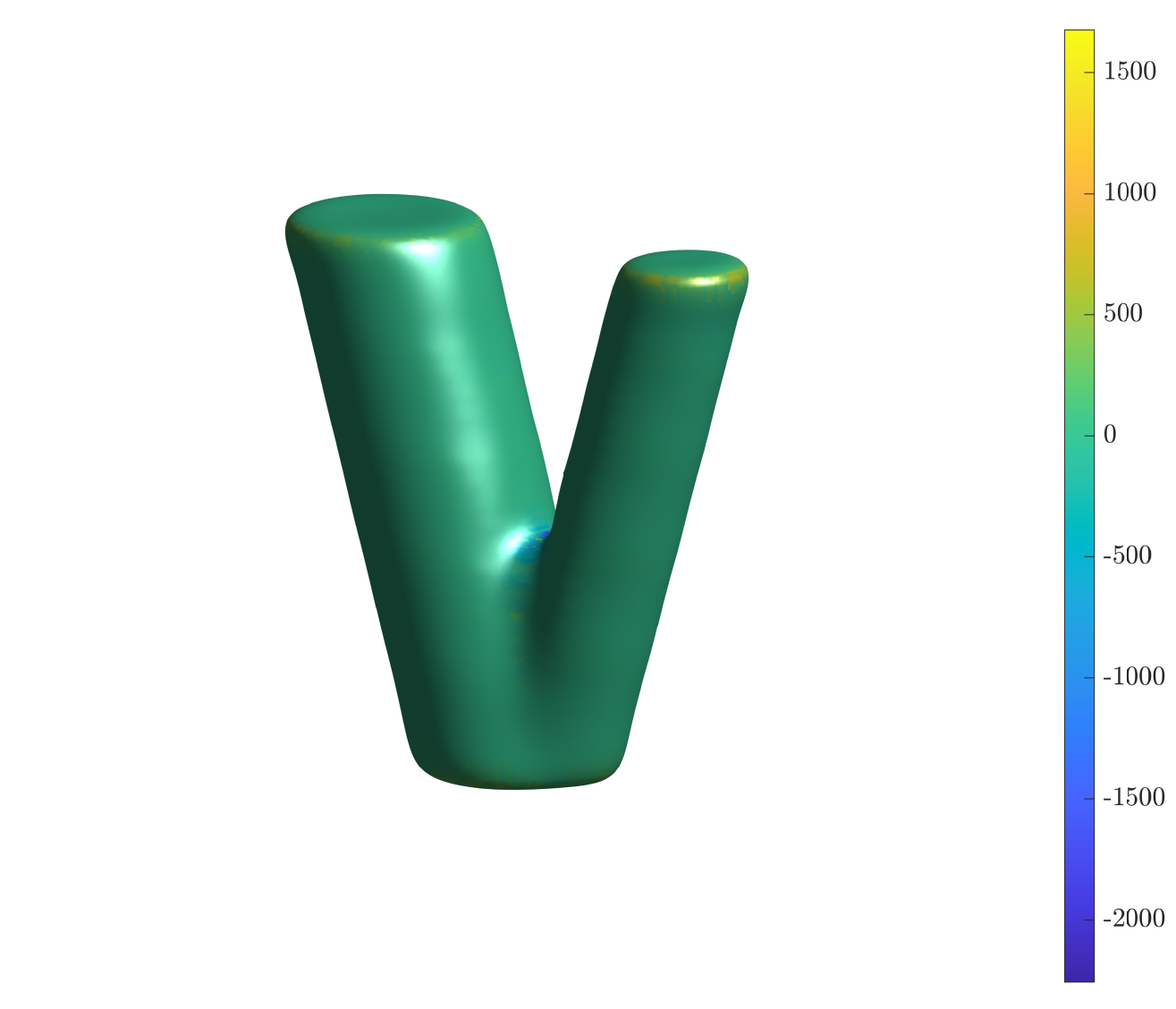}
    \label{p20220215_BC_N128S24_NoTITLE_GCMClim_PWEE_13}
    \end{subfigure}
\\
    \begin{subfigure}[t]{0.32\textwidth}
    \includegraphics[width=\textwidth]{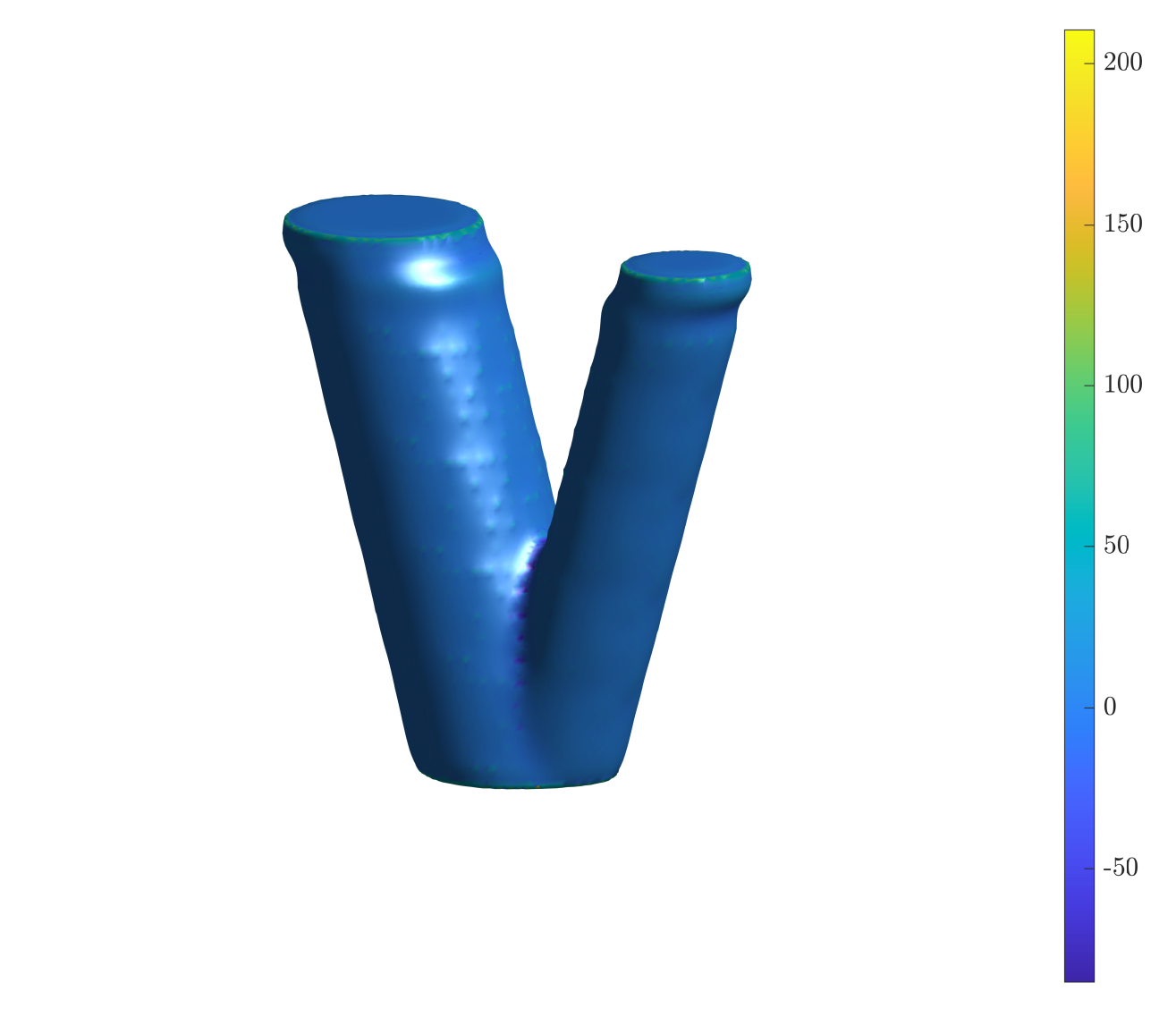}
    \label{p20220215_BC_N128S24_NoTITLE_GCMClim_PWEE_10}
    \end{subfigure}
\hspace{0.1mm}
    \begin{subfigure}[t]{0.32\textwidth}
    \includegraphics[width=\textwidth]{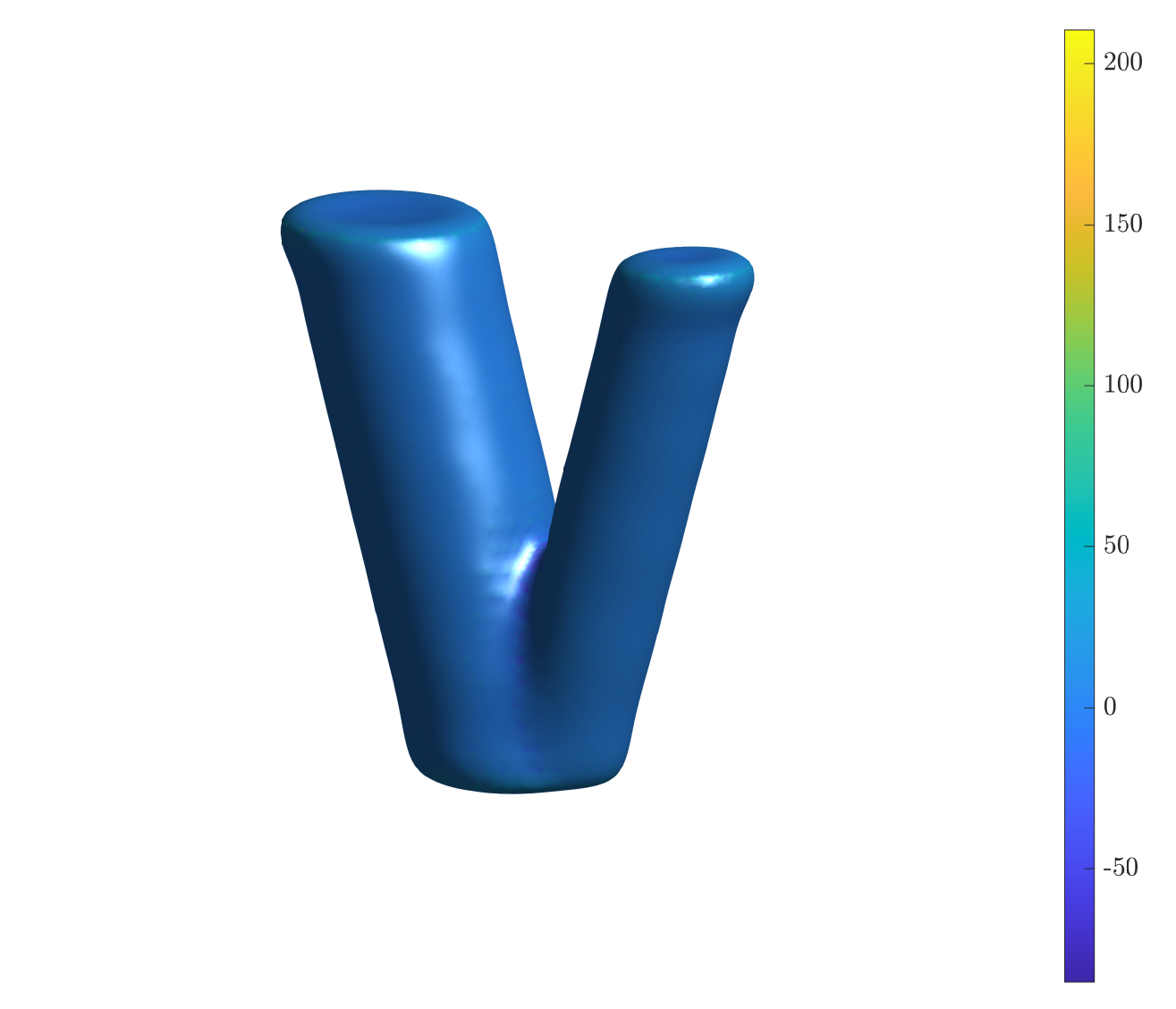}
    \label{p20220215_BC_N128S24_NoTITLE_GCMClim_PWEE_12}
    \end{subfigure}
\hspace{0.1mm}
    \begin{subfigure}[t]{0.32\textwidth}
    \includegraphics[width=\textwidth]{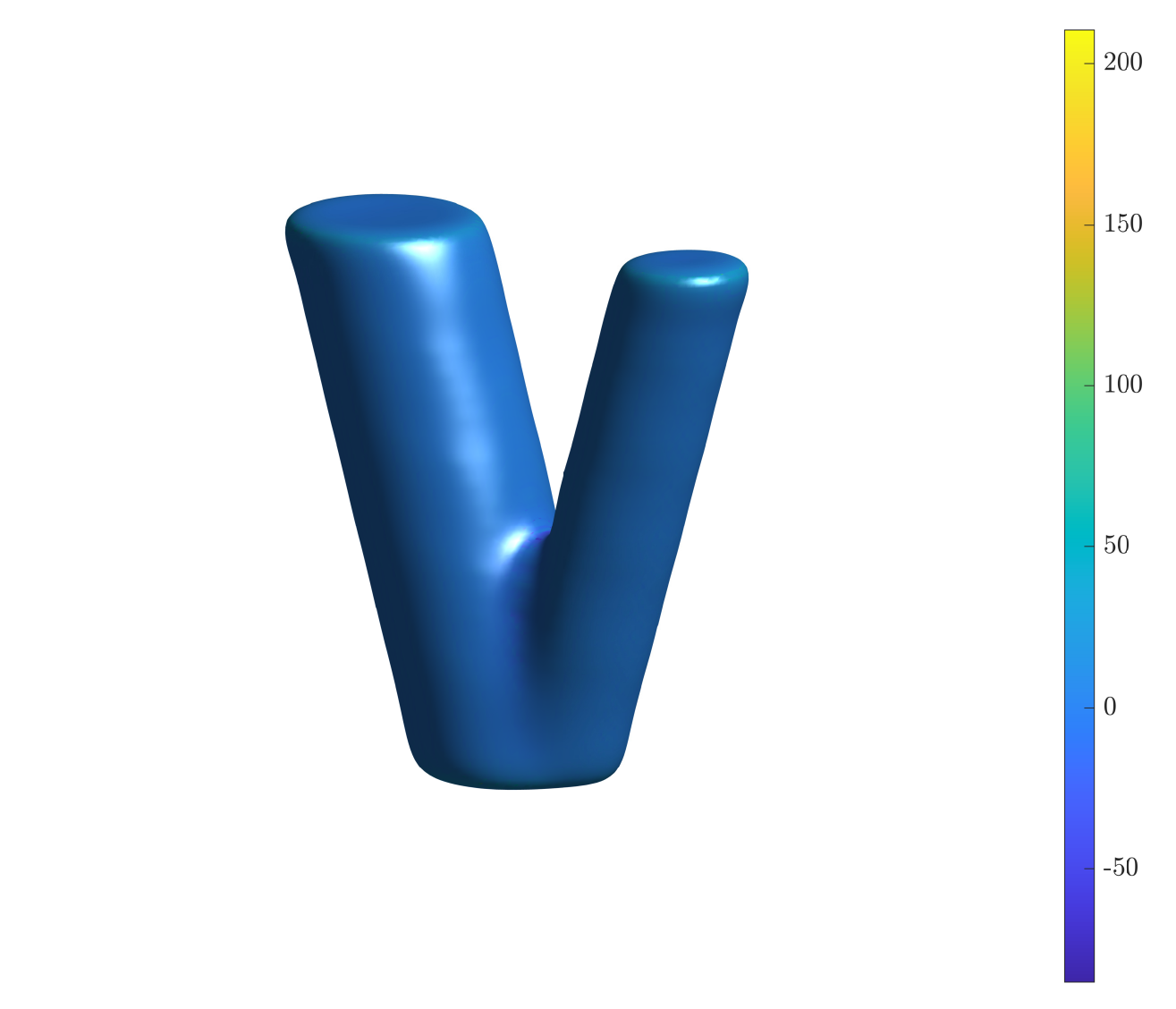}
    \label{p20220215_BC_N128S24_NoTITLE_GCMClim_PWEE_14}
    \end{subfigure}
\vspace{-0.8cm}
\cprotect \caption{Visualisation of Gaussian curvatures (top line) and mean curvatures (bottom line) in final reconstructed results of {\bf Example 2} Branching Cylinders by three formulations: perimeter-based formulation (left); Willmore-based formulation (middle); and Euler-Elastica-based formulation (right). }
\label{p20220215_BC_9_11_13_10_12_14}
\end{figure}

\begin{figure}[H]
    \centering
    \begin{subfigure}[t]{0.48\textwidth}
    \includegraphics[width=\textwidth]{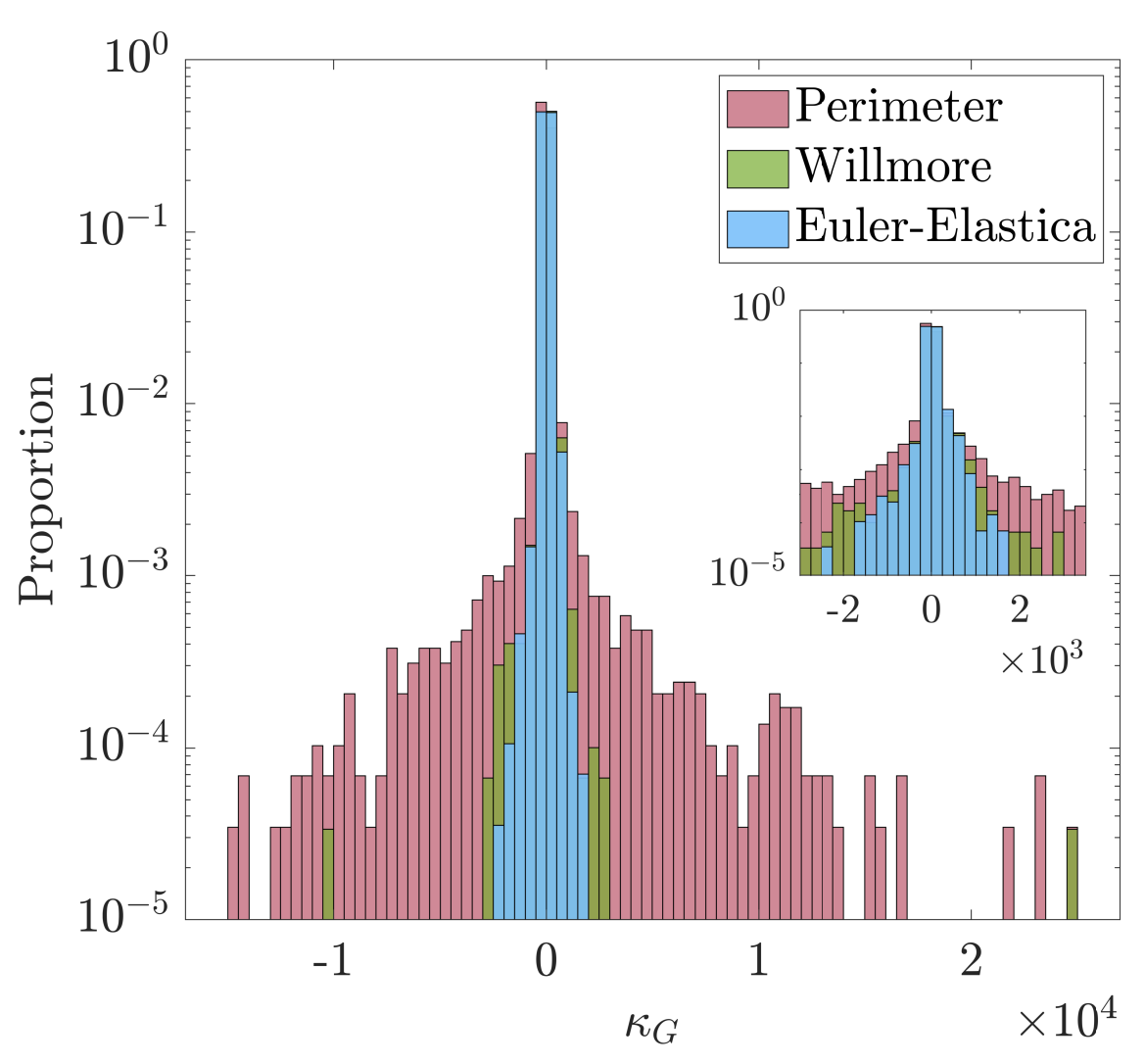}
    \label{p20230319_BC_N128S24_Histogram_GCMClim_PWEE_1}
    \end{subfigure}
\hspace{0.1mm}
    \begin{subfigure}[t]{0.48\textwidth}
    \includegraphics[width=\textwidth]{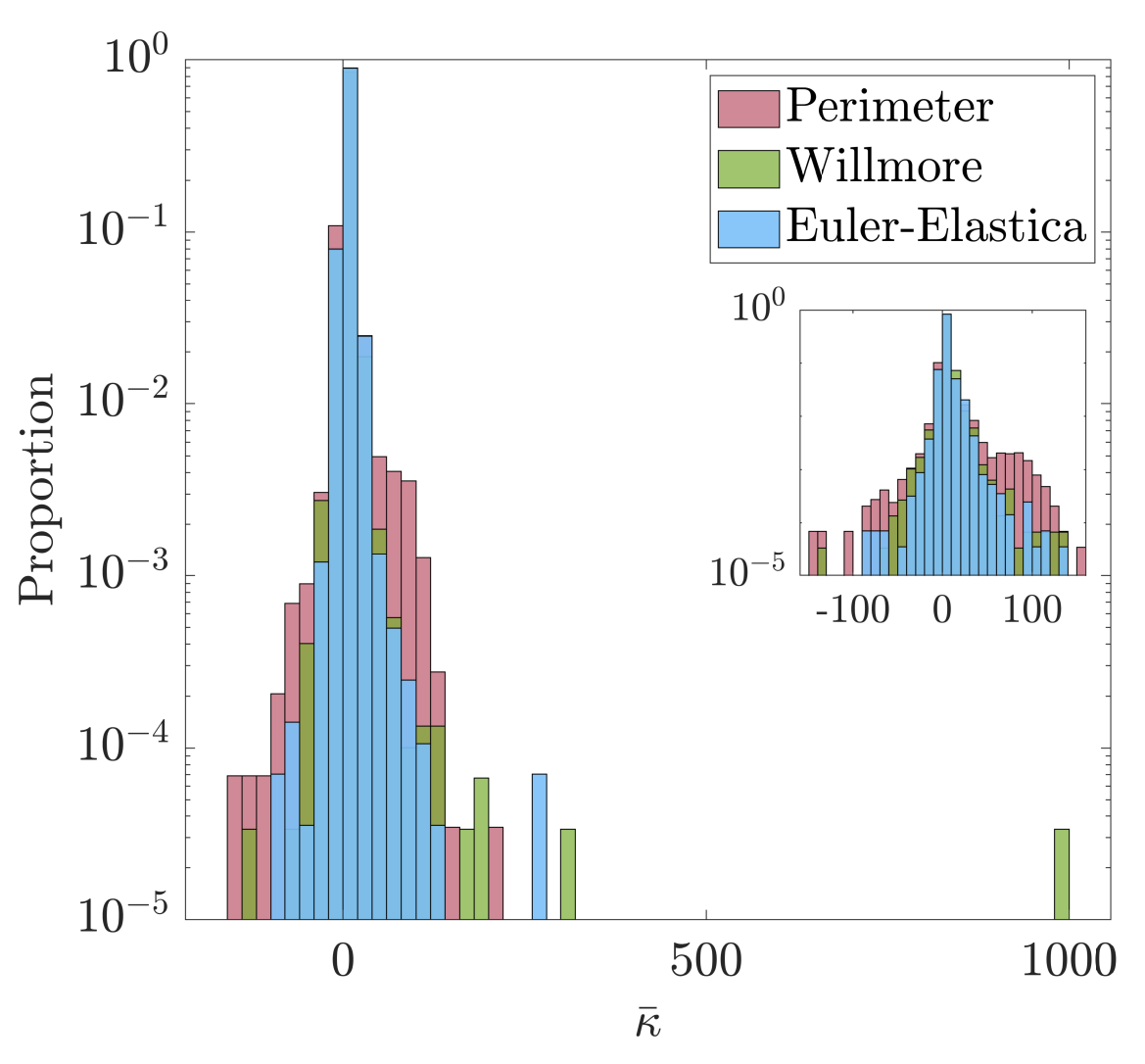}
    \label{p20230319_BC_N128S24_Histogram_GCMClim_PWEE_2}
    \end{subfigure}
\vspace{-0.8cm}
\cprotect \caption{Histograms of all Gaussian curvatures $\kappa_{G}$ (left) and mean curvatures $\bar{\kappa}$ (right) in final reconstructed results of {\bf Example 2} Branching Cylinders by three formulations: perimeter-based formulation (pink); Willmore-based formulation (green); and Euler-Elastica-based formulation (cyan). }
\label{p20230319_BC_N128S24_Histogram_GCMClim_PWEE}
\end{figure}

\begin{figure}[htbp]
    \centering
    \vspace{-0.3cm}
    \includegraphics[width=.7\textwidth]{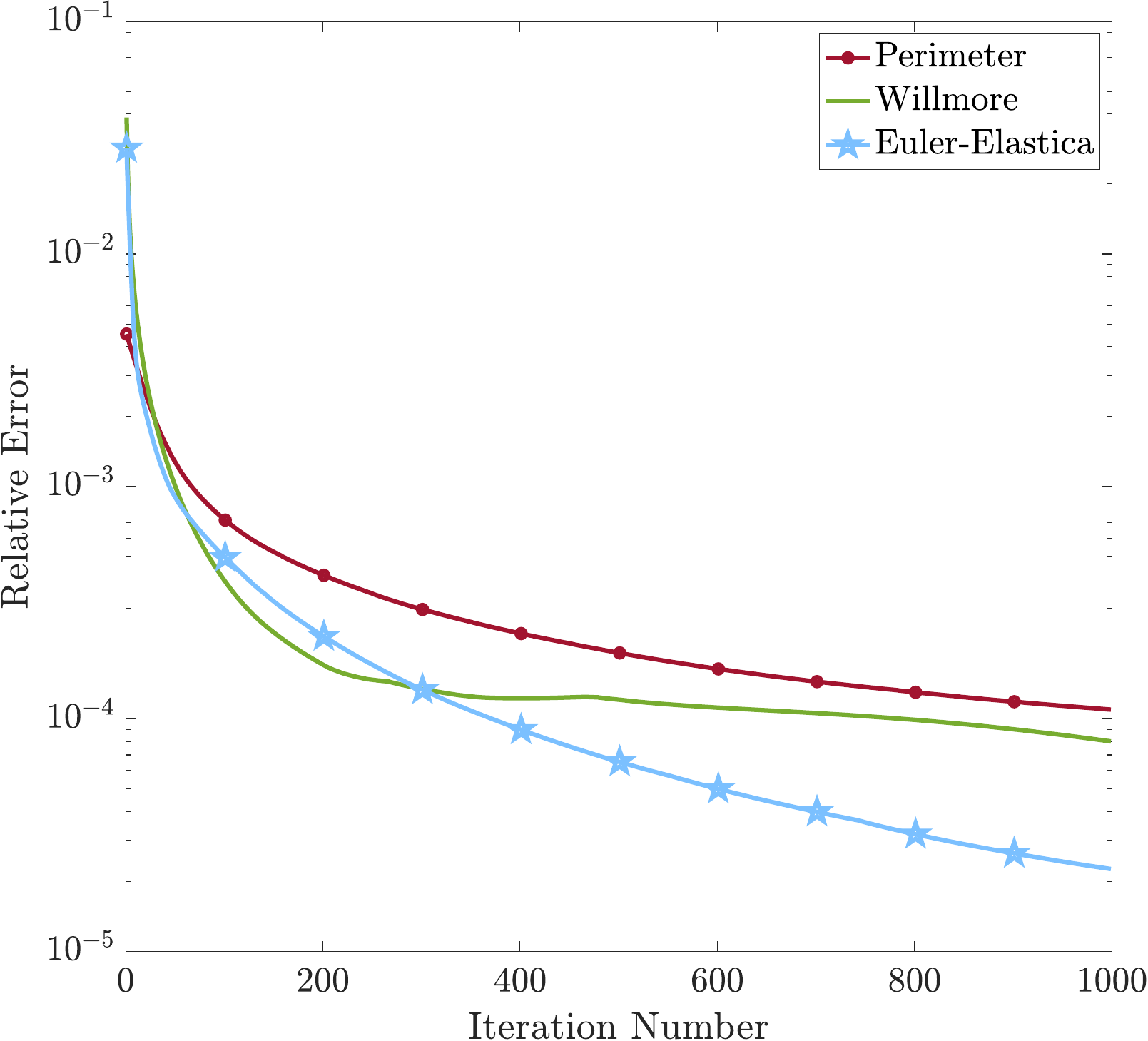}
\vspace{-0.2cm}
\cprotect \caption{Experimental convergence curves by the relative error over the number of iterations for {\bf Example 1} with respect to the perimeter-based (red with circles), the Willmore-based (green), and the Euler-Elastica-based (cyan with pentagrams) formulation. }
\label{p20230327_S_N32S5_RelativeError_PWEE_1}
\vspace{-0.2cm}
\end{figure}

\begin{figure}[htbp]
    \centering
    \includegraphics[width=.7\textwidth]{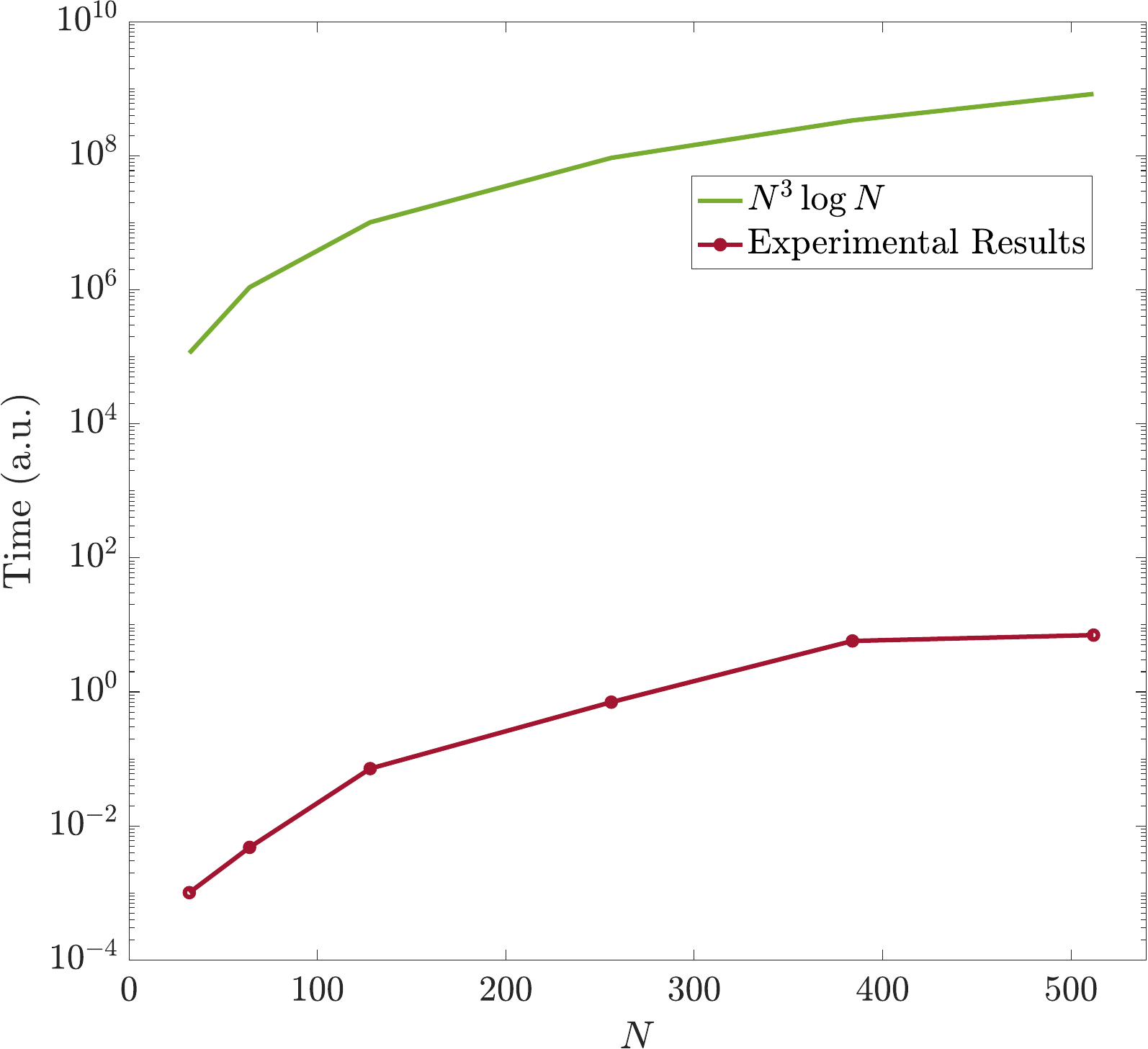}
\cprotect \caption{Computational complexity curves of experimental results (red with circles) and estimated arithmetical values (green) with various inputs $N$ for {\bf Example 1} where the unit of time is set as the arbitrary unit (a.u.). }
\label{p20230327_S_N32S5_ComputationalComplexity_PWEE_1}
\end{figure}

\begin{figure}[htbp]
    \centering
    \begin{subfigure}[t]{0.49\textwidth}
    \includegraphics[width=\textwidth]{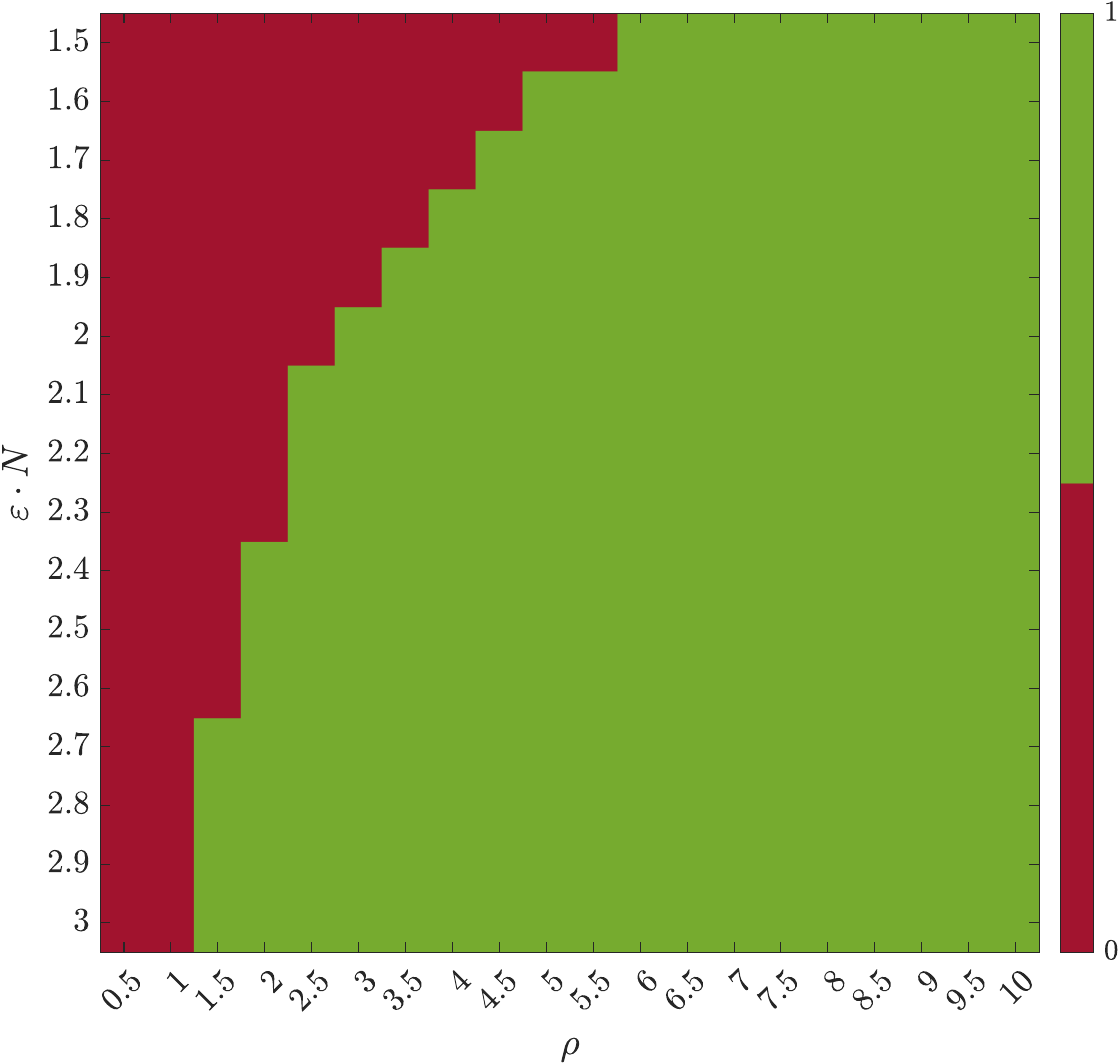}
    \caption{}
    \label{dt3/p20220703_S_N32S5_EE_ADMM_PT_1}
    \end{subfigure}
\hspace{0.1mm}
    \begin{subfigure}[t]{0.49\textwidth}
    \includegraphics[width=\textwidth]{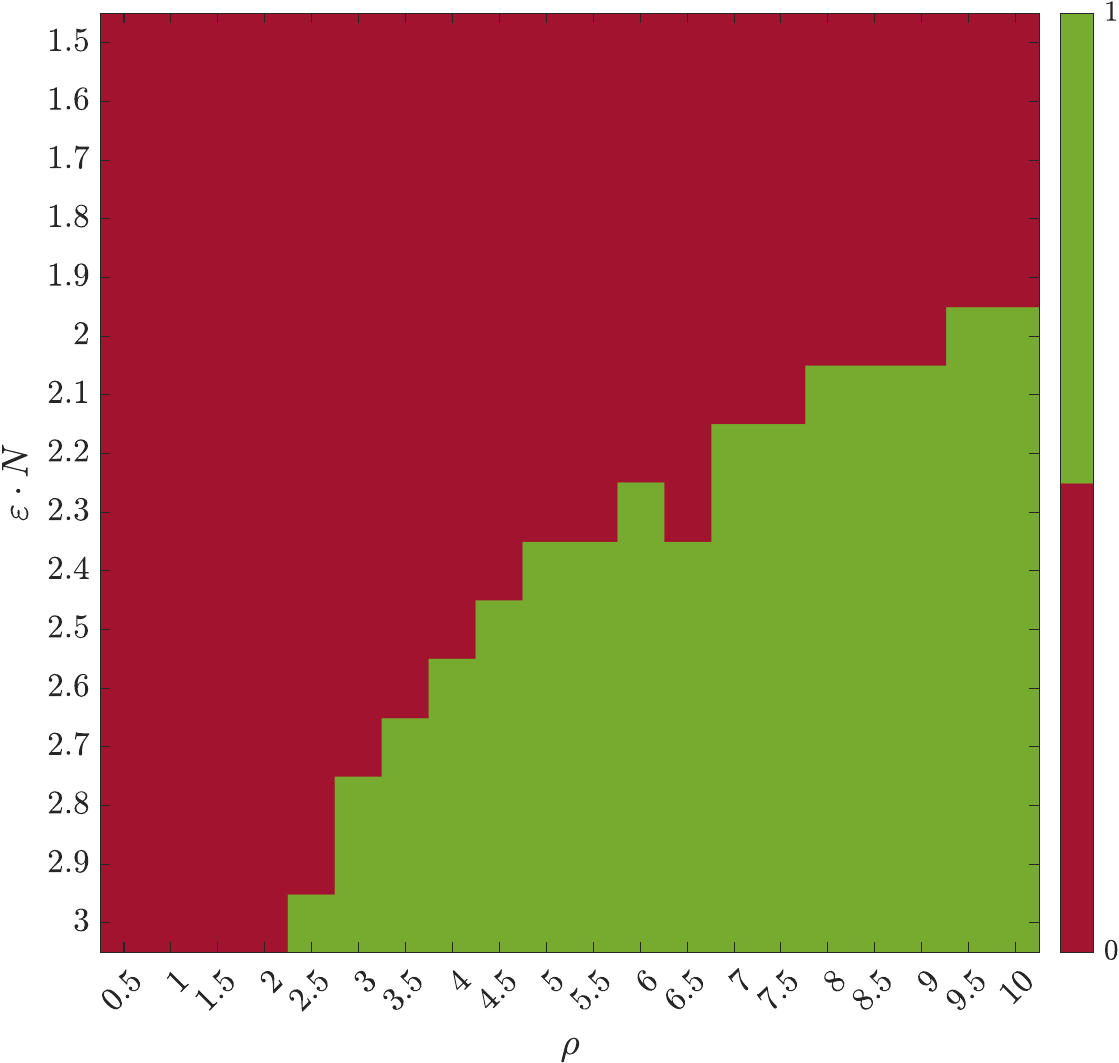}
    \caption{}
    \label{dt3.5/p20220704_S_N32S5_EE_ADMM_PT_1}
    \end{subfigure}
\vspace{-0.2cm}
\cprotect \caption{Binary maps of testing parameters sensitivity with respect to the penalty parameter $\rho \in [0.5, 10]$ with the step of $0.5$ and the representative diffuse interface width $\varepsilon \cdot N \in [1.5, 3]$ with the step of $0.1$ for {\bf Example 1} Sphere under the low resolution $N = 32$ by the numerical algorithm II based on ADMM where $1$ (green) stands for the tolerable results existed and $0$ (red) indicates the unpleasant results during the iterations under the different time steps~\subref{dt3/p20220703_S_N32S5_EE_ADMM_PT_1} $\tau = \varepsilon^{3}$ and~\subref{dt3.5/p20220704_S_N32S5_EE_ADMM_PT_1} $\tau = \varepsilon^{3.5}$. }
\label{p20220703_S_N32S5_EE_ADMM_PT}
\end{figure}

\begin{figure}[htbp]
    \centering
    \begin{subfigure}[t]{0.45\textwidth}
    \includegraphics[width=\textwidth]{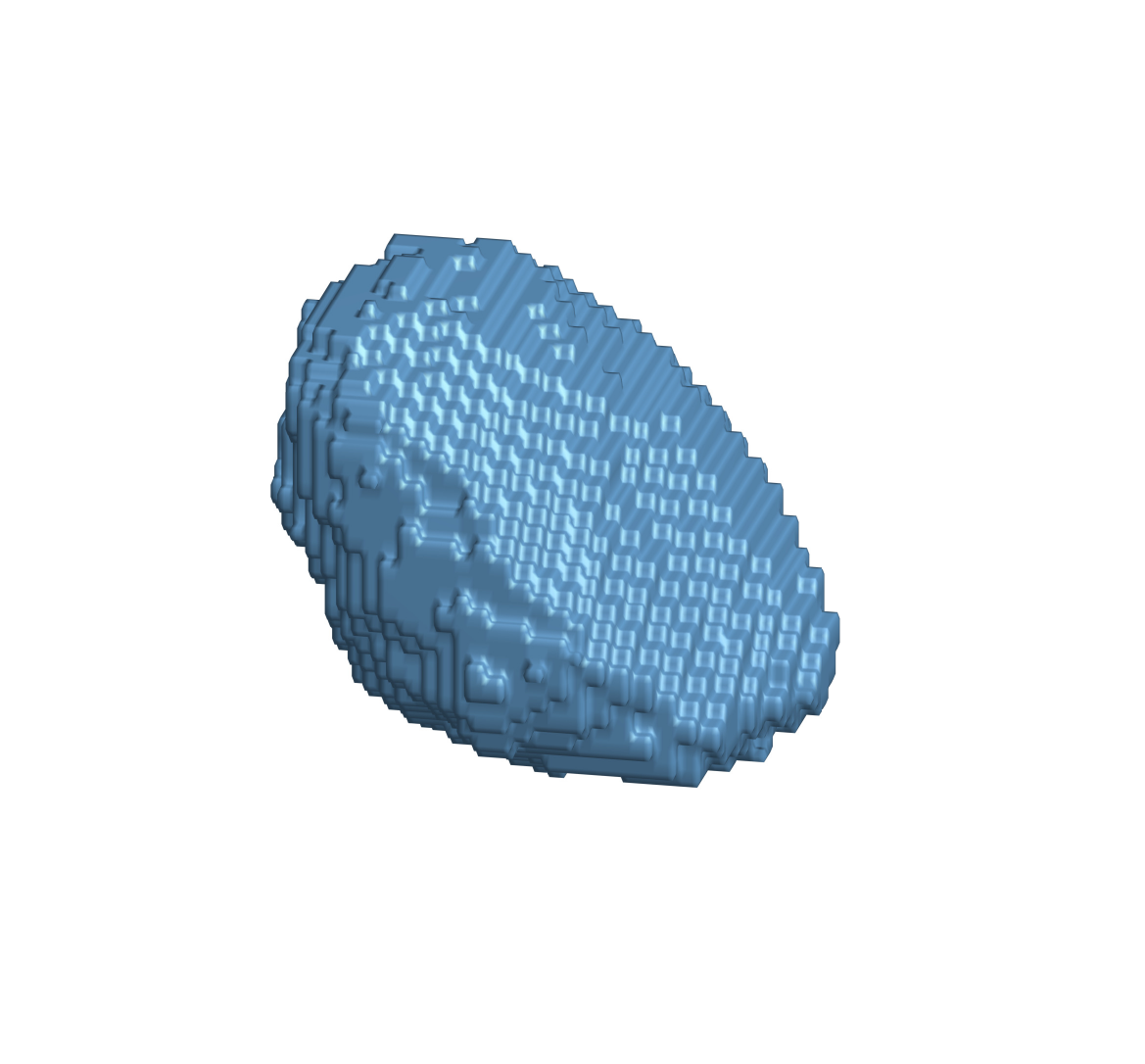}    \label{p20220309_RealData_Cancer_EE_1}
    \end{subfigure}
\hspace{0.1mm}
    \begin{subfigure}[t]{0.45\textwidth}
    \includegraphics[width=\textwidth]{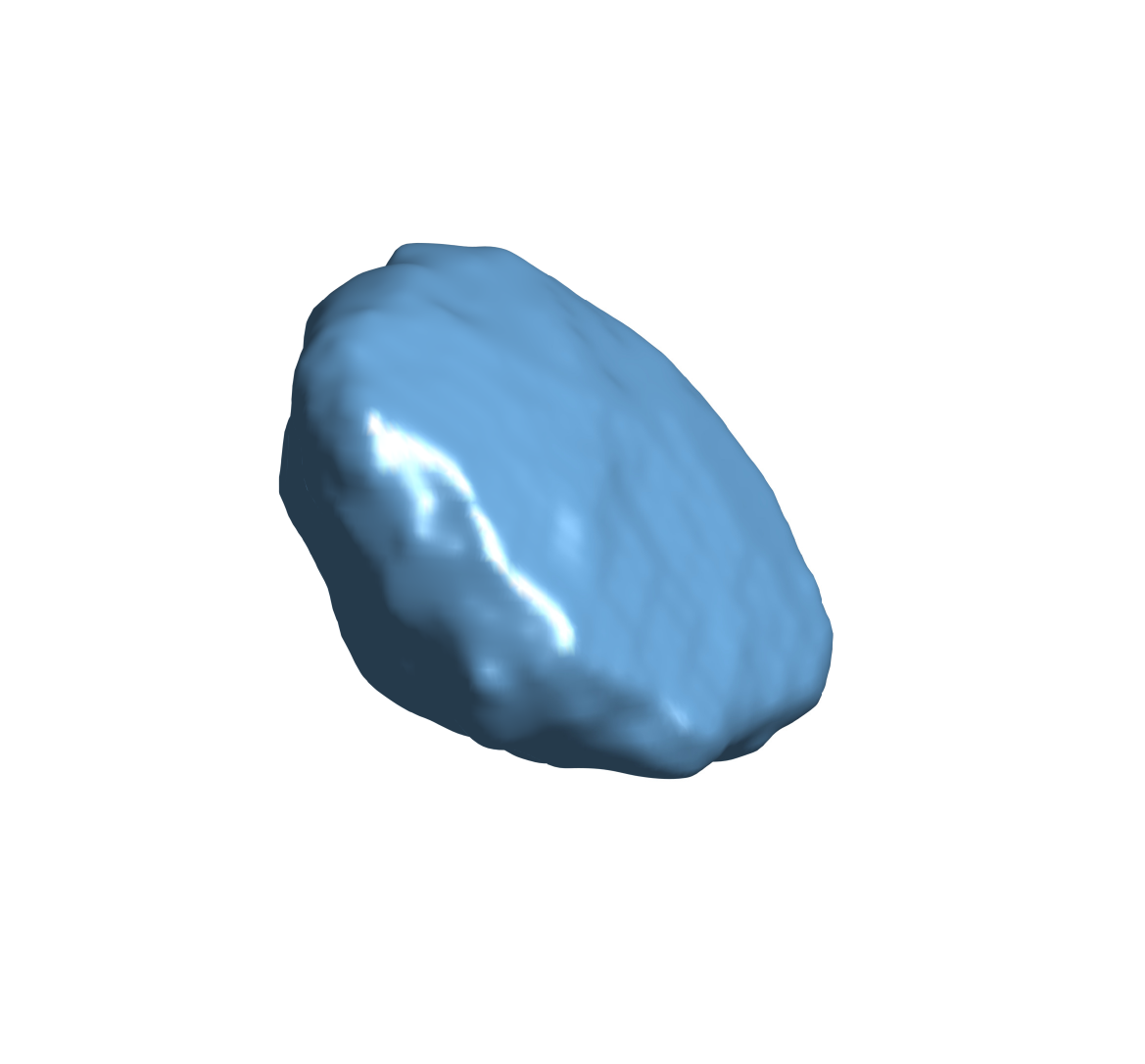}
    \label{p20220309_RealData_Cancer_EE_2}
    \end{subfigure}
\vspace{-1cm}
\cprotect \caption{Visualisation of segmented rough input (left) of {\bf Example 4} Tumour from the given $280$ slices of MRI real data and the smoothed reconstruction results by the Euler-Elastica-based formulation (right). }
\label{p20220309_RealData_Cancer_EE_1-2}
\end{figure}

\begin{figure}[htbp]
    \centering
    \begin{subfigure}[t]{0.45\textwidth}
    \vspace{-.4cm}
    \includegraphics[width=\textwidth]{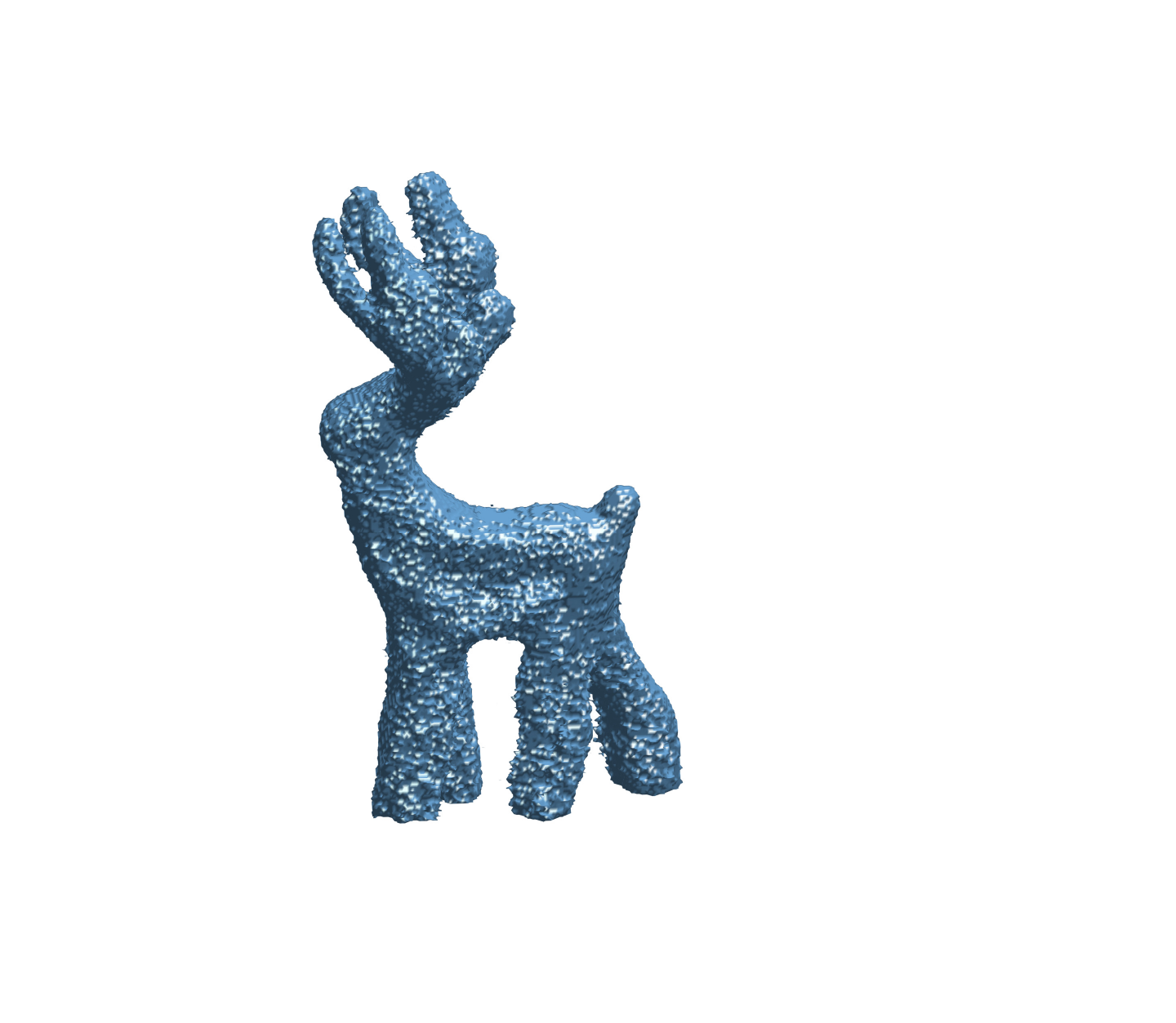}
    \vspace{-1cm}
    \caption{}
    \label{p20230323_THz_Deer_S218G0_EE_Evolution_1}
    \end{subfigure}
\hspace{0.1mm}
    \begin{subfigure}[t]{0.45\textwidth}
    \vspace{-.4cm}
    \includegraphics[width=\textwidth]{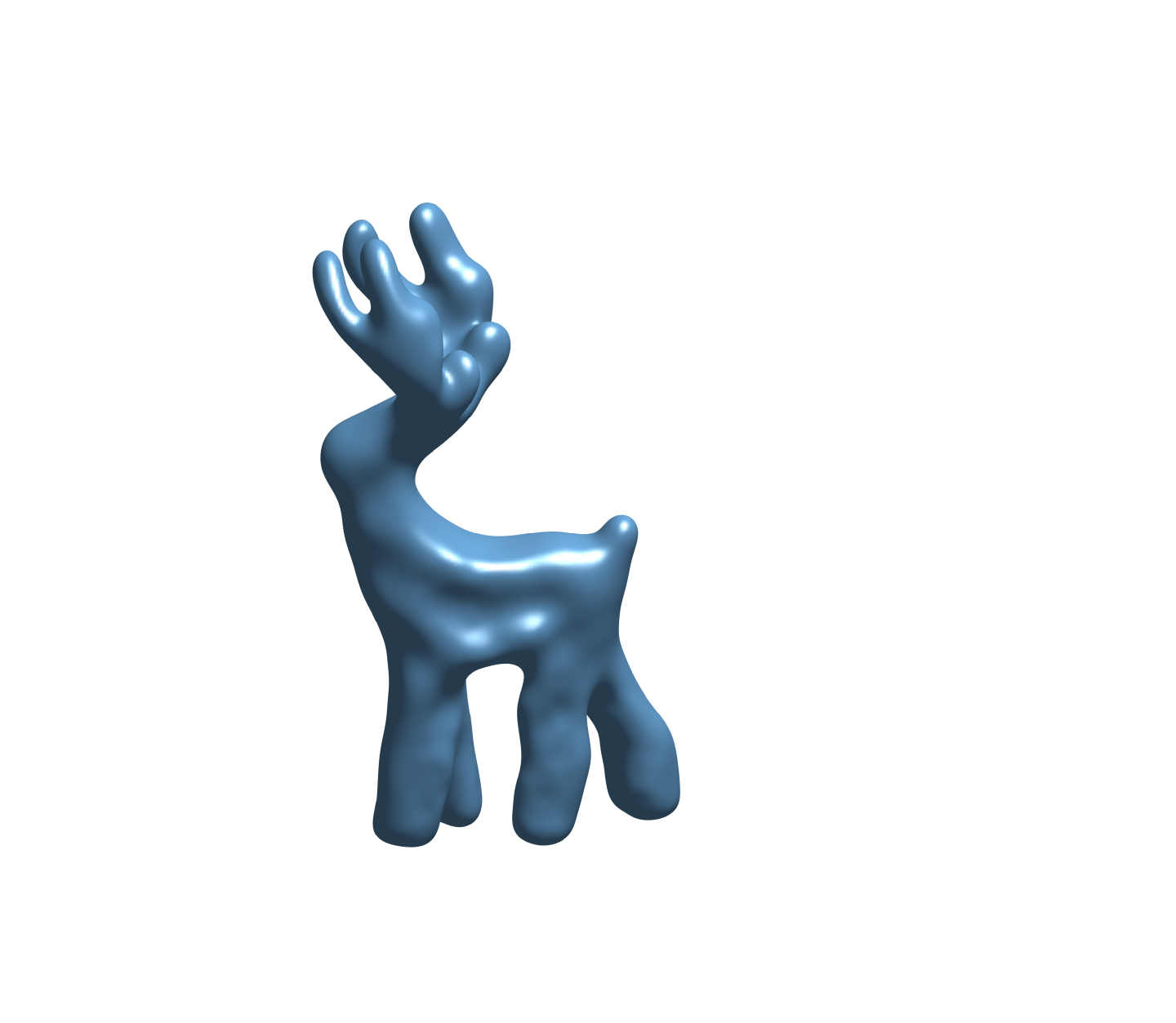}
    \vspace{-1cm}
    \caption{}
    \label{p20230323_THz_Deer_S218G0_EE_Evolution_4}
    \end{subfigure}
    \\
    \vspace{-.5cm}
    \begin{subfigure}[t]{0.45\textwidth}
    \vspace{-.1cm}
    \includegraphics[width=\textwidth]{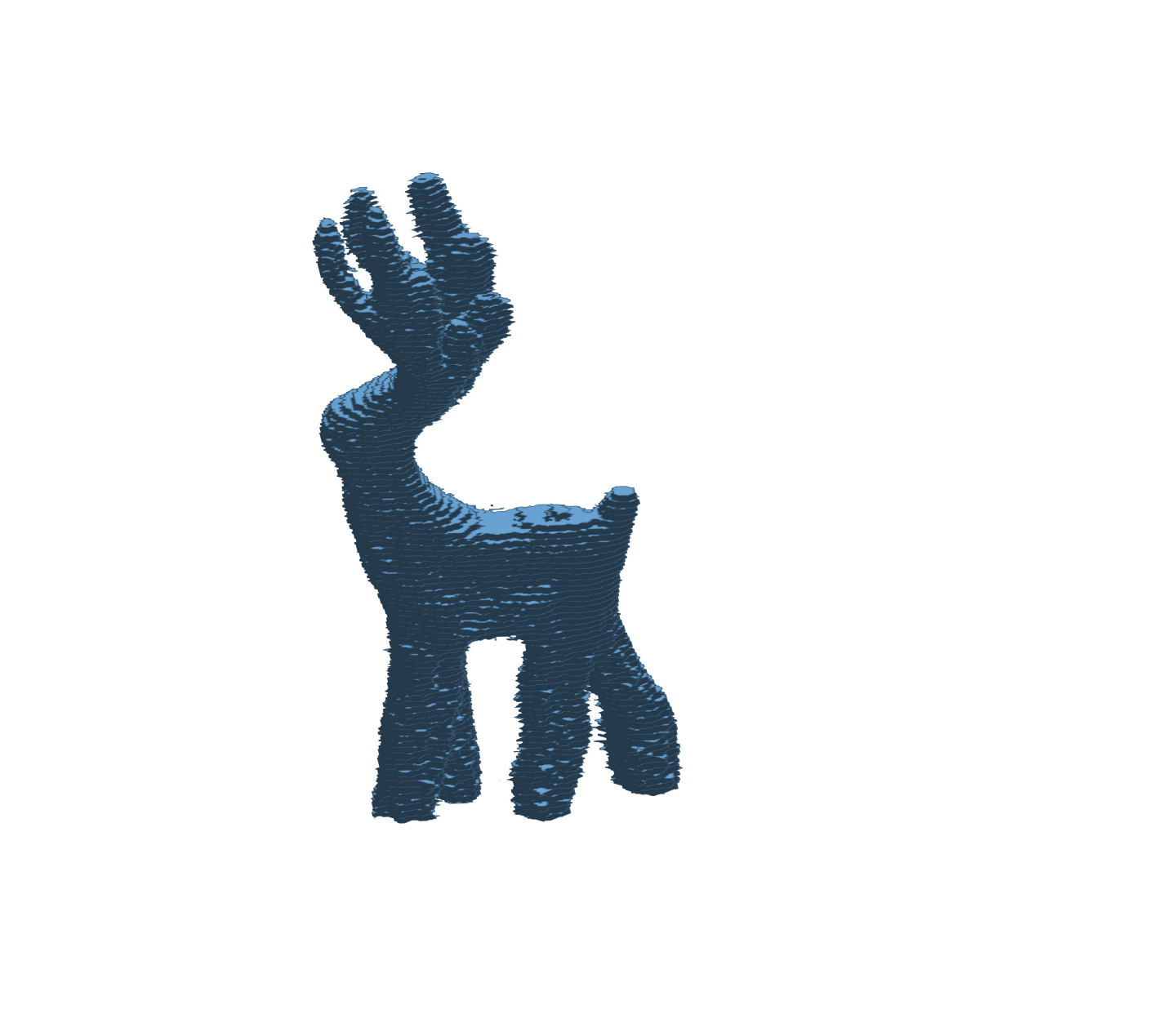}
    \vspace{-1cm}
    \caption{}
    \label{p20230323_THz_Deer_S109G1_EE_Evolution_5}
    \end{subfigure}
\hspace{0.1mm}
    \begin{subfigure}[t]{0.45\textwidth}
    \vspace{-.1cm}
    \includegraphics[width=\textwidth]{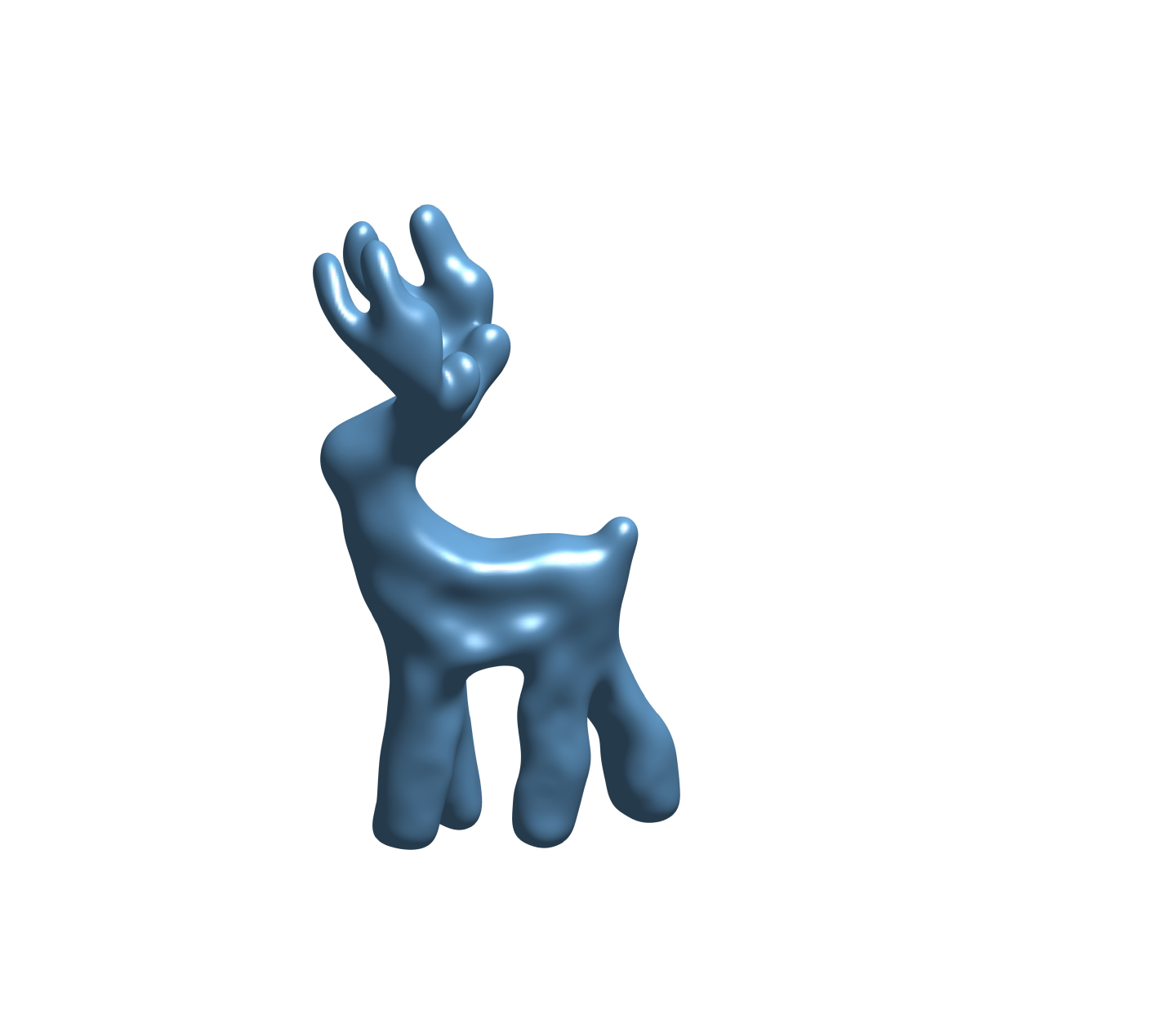}
    \vspace{-1cm}
    \caption{}
    \label{p20230323_THz_Deer_S109G1_EE_Evolution_8}
    \end{subfigure}
    \\\vspace{-.5cm}
    \begin{subfigure}[t]{0.45\textwidth}
    \vspace{-.1cm}
    \includegraphics[width=\textwidth]{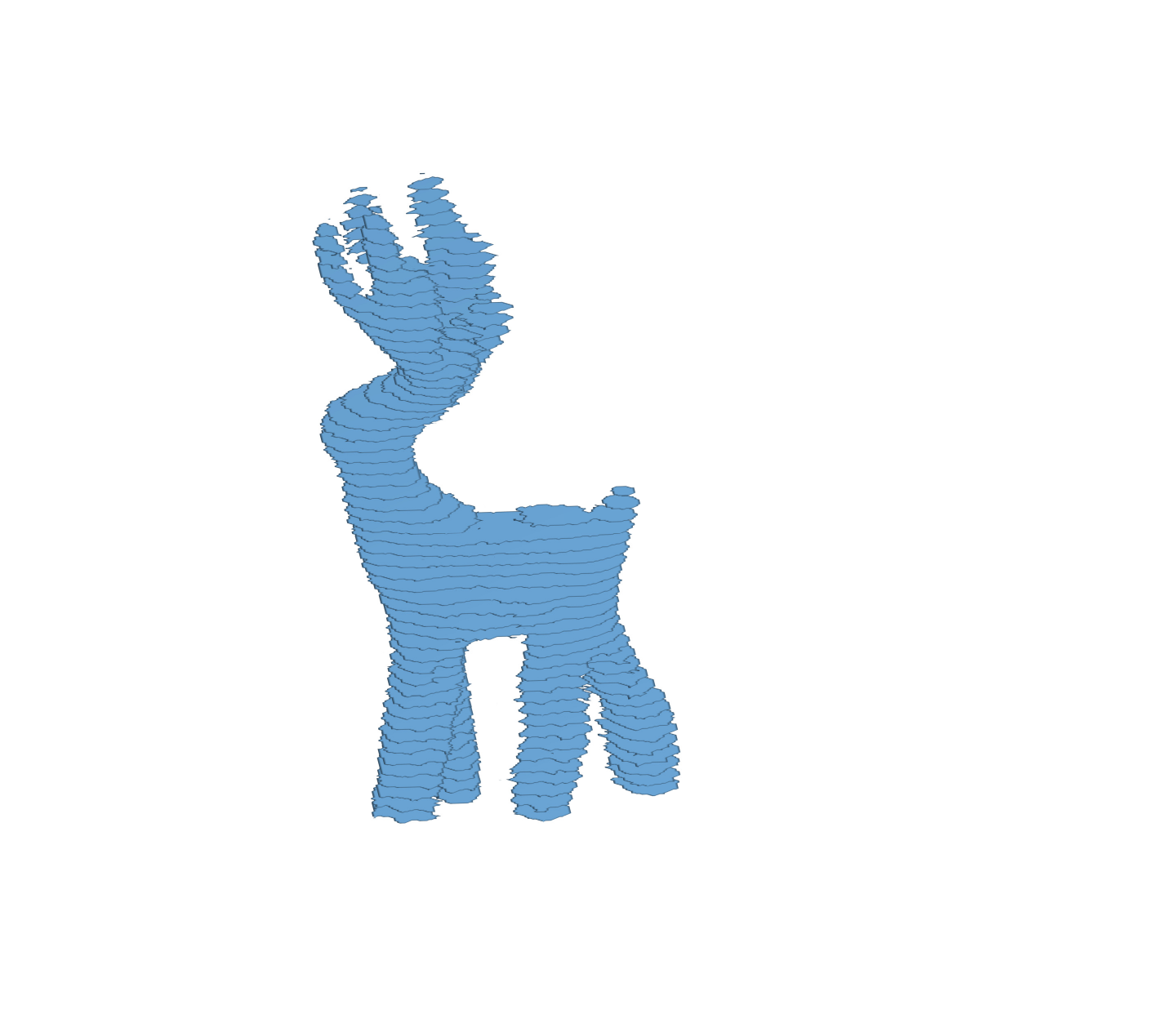}
    \vspace{-1cm}
    \caption{}
    \label{p20230323_THz_Deer_S55G3_EE_Evolution_13}
    \end{subfigure}
\hspace{0.1mm}
    \begin{subfigure}[t]{0.45\textwidth}
    \vspace{-.1cm}
    \includegraphics[width=\textwidth]{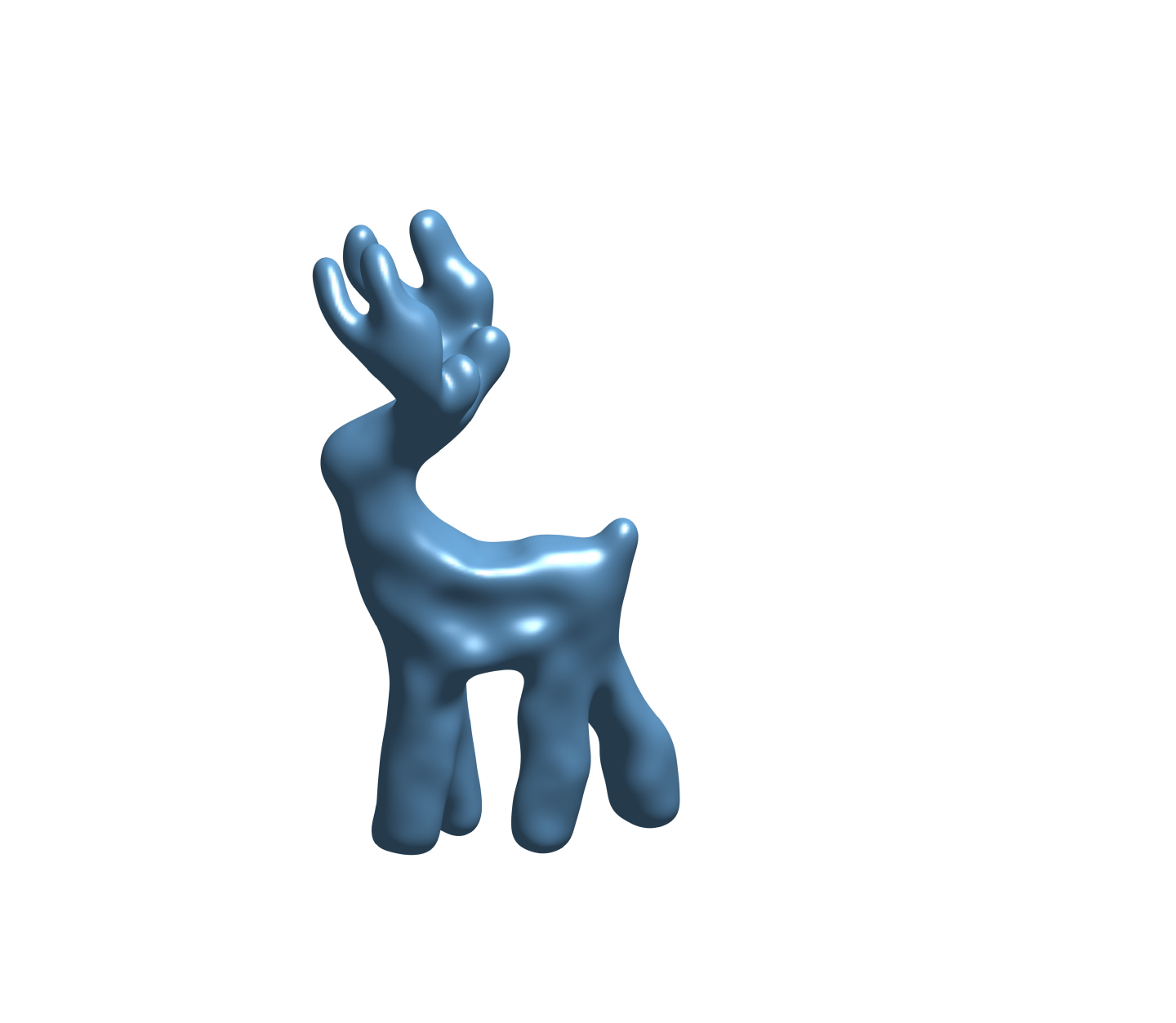}
    \vspace{-1cm}
    \caption{}
    \label{p20230323_THz_Deer_S55G3_EE_Evolution_16}
    \end{subfigure}
    \\\vspace{-.5cm}
    \begin{subfigure}[t]{0.45\textwidth}
    \vspace{-.1cm}
    \includegraphics[width=\textwidth]{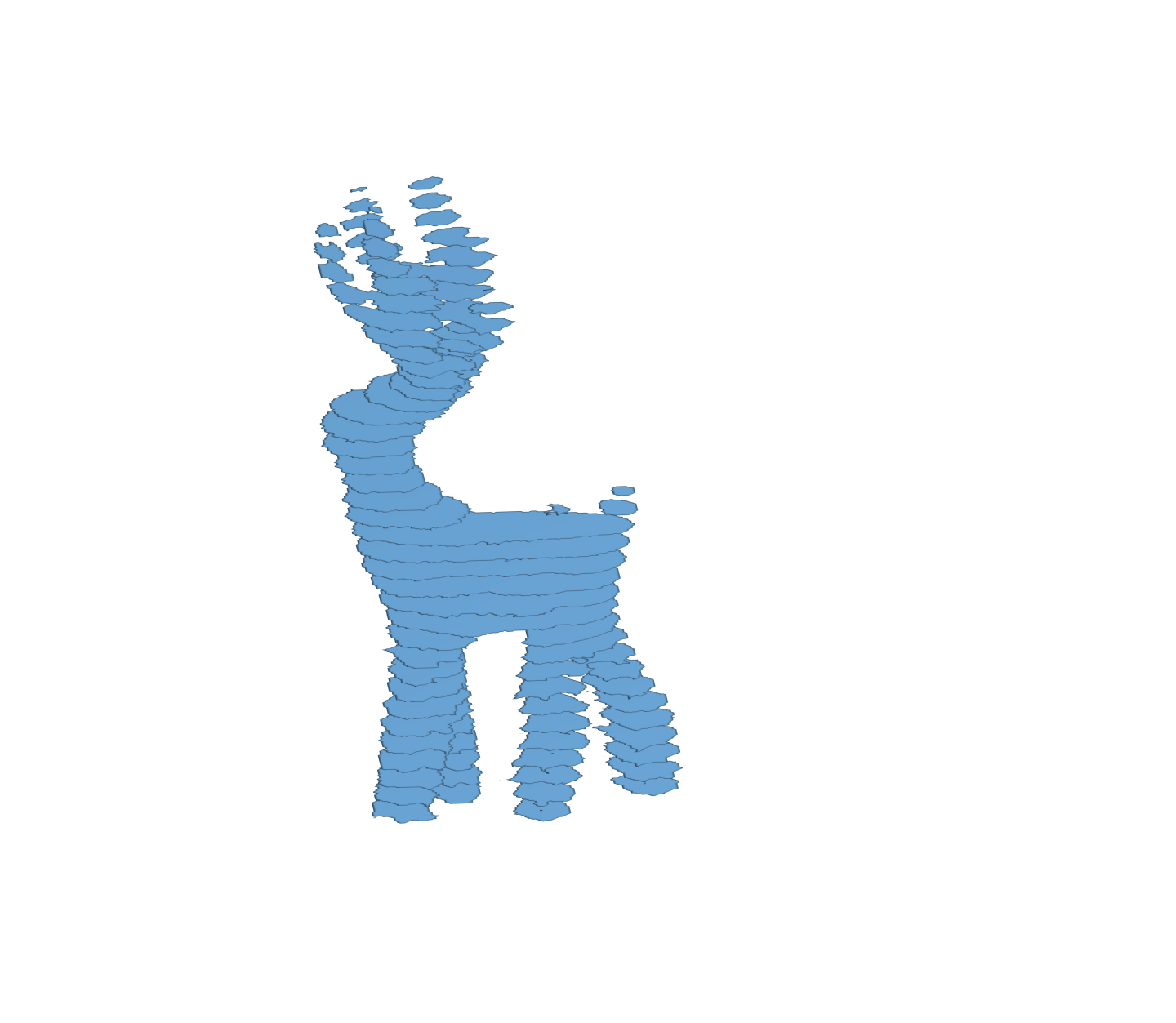}
    \vspace{-1cm}
    \caption{}
    \label{p20230323_THz_Deer_S37G5_EE_Evolution_21}
    \end{subfigure}
    \begin{subfigure}[t]{0.45\textwidth}
    \vspace{-.1cm}
    \includegraphics[width=\textwidth]{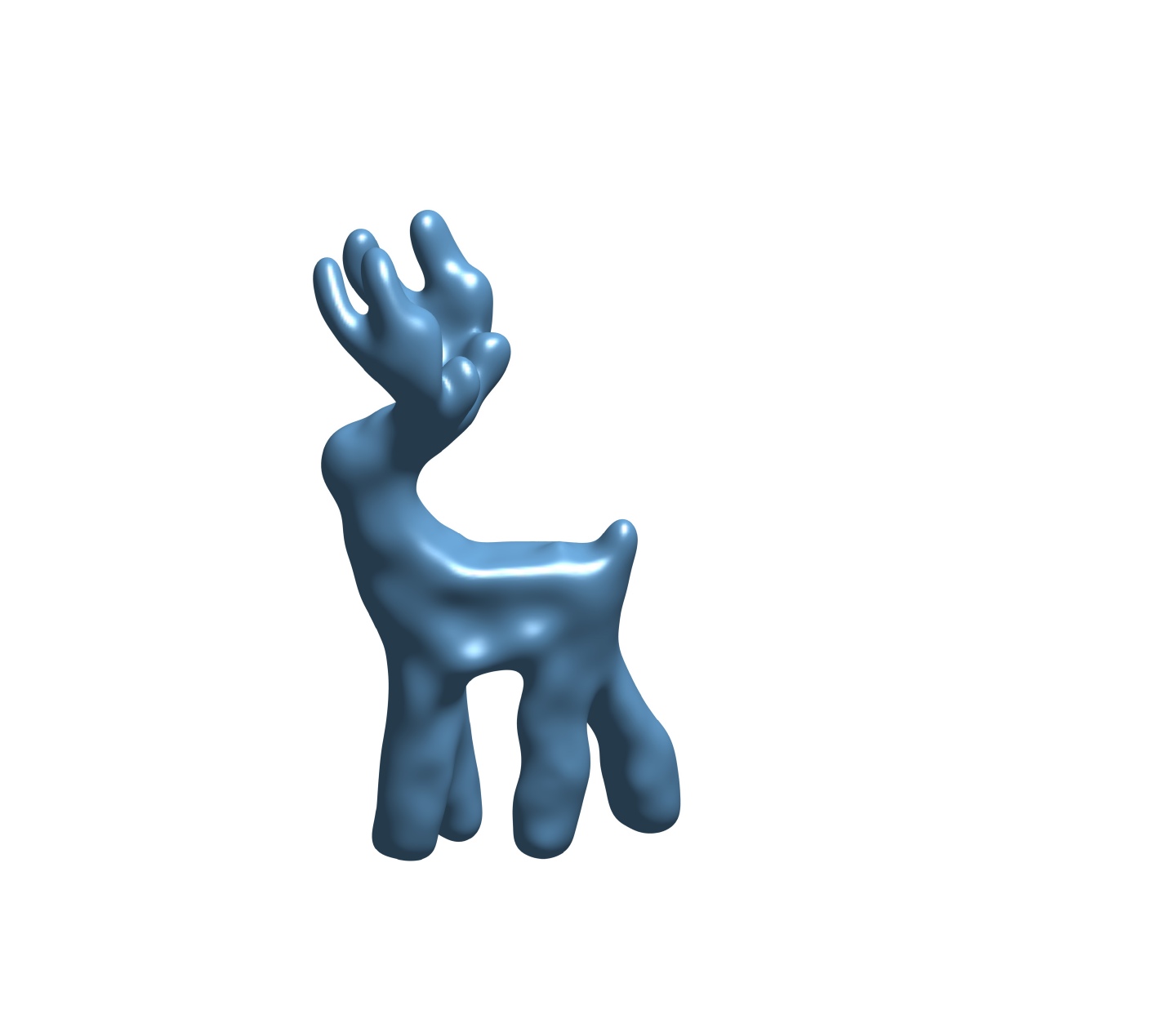}
    \vspace{-1cm}
    \caption{}
    \label{p20230323_THz_Deer_S37G5_EE_Evolution_24}
    \end{subfigure}
\cprotect \caption{Visualisation of {\bf Example 5} Deer from THz imaging:~\subref{p20230323_THz_Deer_S218G0_EE_Evolution_1} rough input with noise from the given all $218$ slices;~\subref{p20230323_THz_Deer_S218G0_EE_Evolution_4} smoothed reconstruction result by the Euler-Elastica-based formulation from~\subref{p20230323_THz_Deer_S218G0_EE_Evolution_1};~\subref{p20230323_THz_Deer_S109G1_EE_Evolution_5} input $109$ slices with $1$ gap;~\subref{p20230323_THz_Deer_S109G1_EE_Evolution_8} output from~\subref{p20230323_THz_Deer_S109G1_EE_Evolution_5};~\subref{p20230323_THz_Deer_S55G3_EE_Evolution_13} input $55$ slices with $3$ gaps;~\subref{p20230323_THz_Deer_S55G3_EE_Evolution_16} output from~\subref{p20230323_THz_Deer_S55G3_EE_Evolution_13};~\subref{p20230323_THz_Deer_S37G5_EE_Evolution_21} input $37$ slices with $5$ gaps;~\subref{p20230323_THz_Deer_S37G5_EE_Evolution_24} output from~\subref{p20230323_THz_Deer_S37G5_EE_Evolution_21}. }
\label{p20230323_THz_Deer_S218G0_EE_Evolution}
\end{figure}

We are now ready to present some experimental results in this section.
3D reconstructed results will be performed in the following using
\begin{itemize}
\item 3D surface inpainting models extended from~\cite{Carola15},
\item the perimeter-based formulation~\eqref{eqn:P_model},
\item the Willmore-based formulation~\eqref{eqn:W_model}, and
\item the new Euler-Elastica-based formulation~\eqref{eqn:EE_model},
\end{itemize}
where four inpainting models are tested as follows
\begin{itemize}
    \item the Cahn-Hilliard model
    (see Section \href{https://doi.org/10.1017/CBO9780511734304.006}{5.3} in~\cite{Carola15}),

    \item the Mumford-Shah model
    (see Chapter \href{https://doi.org/10.1017/CBO9780511734304.008}{7} in~\cite{Carola15}),

    \item the transport model
    (see Section \href{https://doi.org/10.1017/CBO9780511734304.007}{6.1} in~\cite{Carola15}), and

    \item the absolutely minimising Lipschitz extensions
    (see Section \href{https://doi.org/10.1017/CBO9780511734304.005}{4.4} in~\cite{Carola15}).
\end{itemize}
The experimental results consist of two simulated examples, namely ({\bf Example 1} Sphere like tumour-liked simulation and {\bf Example 2} Branching Cylinders as branching blood vessels mimicked), as well as two segmented realistic examples ({\bf Example 3} Stent segmented from real CT images and {\bf Example 4} Tumour segmented from real MRI images), and one realistic example ({\bf Example 5} Deer from THz imaging).
For {\bf Example 1} Sphere, we demonstrate the simulation of the input slices and the gap-filling process (\Cref{p20220215_S_N32S5_NoTITLE_GCMClim_PWEE_2-3}).
The results obtained using the compared models are presented in~\Cref{p20220715_S_N32S5_CH_MS_AMLE_T_H_inpainting_2-5} and~\Cref{p20220215_S_4_5_6}.
However, only the inpainting results are shown in this example as they are unsatisfactory in meeting our expectations, despite some of them being feasible for  gap filling (e.g.~\Cref{p20220715_S_N32S5_CH_MS_AMLE_T_H_inpainting_CH_2} and \subref{p20220715_S_N32S5_CH_MS_AMLE_T_H_inpainting_MS_3}).
In {\bf Example 2} Branching Cylinders, we illustrate the feasibility of concave geometrical morphology for our formulations.
The simulated input and its results using three formulations are depicted in~\Cref{p20220215_BC_2_3} and~\ref{p20220215_BC_4_5_6}.
To compare the results obtained from the three formulations and establish a benchmark as the stopping criterion, we propose an approach for {\bf quantitative comparisons}.
This is done through the mathematical and graphical interpretation (Equations~\eqref{eqn:GC}-\eqref{eqn:MC} and~\Cref{20220716_ViewRing}) from the perspective of discrete geometry.
The variance of the two simulated examples is visualised in~\Cref{p20220215_S_9_11_13_10_12_14} and~\ref{p20220215_BC_9_11_13_10_12_14} with histograms in~\Cref{p20230319_S_N32S5_Histogram_GCMClim_PWEE} and~\ref{p20230319_BC_N128S24_Histogram_GCMClim_PWEE}.
Furthermore, the numeric corroboration of the comparison between the three formulations and {\bf Examples 1-3} is presented in~\Cref{tab:std_GC_MC}, and the {\bf experimental convergence and computational complexity of PGDM for Example 1} is provided in~\Cref{p20230327_S_N32S5_RelativeError_PWEE_1} and~\ref{p20230327_S_N32S5_ComputationalComplexity_PWEE_1}.
Considering the gradient descent method used throughout ADMM, an {\bf experimental analysis of ADMM for Example 1} is provided in~\Cref{p20220703_S_N32S5_EE_ADMM_PT}.
This analysis focuses on the sensitivity of parameters to establish the relationship between the new formulation and parameters for faster and better numerical simulations.
Lastly, we present the results of three realistic examples, ({\bf Example 3} Stent segmented from real CT images, {\bf Example 4} Tumour segmented from real MRI images, and {\bf Example 5} Deer from THz imaging) in~\Cref{p20220302_RealData_Stent_N512S48_NoTITLE_EE_2_3_4_5}-\ref{p20230323_THz_Deer_S218G0_EE_Evolution}.
These results demonstrate the application of the new Euler-Elastica formulation and verify its merits.
Remark that all implementations were coded using the computer programming language:  MATLAB\_R2022a\textsuperscript{\tiny\textregistered} in the operating system: macOS Monterey (Version 12.5) equipped with a 2.3 GHz 8-Core Intel Core i9 Processor and 16GB 2667 MHz DDR4 Memory (some implementations in the revised version were carried out using MATLAB\_R2023a\textsuperscript{\tiny\textregistered} in the operating system: macOS Ventura (Version 13.3.1 (a)) equipped with an Apple M1 Max Chip and 64 GB Memory).


{\bf Example 1 (Sphere tested by all compared models)}.
To simulate the reconstructed problem from a few slices, by setting the low resolution $N = 32$, the rough $5$ slices are collected from a Sphere as tumour-liked simulation being the first example as~\Cref{p20220215_S_N32S5_NoTITLE_GCMClim_PWEE_2} illustrated.
Then, the initially rough surface can be straightforwardly constructed by duplicating the slices to fill the gaps as~\Cref{p20220215_S_N32S5_NoTITLE_GCMClim_PWEE_3} visualised.
Remark that the gaps between the slices range from four to five. To fill the gap between the top and bottom slices, we duplicate half of the top slice and half of the bottom slice. For the remaining slices, we use the slices themselves to fill the half-up and half-down gaps.

First, to compare the variational framework with phase-field approximation, we would like to demonstrate the results (see~\Cref{p20220715_S_N32S5_CH_MS_AMLE_T_H_inpainting_2-5}) by the explored extension of four 3D surface inpainting models (\subref{p20220715_S_N32S5_CH_MS_AMLE_T_H_inpainting_CH_2} the Cahn-Hilliard model,
~\subref{p20220715_S_N32S5_CH_MS_AMLE_T_H_inpainting_MS_3} the Mumford-Shah model,
~\subref{p20220715_S_N32S5_CH_MS_AMLE_T_H_inpainting_T_5} the transport model, and~\subref{p20220715_S_N32S5_CH_MS_AMLE_T_H_inpainting_AMLE_4} the absolute minimising Lipschitz extensions) which are introduced from~\cite{Carola15}.
The extension from 2D image inpainting to 3D surface inpainting is explored by the case of missing slices in the vertical direction to restore one of the planes in two axes other than the vertical axis.
Apparently, such a task is quite challenging for the above inpainting models, which these illustrated results are visually unacceptable with the time-consuming acquisition, and probably the sunken gaps would lead to unexpected variations even though the gaps of two of the results were filled.
The main reason for unsatisfactory results by inpainting models is due to their local and partial inpainting without considering the global variation.
Remark that the stopping criteria for each inpainting model are set as respectively:
the maximum number $500$ of iterations for the Cahn-Hilliard model;
the residual less than the tolerance $10^{-14}$ for the Mumford-Shah model with the maximum number $50$ of iterations;
the residual less than the tolerance $10^{-5}$ for the transport model with the maximum number $50$ of iterations; and
the residual less than the tolerance $10^{-8}$ for the absolutely minimising Lipschitz extensions with the maximum number $50$ of iterations.

Thereupon, for our variational framework with phase-field approximation, the parameters are set as $\varepsilon = 1.5/N$, $\tau = \varepsilon^{4}$, and the stopping criterion by the difference of corresponding energies between new iterative results and previous one less than the preset value.
By the convergence results of~\eqref{equ:UconvergenceResult}, the surface by three formulations (perimeter-based, Willmore-based, and Euler-Elastica-based formulation) can be extracted from the iterative results under the isosurface value at half as~\Cref{p20220215_S_4_5_6} demonstrated.
Under observation of~\Cref{p20220215_S_4_5_6}, the last one by investing the Euler-Elastica-based formulation is the desired result.
Noted that applying the Willmore-based formulation produces a smoother surface with undesired shapes of the top and bottom due to achieving lower energy of mean curvature with the property of Willmore energy.
Moreover, applying the perimeter-based formulation emerges small bulges patently, whereas it maintains the initial shape globally.

 \vspace{-.2cm}

{\bf Example 2 (Branching Cylinders tested by three formulations)}.
Analogously, for the second example: Branching Cylinders as branching blood vessels mimicked,~\Cref{p20220215_BC_2_3} delineates the initially rough surface under the same duplicating idea in the first example, which is constructed from the given $24$ slices under the low resolution $N = 128$.
Remark that the gap of collected slices is intentionally uneven so that the initial surface has a distinct fluctuation waiting to be restored, and the gap range is from three to thirteen.

As the resolution of this example is increased, keeping the same diffuse interface width $\varepsilon = 1.5/N$ and setting the larger time step $\tau = \{\varepsilon^{3}, 10 \varepsilon^{4}\}$ for the faster stable results where $\varepsilon^{3}$ for perimeter-based and Willmore-based formulation, and the latter one for Euler-Elastica-based formulation by the corresponding fixed point iterative schemes, then the reconstructed surfaces by three formulations are extracted as~\Cref{p20220215_BC_4_5_6} performed.
Following the presupposition in~\Cref{sec:Intro}, as the Euler-Elastica-based formulation merges the advantages of the perimeter-based formulation and the Willmore-based formulation to overcome the above deficiencies, the surface by the Euler-Elastica-based formulation is heralded the better reconstruction comparing by the Willmore-based formulation in this scenario either the top plane or the concave between the branching part.

{\bf Quantitative Comparisons}.
Broadly speaking, the reconstruction of Branching Cylinders using the Willmore-based formulation (\Cref{p20220215_BC_N128S24_NoTITLE_GCMClim_PWEE_5}) initially fulfilled the initial requirements, exhibiting minor deficiencies that were difficult to discern through visual inspection alone.
However, when compared to the results obtained using the newly proposed Euler-Elastica-based formulation, these deficiencies become more noticeable.

Accordingly, to indicate the level of smoothness for the above results by three formulations, the quantitative benchmark is considered by computing its standard deviation of Gaussian curvatures (GC) $\sigma_{\mbox{\tiny{GC}}}:=\sigma(\kappa_{G})$ and of mean curvatures (MC) $\sigma_{\mbox{\tiny{MC}}}:=\sigma(\bar{\kappa})$ for surfaces represented by triangular meshes from the viewpoint of discrete geometry.
To recap from the theoretical discrete geometry, the Gaussian curvature of each vertex is given by
\begin{eqnarray}
\label{eqn:GC}
\kappa_{G}(\mathbf{v}_{i})
=
\frac{2\pi - \sum\limits_{k=1}^{\N_{F_{i_{R}}}} \theta_{i_{k}}}{\mathcal{A}_{i_{R}}}
\end{eqnarray}
and the mean curvature of each vertex is given by
\begin{eqnarray}
\label{eqn:MC}
\bar{\kappa}(\mathbf{v}_{i})
=
\frac{1}{2\mathcal{A}_{i_{R}}}
\sum\limits_{j \in _{R}NV_{i}} (\cot{\alpha_{ij}} + \cot{\beta_{ij}}) (\mathbf{v}_{j} - \mathbf{v}_{i})
\end{eqnarray}
where $\mathcal{A}_{i_{R}}$ stands for the appropriately chosen area from the patch within $R$-ring neighbouring vertices $_{R}NV_{i}$ (i.e. the minimum number of edges from $\mathbf{v}_{i}$ to the neighbouring vertex is less than or equal to $R$ where, in this case, $R$ is opted for 1) around the vertex $\mathbf{v}_{i}$, $\theta_{i_{k}}$ denotes the angle of the $k^{\mbox{th}}$ face at the vertex $\mathbf{v}_{i}$, $\N_{F_{i_{R}}}$ is the total number of faces in the set $F_{i_{R}}$ around this vertex $\mathbf{v}_{i}$, as well as $\alpha_{ij}$ and $\beta_{ij}$ are two angles opposite to the sharing edge in the two triangles as~\Cref{20220716_ViewRing} exemplified, which can be consulted~\cite{Discrete_Differential-Geometry_Operators_for_Triangulated_2-Manifolds, 20201014_A_novel_method_for_surface_mesh_smoothing:_Applications_in_biomedical_modeling} for more details.
Remark that, by the concept of  geometrical measure, the Gaussian curvature is intrinsically invariant and relies only on surface-estimated distances, yet the mean curvature embedded surface locally is extrinsic evaluation in some ambient space e.g. Euclidean space.

In accordance with~\Cref{p20220215_S_4_5_6} and~~\Cref{p20220215_BC_4_5_6} for {\bf Examples 1} and {\bf 2}, we present in~\Cref{p20220215_S_9_11_13_10_12_14} and~\Cref{p20220215_BC_9_11_13_10_12_14} the variance of Gaussian curvature and mean curvature at each vertex of the triangular meshes, computed using three different formulations for these two examples.
In order to facilitate the comparison of the results, we set the minimum (maximum) of the colorbar to the maximum (minimum) curvature of all vertices for the three formulations and meshes: $\max{(\min{(\mathscr{P})}, \min{(\mathscr{W})}, \min{(\mathscr{E})})}$ and $\min{(\max{(\mathscr{P})}, \max{(\mathscr{W})}, \max{(\mathscr{E})})}$, respectively.
We also provide histograms in~\Cref{p20230319_S_N32S5_Histogram_GCMClim_PWEE} and~\Cref{p20230319_BC_N128S24_Histogram_GCMClim_PWEE} to show the proportional distribution of all curvatures at each vertex and to compare the curvature distributions across different meshes using three formulations.
We note that, as expected for \textbf{Example 1} and \textbf{2} by $\mathscr{E}$, the majority of the curvatures cluster around zero.


\vspace{-.7cm}

\vspace{-.8cm}

\vspace{-.3cm}

In~\Cref{tab:std_GC_MC}, we present numerical comparisons of the three formulations for {\bf Examples 1-3}, in terms of the standard deviation of Gaussian curvatures ($\sigma_{\mbox{\tiny{GC}}}$), the standard deviation of mean curvatures ($\sigma_{\mbox{\tiny{MC}}}$), and CPU elapsed time.
We observe that, for the same input, the trend of values of $\sigma_{\mbox{\tiny{GC}}}$ and $\sigma_{\mbox{\tiny{MC}}}$ is decreasing, indicating that the level of smoothness is improving as the amount of variation is reduced, thus indicating a better reconstruction by the new proposed formulation.
Additionally, we report the average elapsed time of each iteration by the three formulations, highlighting the computational efficiency of our numerical approach.

\begin{table}[htbp]
\cprotect \caption{Comparisons of the standard deviation of Gaussian curvatures $\sigma_{\mbox{\tiny{GC}}}$ and of mean curvatures $\sigma_{\mbox{\tiny{MC}}}$, and of the average elapsed time of each iteration (seconds/iteration or s/iter) by three formulations: perimeter-based ($\mathscr{P}$), Willmore-based ($\mathscr{W}$), and Euler-Elastica-based ($\mathscr{E}$) formulation to three examples: {\bf Example 1} (\Cref{p20220215_S_9_11_13_10_12_14}), {\bf Example 2} (\Cref{p20220215_BC_9_11_13_10_12_14}) and {\bf Example 3} (\Cref{p20220302_RealData_Stent_N512S48_NoTITLE_EE_2_3_4_5}). }
\vspace{-0.3cm}
\label{tab:std_GC_MC}
\begin{center}
\begin{tabular}{cc|ccc}
  \toprule
  & &
  Model
  $\mathscr{P}$
  & Model
  $\mathscr{W}$
  & New model
  $\mathscr{E}$
  \\
  \hline
  \begin{tabular}{c}
 {\bf Example 1} \\ Sphere \\ ($N = 32$)
  \end{tabular}
  &
  \begin{tabular}{c}
   $\sigma_{\mbox{\tiny{GC}}}$\\ $\sigma_{\mbox{\tiny{MC}}}$\\  s/iter
  \end{tabular}
  &
  \begin{tabular}{c}
   85.6522\\ 4.9897\\ 0.0009
  \end{tabular}
  &
  \begin{tabular}{c}
 10.5005\\ 1.2315\\ 0.0030
  \end{tabular}
  &
  \begin{tabular}{c}
\textbf{4.3289}\\ \textbf{0.9789}\\ 0.0033
  \end{tabular}
  \\
  \hline
  \begin{tabular}{c}
 {\bf Example 2} \\ Branching Cylinders \\ ($N = 128$)
  \end{tabular}
  &
  \begin{tabular}{c}
    $\sigma_{\mbox{\tiny{GC}}}$\\ $\sigma_{\mbox{\tiny{MC}}}$\\  s/iter
  \end{tabular}
  &
  \begin{tabular}{c}
 772.2918\\ 10.9077\\ 0.0835
  \end{tabular}
  &
  \begin{tabular}{c}
 191.8314\\ 8.9567\\ 0.2791
  \end{tabular}
  &
  \begin{tabular}{c}
 \textbf{89.1912}\\ \textbf{6.3708}\\ 0.3108
  \end{tabular}
  \\
  \hline
   \begin{tabular}{c}
 {\bf Example 3} \\ Stent \\ ($N = 512$)
  \end{tabular}
  &
  \begin{tabular}{c}
    $\sigma_{\mbox{\tiny{GC}}}$\\ $\sigma_{\mbox{\tiny{MC}}}$\\  s/iter
  \end{tabular}
  &
  \begin{tabular}{c}
 16372.1718 \\ 100.5420 \\ 0.4761
  \end{tabular}
  &
  \begin{tabular}{c}
 1734.3918 \\ 35.2987 \\ 1.4679
  \end{tabular}
  &
  \begin{tabular}{c}
 \textbf{1628.0449} \\ \textbf{33.3259} \\ 1.4778
  \end{tabular}
  \\
  \bottomrule
\end{tabular}
\end{center}
\end{table}


{\bf Experimental Convergence and Computational Complexity of PGDM for Example 1}.
To evaluate the convergence of the numerical algorithm I, in~\Cref{p20230327_S_N32S5_RelativeError_PWEE_1}, we plot the relative error between the current and previous iterations over the number of iterations for {\bf Example 1} with respect to three formulations.
From the convergence plot, we observe that the algorithm with Euler-Elastica-based formulation converges rapidly within the first $100$ iterations and reaches a relative error of $10^{-4}$ after approximately $300$ iterations.

Additionally, we evaluate the computational complexity of numerical algorithm I in terms of the number of iterations and the time required to run the algorithm.
To estimate the computational complexity of an algorithm, we analyse the number of operations or steps the algorithm takes as the size $N$ of the input increases.
In our 3D scenario, the main operations inside the loop are the fast Fourier transform (FFT) and its inverse, which have a complexity of $\mathcal{O}(N\log{N})$ for each axis.
Besides, the Laplacian operator has a complexity of $\mathcal{O}(N^{3})$, while the other operations have a lower complexity.
Therefore, the computational complexity for each iteration of the algorithm I can be estimated as $\mathcal{O}(N^{3}\log{N})$.
Here, we measured the running time for {\bf Example 1} with various inputs $N$ on the same hardware and software environment and used the unit of time set as the arbitrary unit (a.u.) in~\Cref{p20230327_S_N32S5_ComputationalComplexity_PWEE_1}.
Then, we observe that the trend of experimental results and estimated arithmetical values are semblable.
Note that the actual running time of an algorithm depends not only on its computational complexity but also on the specific hardware and software environment in which it is executed.

{\bf Experimental Analysis of ADMM for Example 1}.
For the experimental results by the numerical algorithm II, the results are similar to the key computing gradient descent method used throughout ADMM.
Meanwhile, the result by this algorithm is more sensitive than by the first~\Cref{alg:Projected_Gradient_Descent_Method_GDM}, which means the result is strongly influenced by the related parameters (the penalty parameter $\rho$, the diffuse interface width $\varepsilon$, and the time step $\tau$), even though the speed of acquiring expected results is faster associated with fewer iterations.


Due to the sensitivity of parameters,~\Cref{p20220703_S_N32S5_EE_ADMM_PT} provides the binary maps of testing in relation to the penalty parameter $\rho \in [0.5, 10]$ with the step of $0.5$, and the representative diffuse interface width $\varepsilon \cdot N \in [1.5, 3]$ with the step of $0.1$ for {\bf Example 1} Sphere under the low resolution $N = 32$ in two time steps $\tau = \varepsilon^{3}$ and $\tau = \varepsilon^{3.5}$.
After setting the criterion by our proposed benchmark in the case of the standard deviation of Gaussian curvatures $\sigma_{\mbox{\tiny{GC}}}$ less than the value by the Willmore-based formulation which is $\sigma_{\mbox{\tiny{GC}}} < 10.5005$ in~\Cref{tab:std_GC_MC}, then the binary maps can be ascertained the suitable range of parameters for the reliable results where the binary value $1$ stands for the tolerable results existed, and $0$ indicates the unpleasant results during the iterations.

\begin{figure}[htbp]
    \centering
    \includegraphics[width=0.7\textwidth]{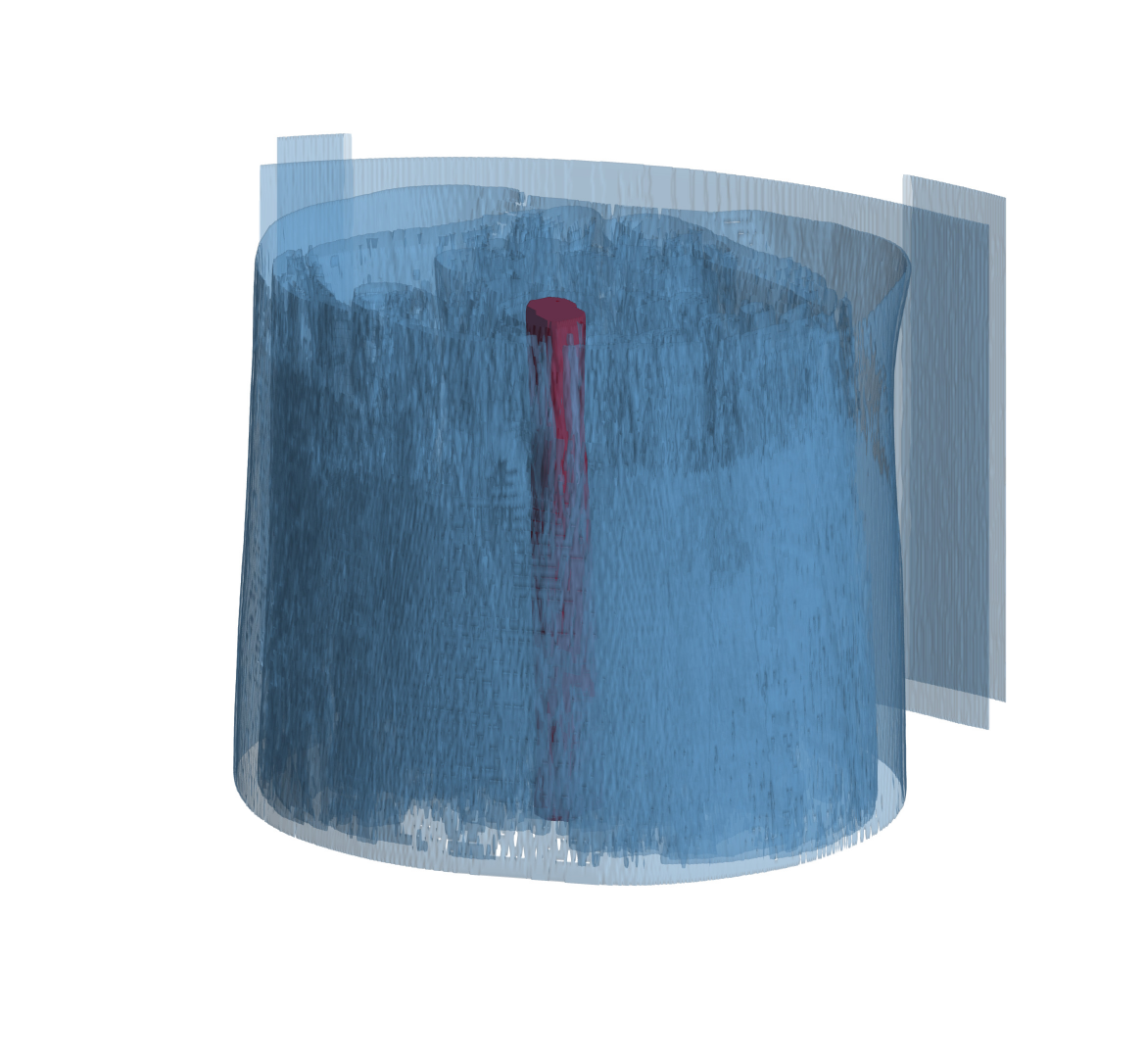}
\vspace{-0.5cm}
\cprotect \caption{Visualisation of initial CT real data of the Chest (blue) under the low resolution $N = 512$ from $48$ slices where {\bf Example 3} segmented Stent is indicated as red. }
\label{p20220302_RealData_Stent_N512S48_NoTITLE_EE_1}
\vspace{-0.3cm}
\end{figure}

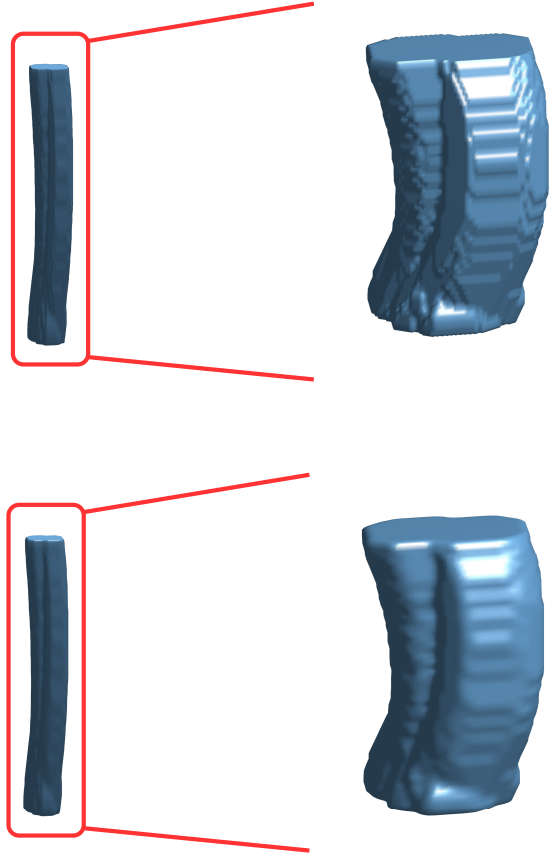
\begin{figure}[htbp]
    \centering
    \begin{subfigure}[t]{0.6\textwidth}
    \begin{tikzpicture}
    \hspace{-3cm}
    \node at (-3, 0)
    {\includegraphics[width=\textwidth]
    {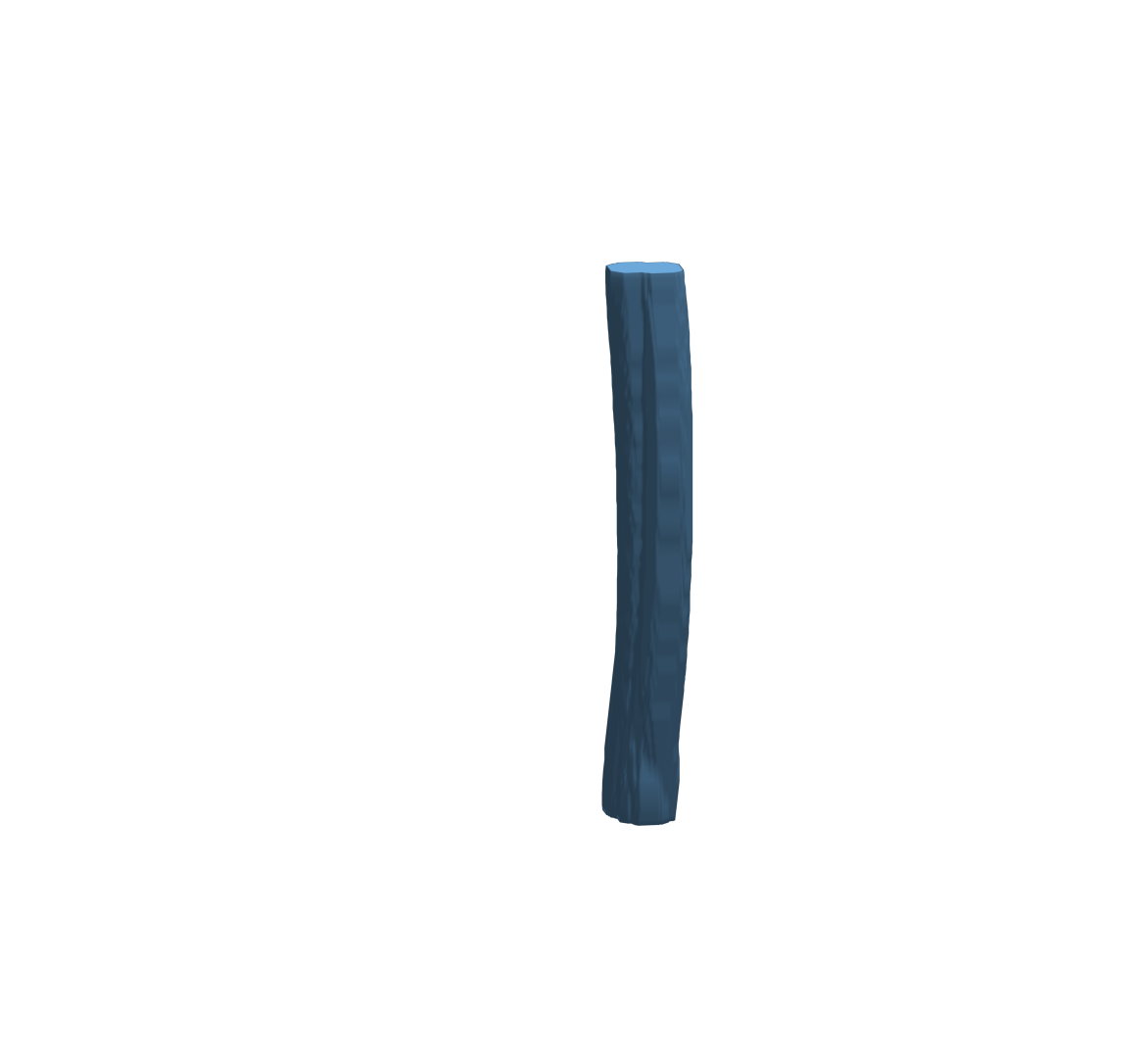}
    };

    \node at (3, 0) {\includegraphics[width=\textwidth]
    {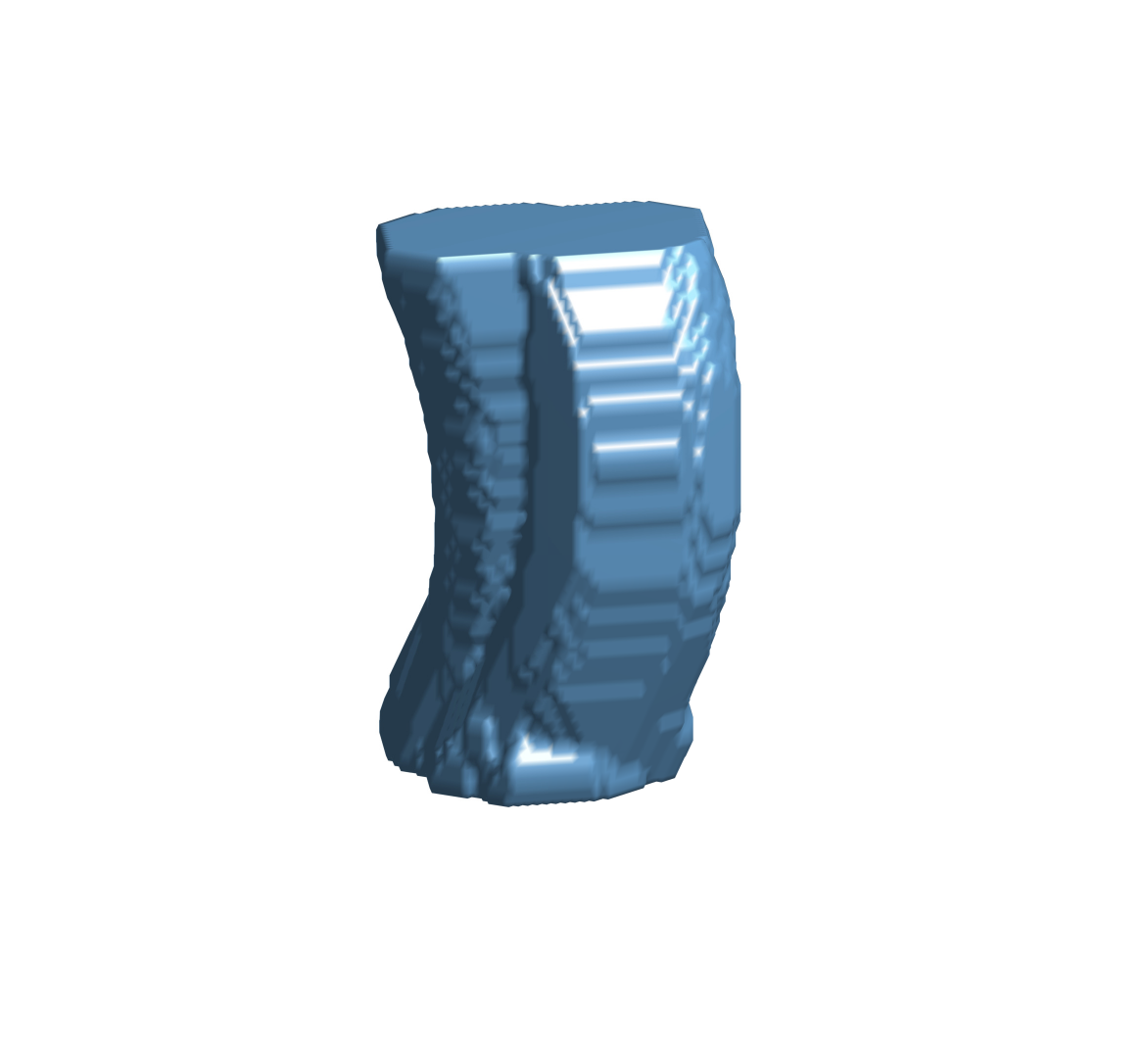}
    };

    \draw[red!80, ultra thick, rounded corners]
    (-3,-2.2) rectangle (-2,2.2);
    \draw[red!80, ultra thick]
    (-2.02,2.1) -- (1,2.6);
    \draw[red!80, ultra thick]
    (-2.02,-2.1) -- (1,-2.4);
    \end{tikzpicture}
        \label{p20220302_RealData_Stent_N512S48_NoTITLE_EE_2-3}
    \end{subfigure}
\\
\vspace{-1.5cm}
    \begin{subfigure}[t]{0.6\textwidth}
    \begin{tikzpicture}
    \hspace{-3cm}
    \node at (-3, 0)
    {\includegraphics[width=\textwidth]
    {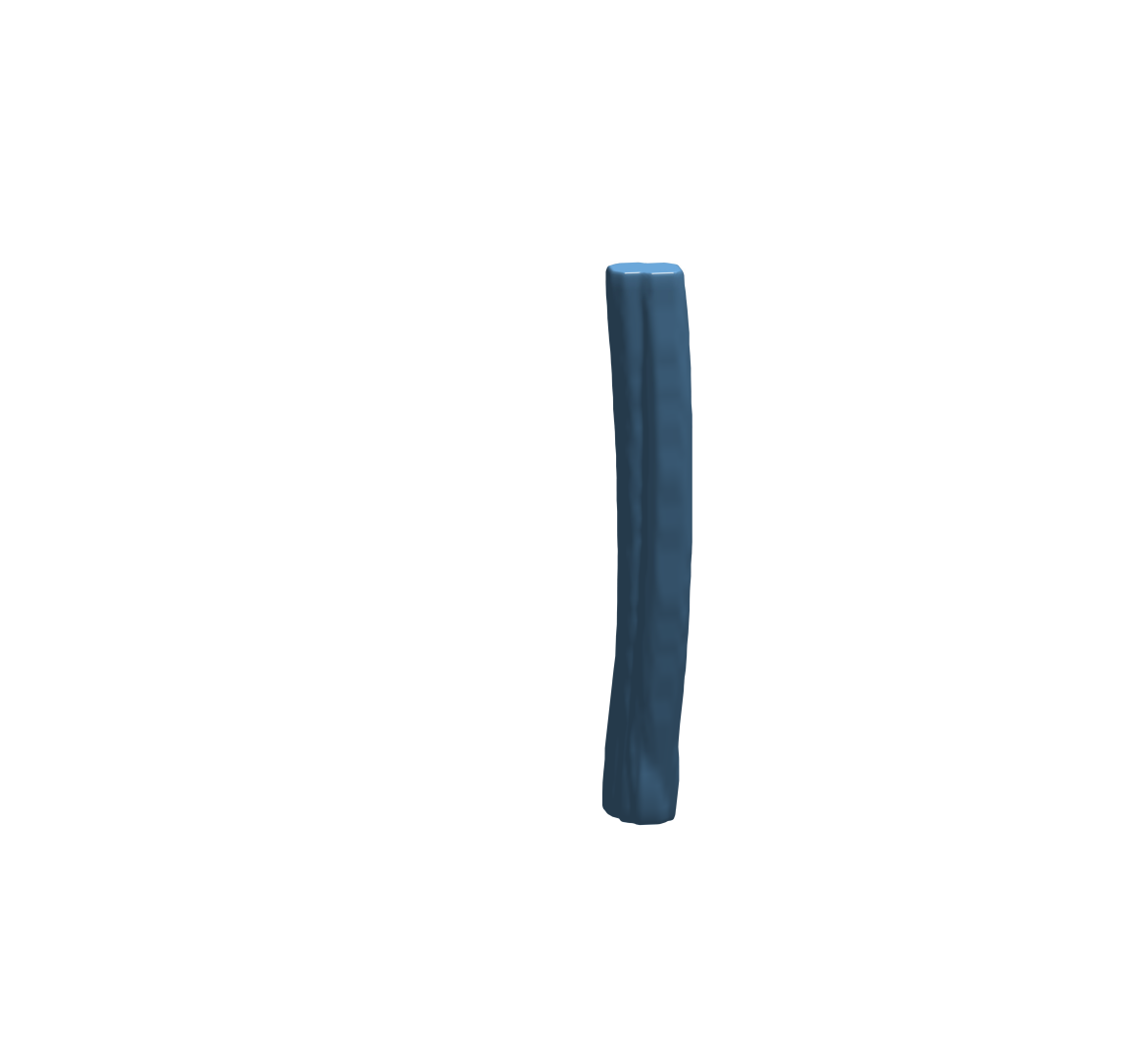}
    };

    \node at (3, 0) {\includegraphics[width=\textwidth]
    {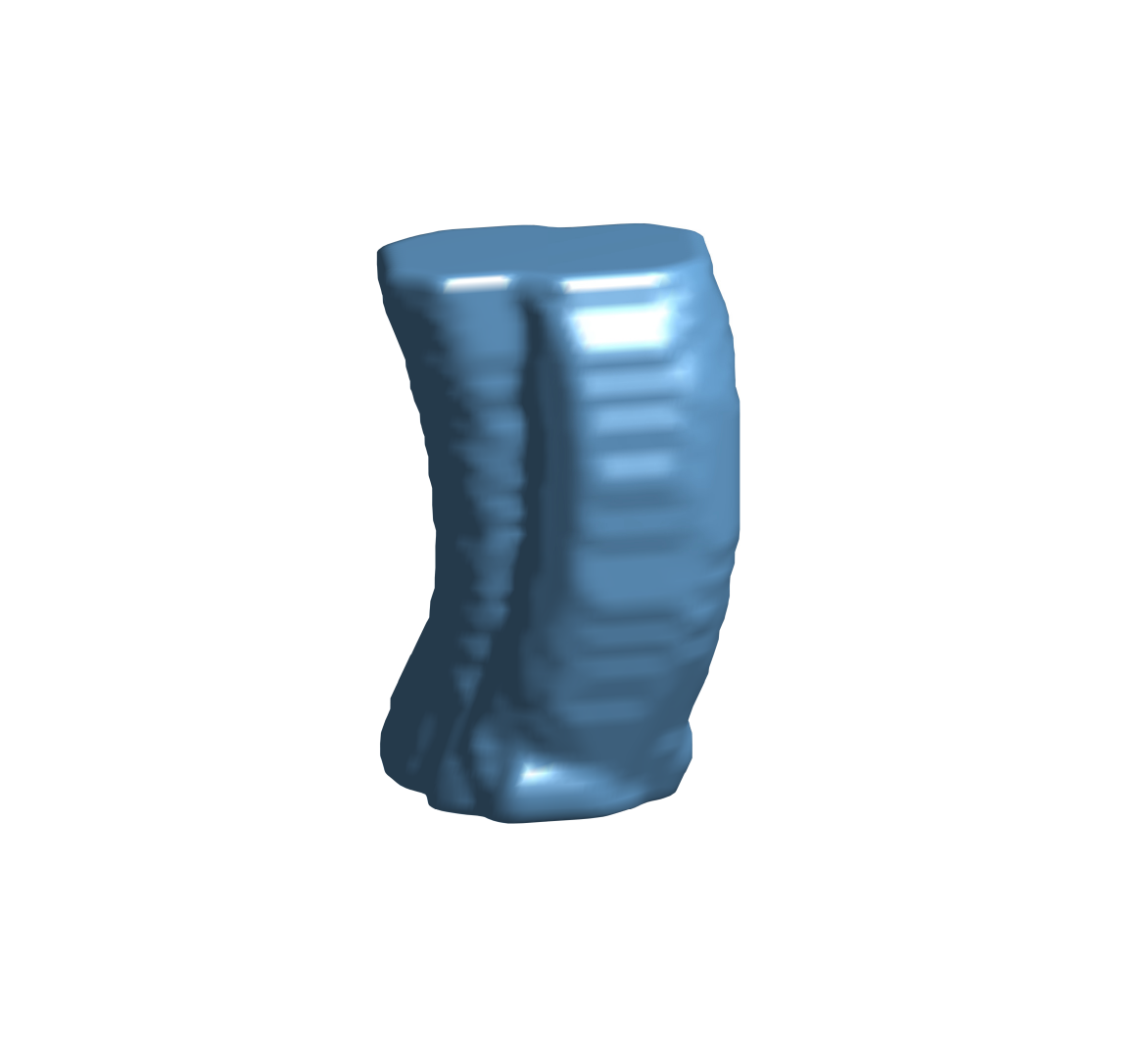}
    };

    \draw[red!80, ultra thick, rounded corners]
    (-3,-2.2) rectangle (-2,2.2);
    \draw[red!80, ultra thick]
    (-2.02,2.1) -- (1,2.6);
    \draw[red!80, ultra thick]
    (-2.02,-2.1) -- (1,-2.4);
    \end{tikzpicture}
    \label{p20220302_RealData_Stent_N512S48_NoTITLE_EE_4-5}

    \end{subfigure}
\vspace{-.8cm}
\cprotect \caption{Visualisation of segmented rough input (top) of {\bf Example 3} Stent from the given $48$ slices of CT real data under the low resolution (\Cref{p20220302_RealData_Stent_N512S48_NoTITLE_EE_1}) and the smoothed reconstruction results by the Euler-Elastica-based formulation (bottom) where the results on the right column are the enlarged view of the left column. }
\label{p20220302_RealData_Stent_N512S48_NoTITLE_EE_2_3_4_5}
\end{figure}

\begin{figure}[htbp]
    \centering
    \begin{subfigure}[t]{0.6\textwidth}
    \begin{tikzpicture}
    \hspace{-3cm}
    \node at (-3, 0)
    {\includegraphics[width=\textwidth]
    {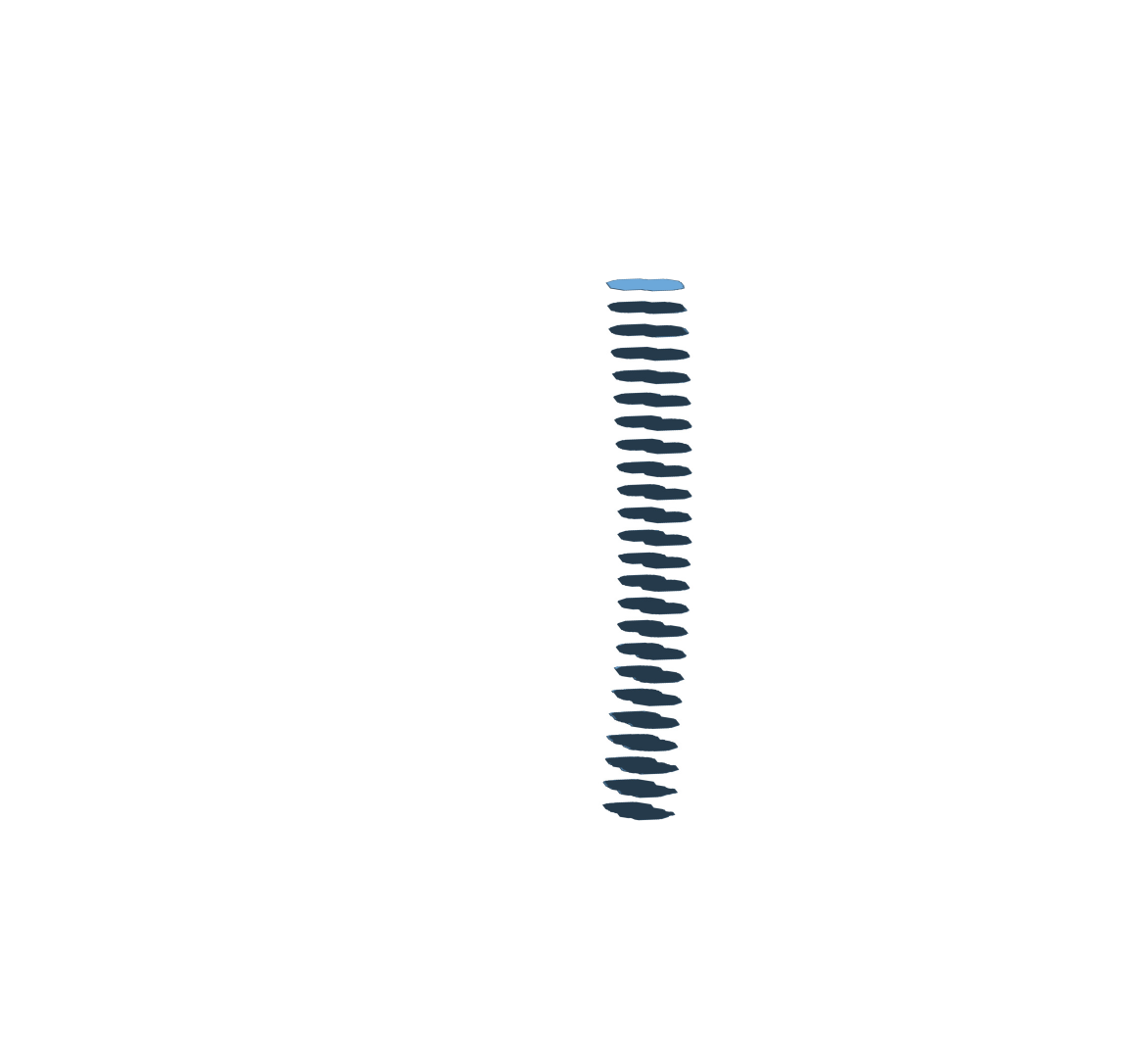}
    };

    \node at (3, 0) {\includegraphics[width=\textwidth]
    {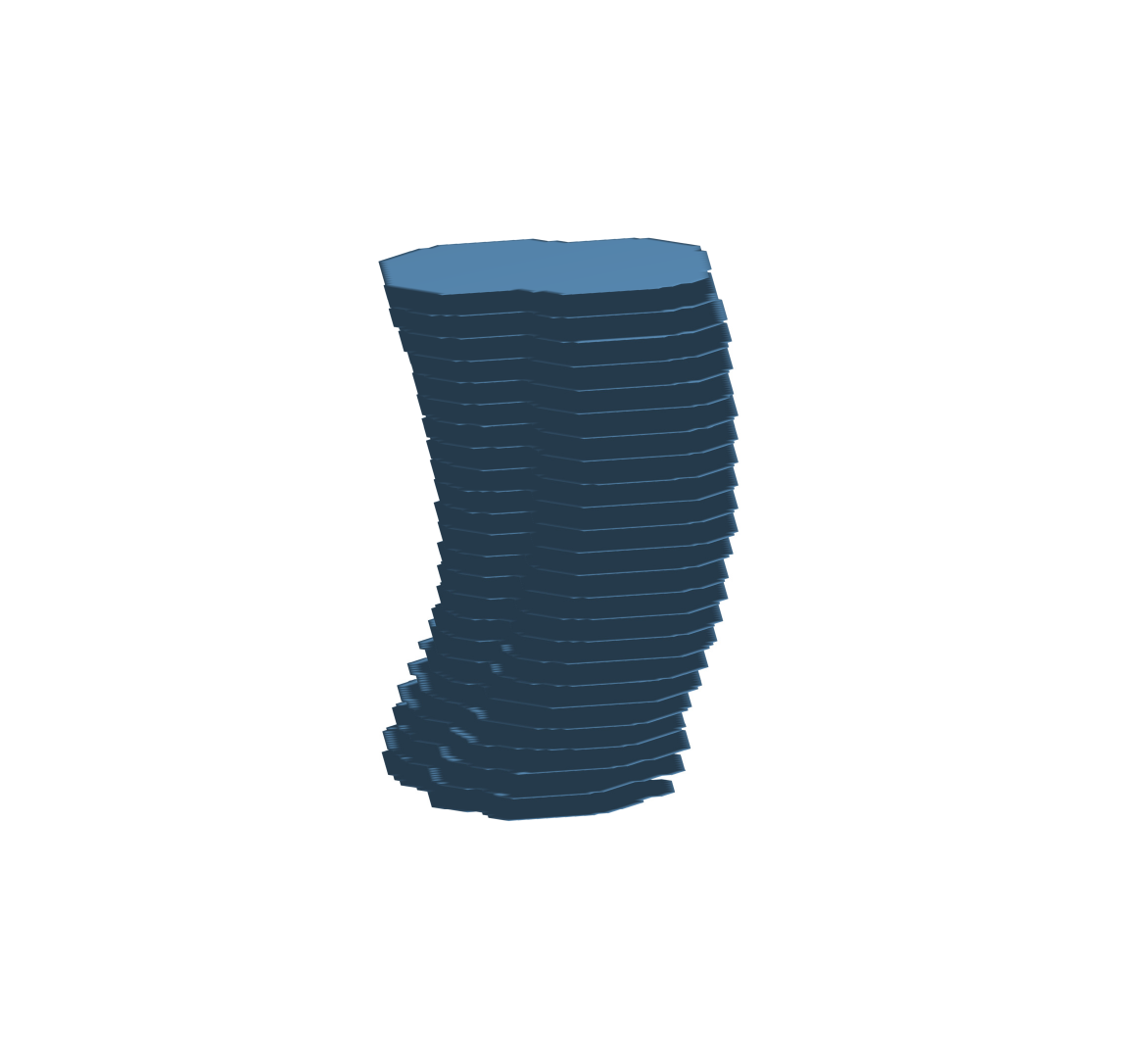}
    };

    \draw[red!80, ultra thick, rounded corners]
    (-3,-2.2) rectangle (-2,2.2);
    \draw[red!80, ultra thick]
    (-2.02,2.1) -- (1,2.6);
    \draw[red!80, ultra thick]
    (-2.02,-2.1) -- (1,-2.4);
    \end{tikzpicture}
    \label{p20220302_RealData_Stent_N512S24_NoTITLE_EE_1-2}
    \end{subfigure}
\\
\vspace{-1.5cm}
    \begin{subfigure}[t]{0.6\textwidth}
    \begin{tikzpicture}
    \hspace{-3cm}
    \node at (-3, 0)
    {\includegraphics[width=\textwidth]
    {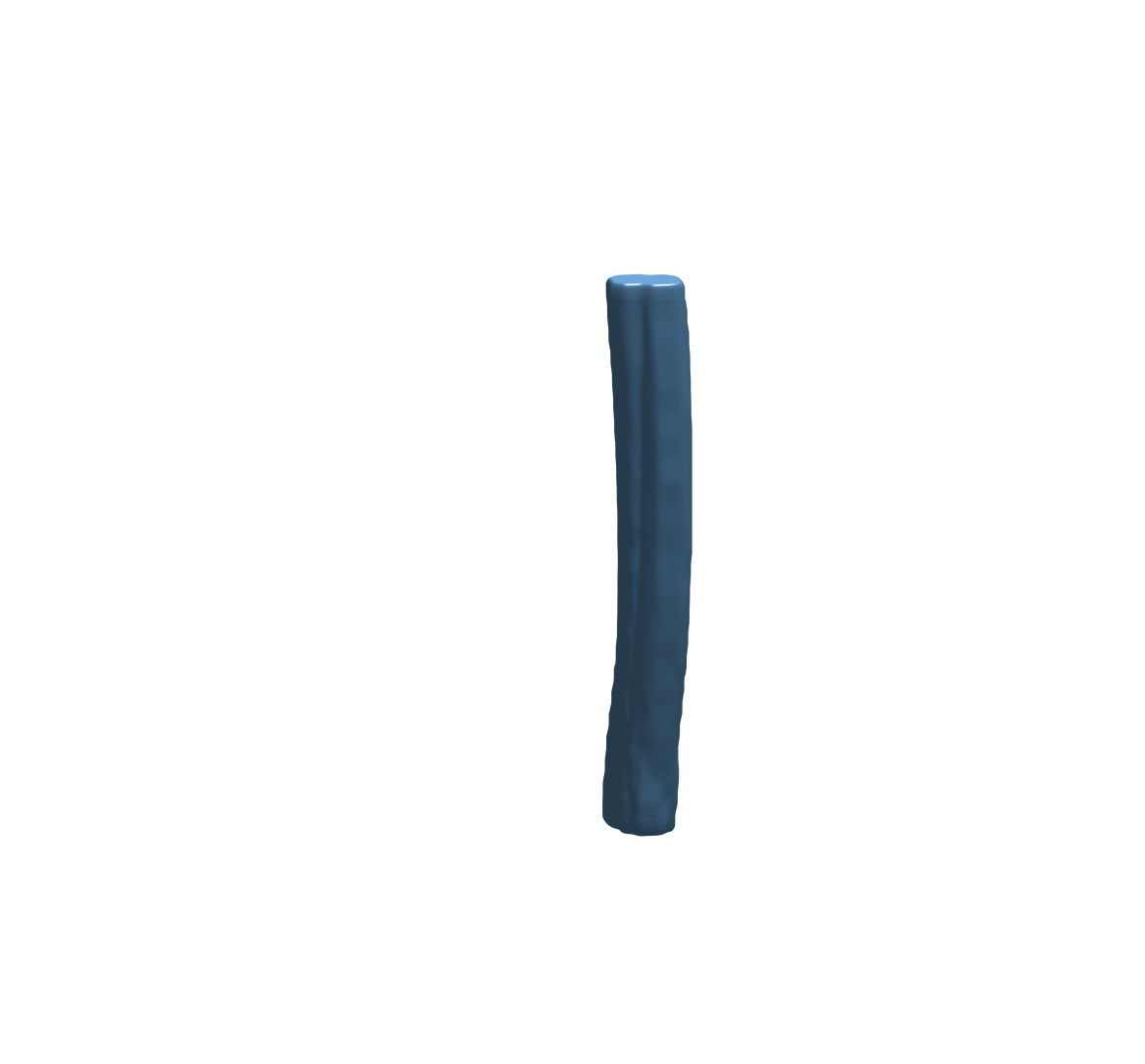}
    };

    \node at (3, 0) {\includegraphics[width=\textwidth]
    {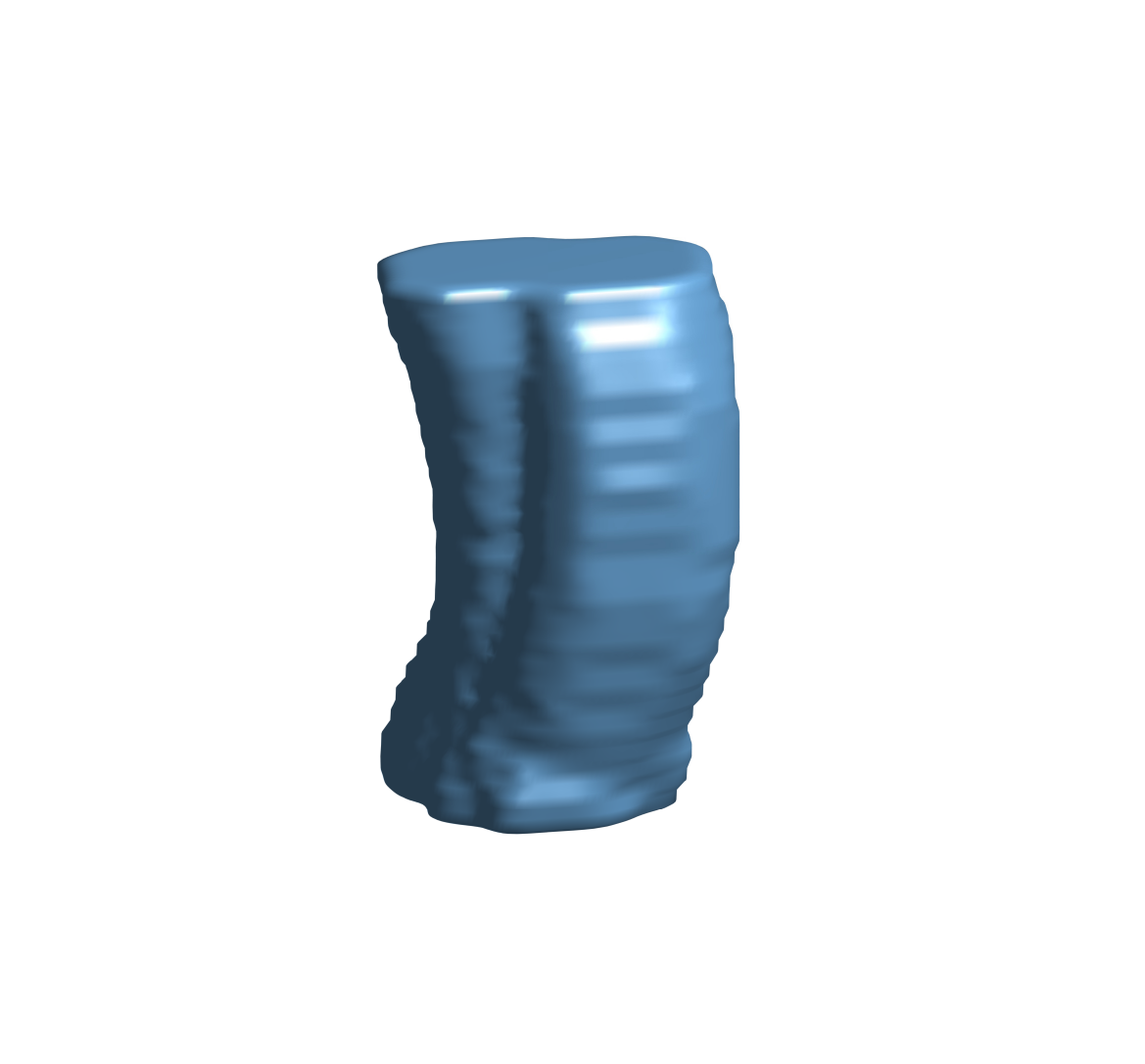}
    };

    \draw[red!80, ultra thick, rounded corners]
    (-3,-2.2) rectangle (-2,2.2);
    \draw[red!80, ultra thick]
    (-2.02,2.1) -- (1,2.6);
    \draw[red!80, ultra thick]
    (-2.02,-2.1) -- (1,-2.4);
    \end{tikzpicture}
        \label{p20220302_RealData_Stent_N512S24_NoTITLE_EE_5-6}
    \end{subfigure}
\vspace{-.5cm}
\cprotect \caption{Visualisation of segmented rough input slices (top) of {\bf Example 3} Stent by less $24$ slices from CT real data under the low resolution (\Cref{p20220302_RealData_Stent_N512S48_NoTITLE_EE_1}) and the smoothed reconstruction results by the Euler-Elastica-based formulation (bottom) where the results on the right column are the enlarged view of the left column. }
\label{p20220302_RealData_Stent_N512S24_NoTITLE_EE_1_2_5_6}
\vspace{-0.3cm}
\end{figure}

{\bf Example 3 (Stent segmented from real CT images)}.
Last but not least, the realistic examples are procured from the variational segmentation work by Dr Liam Burrows~\cite{Reproducible_kernel_Hilbert_space_based_global_and_local_image_segmentation} for the application of 3D reconstruction in medical imaging from 2D X-ray computed tomography (CT) scans and magnetic resonance imaging (MRI).
\Cref{p20220302_RealData_Stent_N512S48_NoTITLE_EE_1} visualises the direct construction by stacking all $48$ 2D low-resolution CT slices ($N = 512$) of the Chest.
Then, the segmented Stent is constructed in the top line of~\Cref{p20220302_RealData_Stent_N512S48_NoTITLE_EE_2_3_4_5} and following the smoothed results by the Euler-Elastica-based formulation.
Returning to the original objective of surface reconstruction from a reduced number of slices,~\Cref{p20220302_RealData_Stent_N512S24_NoTITLE_EE_1_2_5_6} demonstrates the efficacy of the new Euler-Elastica-based formulation by reconstructing the surface from only half of the available slices.
This approach not only reduces the time required for data collection in clinical imaging, thereby minimising patient exposure to uncontrollable high-dose radiation, but also improves the quality of reconstructed objects by incorporating super-resolution techniques.

{\bf Example 4 (Tumour segmented from real MRI images)}.
In this example, we demonstrate the efficacy of our Euler-Elastica-based formulation in segmenting tumours from high-resolution MRI images.
Specifically, we focus on a small region of interest in a large collection of 280 MRI images, resizing the area of interest to $150 \times 150$ pixels from the original size of $1210 \times 2378$.
This scenario poses a subtle challenge, as the small size of the region of interest makes it difficult to accurately segment the tumour from the surrounding tissue.
However, by applying our new proposed model, we are able to achieve highly accurate results, as shown in~\Cref{p20220309_RealData_Cancer_EE_1-2}.
Remark that we use the similar parameter settings as in {\bf Example 3}, with the same time step of $\tau_t = \varepsilon^4$ and $\varepsilon = 2/N$ for the stent and $\varepsilon = 2.5/\bar{N}$ with the average pixel number $\bar{N}$ of three axes for the tumour.

{\bf Example 5 (Deer from real THz imaging)}.
Terahertz (THz) imaging has the potential to revolutionise medical imaging due to its non-ionising nature and ability to penetrate through certain materials.
However, the images produced by THz imaging often suffer from low resolution with noise, and require significant time for acquisition~\cite{Yiyao_2022}.
The proposed Euler-Elastica-based formulation is also applicable for the 3D reconstruction from THz imaging data, as demonstrated by the example of a deer in~\Cref{p20230323_THz_Deer_S218G0_EE_Evolution}.
By using the proposed formulation, we are able to effectively address the challenges posed by THz imaging and reconstruct smooth 3D models of the deer.
The results, shown in~\Cref{p20230323_THz_Deer_S218G0_EE_Evolution}~\subref{p20230323_THz_Deer_S218G0_EE_Evolution_4},~\subref{p20230323_THz_Deer_S109G1_EE_Evolution_8},~\subref{p20230323_THz_Deer_S55G3_EE_Evolution_16},~\subref{p20230323_THz_Deer_S37G5_EE_Evolution_24}, illustrate the successful reconstruction from the full input of 218 slices~\subref{p20230323_THz_Deer_S218G0_EE_Evolution_1} and fewer inputs~\subref{p20230323_THz_Deer_S109G1_EE_Evolution_5},~\subref{p20230323_THz_Deer_S55G3_EE_Evolution_13},~\subref{p20230323_THz_Deer_S37G5_EE_Evolution_21}.
The parameter settings used for this example are $\tau = \varepsilon^{3.5}$ and $\varepsilon = 3/\max{(N_{x}, N_{y}, N_{z})}$ with the maximum pixel number of three axes.
These results demonstrate the potential of the proposed formulation for improving the quality and speeding up of 3D reconstructions from THz imaging data with low resolution and fewer slices, which can have important applications in medical imaging and other fields.

\section{Conclusions}
\label{sec:conclusions}

The problem of reconstructing a high-quality 3D surface and achieving super-resolution is considered from a limited collection of low-resolution 2D slices.
We proposed an Euler-Elastica-based formulation in the phase-field framework, which allows for improved construction quality by capturing both local edge features and global surface smoothness.
Two numerical algorithms are developed for the numerical implementations.
Besides visual comparisons with existing methods, we have compared construction qualities by measuring Gaussian curvatures and mean curvatures, showing that the proposed model outperforms previous works.
The presented findings validate the effectiveness of addressing the challenges and offer promising prospects for various applications in medical imaging, computer vision, and other fields where high-quality surface reconstruction is essential.







\section*{Acknowledgments}
All authors would like to express their sincerest gratitude for the time and effort all referees and editors have dedicated to improving the quality of our work.
The first author is grateful for partial support from the UoL-NTHU Dual PhD Programme and would like to express particular appreciation to Dr Elie Bretin for the clear explanations and implementations of their work in publications, as well as to Dr Liam Burrows for sharing his segmented real data.
S.-H. Yang expresses thanks for the support from the Ministry of Science and Technology, Taiwan (MOST 110-2636-E-007-017).

%

\begin{spacing}{1}
\begin{table}[htbp]

\appendix
\section{Summary of notations}
\label{sec:appendix}

\renewcommand{\arraystretch}{1.3}
\begin{center}
{\small \begin{tabular}{c l}
  \toprule
  \textbf{Notations} & \textbf{Implications} \\
  \hline
  $D$
  &
  Dimension of the space.
  \\

  \hline
  $\R^{D}$
  &
  $D$-dimensional Euclidean space.
  \\

  \hline
  $\tilde{E}$ ($E_{0}$)
  &
  Target (Initial) set.
  \\

  \hline
  $E^{*}$ ($E$)
  &
  Final (Potential) result.
  \\

  \hline
  $\{\Pi_{i = 1, \dots, s}\}$
  &
  Set of the given $s$ parallel cross-sections/slices/hyperplanes $\Pi_{i}$.
  \\

  \hline
  $\omega^{in}$ ($\omega^{ex}$)
  &
  Set of the interior (exterior) restriction $\omega_{i}^{in}$ ($\omega_{i}^{ex}$) for all slices.
  \\

  \hline
  $\Omega^{in}$ ($\Omega^{ex}$)
  &
  \tabincell{l}{
  Set of the fattened interior (exterior) restriction \\ $\bigcup\limits_{i = 1}^{s} \Omega_{i, h}^{in}$ ($\bigcup\limits_{i = 1}^{s} \Omega_{i, h}^{ex}$) for all slices with the thickness $h$.
  }
  \\

  \hline
  $\mathpzc{d}_{i}(\xi, \pi_{i})$
  &
  \tabincell{l}{
  Signed distance function (oriented distance function)\\
  that identifies the distance between the given point $\xi$\\
  to an arbitrary subset $\pi_{i}$ of hyperplanes $\Pi_{i}$.
  }

  \\

    \hline
    \tabincell{c}{
    $q(\cdot)$
    }
    &
    \tabincell{l}{
    Profile function.
    }

    \\

    \hline
    \tabincell{c}{
    $W(\cdot)$
    }
    &
    \tabincell{l}{
    Double-well potential function.
    }

    \\

  \hline
    \tabincell{c}{
    $u_{\varepsilon}(\cdot)$
    }
    &
    \tabincell{l}{
    Phase-field function.
    }

    \\

    \hline
    \tabincell{c}{
    $u^{in}_{E_{0}} := u_{\varepsilon}^{in}$ \\
    ($u^{ex}_{E_{0}} := u_{\varepsilon}^{ex}$)
    }
    &
    \tabincell{l}{
    Phase-field profiles/approximations of the interior (exterior)\\
    restriction $\Omega^{in}$ ($\Omega^{ex}$) with the thickness control $\varepsilon > 0$ for $h = \varepsilon^{\alpha}$.
    }

    \\

    \hline
    \tabincell{c}{
    $\mathbbm{1}_{\Omega} (\chi_{\Omega})$
    }
    &
    \tabincell{l}{
    Indicator function (characteristic function) of set $\Omega$.
    }

    \\

    \hline
    \tabincell{c}{
    $\mathcal{H}^{D}$
    }
    &
    \tabincell{l}{
    $D$-dimensional Hausdorff measure.
    }

    \\

    \hline
    \tabincell{c}{
    $H$
    }
    &
    \tabincell{l}{
    Mean curvatures.
    }

    \\

  \hline
  $\mathcal{E} = \mathscr{P}, \mathscr{W}, \mathscr{E}$
  &
  \tabincell{l}{
  Variational energy, which can be the Perimeter-based $\mathscr{P}$, \\
  the Willmore-based $\mathscr{W}$, or the Euler-Elastica-based $\mathscr{E}$ energy.
  }

  \\

    \hline
    \tabincell{c}{
    $\mathscr{P}_{\varepsilon}, \mathscr{W}_{\varepsilon}, \mathscr{E}_{\varepsilon}$
    }
    &
    \tabincell{l}{
    Perimeter-based formulation, Willmore-based formulation, Euler-\\
    Elastica-based formulation.
    }

    \\

    \hline
    \tabincell{c}{
    $\Gamma - \lim$
    }
    &
    \tabincell{l}{
    $\Gamma$-convergence.
    }

    \\

    \hline
    \tabincell{c}{
    $\tau$
    }
    &
    \tabincell{l}{
    Synthetic time step.
    }

    \\

    \hline
    \tabincell{c}{
    $\rho$
    }
    &
    \tabincell{l}{
    Penalty parameter.
    }

    \\

    \hline
    \tabincell{c}{
    $\grad u$
    }
    &
    \tabincell{l}{
    Gradient operator:
    $\grad u = \left( \frac{\partial u}{\partial \xi_{1}}, \ldots, \frac{\partial u}{\partial \xi_{D}} \right)$.
    }

    \\

    \hline
    \tabincell{c}{
    $\operatorname{div} u$
    }
    &
    \tabincell{l}{
    Divergence operator: $\operatorname{div} u = \sum\limits_{i = 1}^{D}\frac{\partial u}{\partial \xi_{i}}$.
    }

    \\

    \hline
    \tabincell{c}{
    $\laplace u$
    }
    &
    \tabincell{l}{
    Laplacian operator: $\laplace u = \sum\limits_{i = 1}^{D}\frac{\partial^{2} u}{\partial \xi_{i}^{2}}$.
    }

    \\

    \hline
    \tabincell{c}{
    $\sigma_{\mbox{\tiny{GC}}}:=\sigma(\kappa_{G})$ \\
    $(\sigma_{\mbox{\tiny{MC}}}:=\sigma(\bar{\kappa}))$
    }
    &
    \tabincell{l}{
    Standard deviation of Gaussian (mean) curvatures.
    }

    \\

  \bottomrule
 \end{tabular}}
\end{center}
\end{table}
\end{spacing}

\medskip
Received August 2022; 1st revision May 2023; 2nd revision August 2023; early access September 2023.
\medskip

\end{document}

%% file: 00_Fig/0_Tikz/20220716_ViewRing.tex
\begin{tikzpicture}
[scale = 1.5, every node/.style={scale=1}]
\fill[blue] (0.0627,0.0187) circle[radius=1.5pt];

\fill[fill={rgb:red,255;green,0;blue,0}] (0.5,0.4) circle[radius=1.5pt];
\fill[fill={rgb:red,255;green,0;blue,0}] (0.4,-0.3) circle[radius=1.5pt];
\fill[fill={rgb:red,255;green,0;blue,0}] (-0.1,0.5) circle[radius=1.5pt];
\fill[fill={rgb:red,255;green,0;blue,0}] (-0.6,0.2) circle[radius=1.5pt];
\fill[fill={rgb:red,255;green,0;blue,0}]  (-0.3,-0.3)circle[radius=1.5pt];

\draw (0.5,0.4) -- (0.0627,0.0187);
\draw (0.4,-0.3) -- (0.0627,0.0187);
\draw (0.4,-0.3) -- (0.5,0.4);
\draw (-0.1,0.5) -- (0.0627,0.0187);
\draw  (-0.1,0.5) -- (0.5,0.4) ;
\draw (-0.6,0.2) -- (0.0627,0.0187);
\draw (-0.6,0.2)  -- (-0.1,0.5) ;
\draw (-0.3,-0.3) -- (0.0627,0.0187);
\draw (-0.3,-0.3) -- (-0.6,0.2);
\draw (-0.3,-0.3) -- (0.4,-0.3);

\fill[fill={rgb:red,119;green,255;blue,48}] (-1,0.7)  circle[radius=1.5pt];
\fill[fill={rgb:red,119;green,255;blue,48}] (-0.5,1)  circle[radius=1.5pt];
\fill[fill={rgb:red,119;green,255;blue,48}]  (0.2,1) circle[radius=1.5pt];
\fill[fill={rgb:red,119;green,255;blue,48}] (0.8,0.8) circle[radius=1.5pt];
\fill[fill={rgb:red,119;green,255;blue,48}] (1,0) circle[radius=1.5pt];
\fill[fill={rgb:red,119;green,255;blue,48}] (0.7,-0.6) circle[radius=1.5pt];
\fill[fill={rgb:red,119;green,255;blue,48}] (-0.1,-0.8) circle[radius=1.5pt];
\fill[fill={rgb:red,119;green,255;blue,48}] (-0.8,-0.4) circle[radius=1.5pt];
\fill[fill={rgb:red,119;green,255;blue,48}] (-1.2,0.3) circle[radius=1.5pt];

\draw (-1,0.7)  -- (-0.5,1)  ;
\draw (-1,0.7)  -- (-0.6,0.2) ;
\draw (-0.6,0.2)  -- (-0.5,1)  ;
\draw(-0.1,0.5)-- (-0.5,1)  ;
\draw (-1,0.7)  -- (-1.2,0.3) ;
\draw (-0.6,0.2)  -- (-1.2,0.3) ;
\draw (-0.6,0.2)  --(-0.8,-0.4) ;
\draw (-1.2,0.3) --(-0.8,-0.4) ;
\draw (-0.8,-0.4)  -- (-0.3,-0.3) ;
 \draw  (-0.8,-0.4) -- (-0.1,-0.8) ;
 \draw (-0.1,-0.8)  -- (-0.3,-0.3) ;
  \draw (-0.1,-0.8)  --(0.4,-0.3) ;
   \draw  (0.4,-0.3) -- (0.7,-0.6) ;
 \draw (-0.1,-0.8) -- (0.7,-0.6) ;
  \draw (1,0) -- (0.7,-0.6) ;
    \draw (1,0) -- (0.4,-0.3);
\draw (1,0) --(0.8,0.8);
\draw (0.5,0.4) --(0.8,0.8);
\draw (0.2,1) --(0.8,0.8);
\draw (0.2,1) --(0.5,0.4);
\draw (0.2,1) --(-0.1,0.5);
\draw (0.2,1) --(-0.5,1);
\draw (0.5,0.4) --(1,0);

\coordinate (A) at (0.0627,0.0187) {};
\coordinate (B) at (-0.6,0.2);
\coordinate (C) at (-0.1,0.5);
\draw (0.0627,0.0187) -- (-0.6,0.2) -- (-0.1,0.5)  pic [draw, angle radius=2mm] {angle};

\coordinate (A) at (0.0627,0.0187) {};
\coordinate (B) at (0.5,0.4);
\coordinate (C) at (-0.1,0.5);
\draw (A) -- (B) -- (C)  pic [draw, angle radius=2mm] {angle = C--B--A};

\node at (-0.3193,0.2036) {\scriptsize $\alpha_{ij}$};

\node at (0.2408,0.3316) {\scriptsize $\beta_{ij}$};

\node at (0.2181,-0.0018) {\scriptsize $\mathbf{v}_{i}$};

\node at (-0.1264,0.6492) {\scriptsize $\mathbf{v}_{j}$};

\coordinate (A) at (0.0627,0.0187) {};
\coordinate (B) at (-0.3,-0.3);
\coordinate (C) at (-0.6,0.2);
\draw (B) -- (A) -- (C)  pic [draw, angle radius=2mm] {angle = C--A--B};

\node at (-0.1927,-0.0294) {\scriptsize $\theta_{i_{k}}$};

\end{tikzpicture}